\documentclass[a4paper,11pt]{amsart}
\usepackage[english]{babel}
\usepackage[dvips]{graphicx}
\usepackage{amsmath,amsfonts,amsthm,amssymb,bbm,comment}
\usepackage[usenames, dvipsnames]{color}
\newcommand{\ind}{\mathbbm{1}}

\usepackage[top=1in, bottom=1in, left=0.75in, right=0.75in]{geometry}
\usepackage{hyperref}

\newtheorem{theo}{Theorem}[section]

\newtheorem{ass}{Assumption}[section]
\newtheorem{lem}{Lemma}[section]

\newtheorem{rem}{Remark}[section]
\newtheorem{example}{Example}[section]

\allowdisplaybreaks

\title[Weighted empirical means of interacting RSPs]{\bf Interacting
  Reinforced Stochastic Processes:\\ Statistical Inference based on
  the\\ Weighted Empirical Means}

\author[G. Aletti]{Giacomo Aletti}
\address{ADAMSS Center,
  Universit\`a degli Studi di Milano, Milan, Italy}
\email{giacomo.aletti@unimi.it}

\author[I. Crimaldi]{Irene Crimaldi}
\address{IMT School for Advanced Studies, Lucca, Italy}
\email{irene.crimaldi@imtlucca.it (Corresponding author)}

\author[A. Ghiglietti]{Andrea Ghiglietti}
\address{LivaNova, Milan, Italy}
\email{andrea.ghiglietti@livanova.com}

\date{\today}

\begin{document}

\maketitle

\begin{abstract}
This work deals with a system of {\em interacting reinforced
  stochastic processes}, where each process $X^j=(X_{n,j})_n$ is
located at a vertex $j$ of a finite weighted direct graph, and it can
be interpreted as the sequence of ``actions'' adopted by an agent $j$
of the network.  The interaction among the dynamics of these processes
depends on the weighted adjacency matrix $W$ associated to the
underlying graph: indeed, the probability that an agent $j$ chooses a
certain action depends on its personal ``inclination'' $Z_{n,j}$ and
on the inclinations $Z_{n,h}$, with $h\neq j$, of the other agents
according to the entries of $W$.  The best known example of reinforced
stochastic process is the P\`{o}lya urn.\\ \indent The present paper
characterizes the asymptotic behavior of the {\em weighted empirical
  means} $N_{n,j}=\sum_{k=1}^n q_{n,k} X_{k,j}$, proving their almost
sure synchronization and some central limit theorems in the sense of
stable convergence. By means of a more sophisticated decomposition of
the considered processes adopted here, these findings complete and
improve some asymptotic results for the personal inclinations
$Z^j=(Z_{n,j})_n$ and for the empirical means
$\overline{X}^j=(\sum_{k=1}^n X_{k,j}/n)_n$ given in recent papers
(e.g.~\cite{ale-cri-ghi-MEAN, ale-cri-ghi, cri-dai-lou-min}). Our work
is motivated by the aim to understand how the different rates of
convergence of the involved stochastic processes combine and, from an
applicative point of view, by the construction of confidence intervals
for the common limit inclination of the agents and of a test
statistics to make inference on the matrix $W$, based on the weighted
empirical means.  In particular, we answer a research question posed
in \cite{ale-cri-ghi-MEAN}.
\end{abstract}

\\

\smallskip
\noindent {\em 2010 AMS classification:} 60F05, 60F15, 60K35;
62P35, 91D30.

\section{Framework, model and motivations}
The stochastic evolution of systems composed by elements which
interact among each other has always been of great interest in several
scientific fields. For example, economic and social sciences deal with
agents that take decisions under the influence of other agents. In
social life, preferences and beliefs are partly transmitted by means
of various forms of social interaction and opinions are driven by the
tendency of individuals to become more similar when they
interact. Hence, a collective phenomenon, that we call
``synchronization'', reflects the result of the interactions among
different individuals. The underlying idea is that individuals have
opinions that change according to the influence of other individuals
giving rise to a sort of collective behavior.  \\ \indent In
particular, there exists a growing interest in systems of {\em
  interacting urn models} (e.g.~\cite{ale-ghi, ben, che-luc, cir,
  cri-dai-min, dai-lou-min, fortini, ieee-paper, lima, pag-sec}) and their
variants and generalizations (e.g.~\cite{ale-cri-ghi-MEAN,
  ale-cri-ghi, cri-dai-lou-min}). Our work is placed in the stream of
this scientific literature.  Specifically, it deals with the class of
the so-called {\em interacting reinforced stochastic processes}
considered in~\cite{ale-cri-ghi-MEAN, ale-cri-ghi} with a general
network-based interaction and in~\cite{cri-dai-lou-min} with a
mean-field interaction.  Generally speaking, by reinforcement in a
stochastic dynamics we mean any mechanism for which the probability
that a given event occurs has an increasing dependence on the number
of times that the same event occurred in the past. This {\em
  ``reinforcement mechanism''}, also known as {\em ``preferential
  attachment rule''} or {\em ``Rich get richer rule''} or {\em
  ``Matthew effect''}, is a key feature governing the dynamics of many
biological, economic and social systems (see, e.g. \cite{pem}).  The
best known example of reinforced stochastic process is the standard
Eggenberger-P\`{o}lya urn \cite{egg-pol,mah}, which has been widely
studied and generalized (some recent variants can be found in
\cite{ale-ghi-ros, ale-ghi-vid, ber-cri-pra-rig-barriere, chen-kuba,
  collevecchio, cri-ipergeom, ghi-pag14, ghi-vid-ros, laru-page}).
\\ \indent A {\em Reinforced Stochastic Process} (RSP) can be defined
as a stochastic process in which, along the time-steps, an agent
performs an action chosen in the set $\{0,1\}$ in such a way that the
probability of adopting ``action 1'' at a certain time-step has an
increasing dependence on the number of times that the agent adopted
``action 1'' in the previous actions. Formally, it is a stochastic
process $X=\{X_n: n\geq 1\}$ taking values in $\{0,1\}$ and such that
\begin{equation}\label{reinforced-1}
P(X_{n+1}=1\, |\, Z_0,\, X_1,....,X_n)=Z_n\,,
\end{equation}
with 
\begin{equation}\label{reinforced-2}
Z_{n}=(1-r_{n-1})Z_{n-1}+r_{n-1}X_{n}\,,
\end{equation}
where $Z_{0}$ is a random variable with values in $[0,1]$, ${\mathcal
  F}_n:=\sigma(Z_{0})\vee \sigma(X_{k}:\, 1\leq k\leq n)$ and
$(r_n)_{n\geq 0}$ is a sequence of real numbers in $(0,1)$ such that
\begin{equation}\label{ass-r-intro}
\lim_n n^{\gamma} r_n=c>0\qquad\hbox{with } 1/2<\gamma\leq 1.
\end{equation}
(We refer to \cite{cri-dai-lou-min} for a discussion on the case
$0<\gamma \leq 1/2$, for which there is a different asymptotic
behavior of the model that is out of the scope of this research work.)
The process $X$ describes the sequence of actions along the time-steps
and, if at time-step $n$, the ``action 1'' has taken place, that is
$X_n=1$, then for ``action 1'' the probability of occurrence at
time-step $(n+1)$ increases.  Therefore, the larger $Z_{n-1}$, the
higher the probability of having $X_n=1$ and so the higher the
probability of having $Z_n$ greater than $Z_{n-1}$. This means the
larger the number of times in which $X_k=1$ with $1\leq k\leq n$, the
higher the probability $Z_n$ of observing $X_{n+1}=1$.  \\ \indent As
told before, the best known example of reinforced stochastic process
is the standard Eggenberger-P\`{o}lya urn, where an urn contains $a$
red and $b$ white balls and, at each discrete time, a ball is drawn
out from the urn and then it is put again inside the urn together with
one additional ball of the same color. In this case, we have
\begin{equation*}
Z_n=\frac{a+\sum_{k=1}^{n}X_k}{a+b+n}.
\end{equation*}
It is immediate to verify that
$$
Z_{0}=\frac{a}{a+b}\quad\mbox{and}\quad Z_{n+1}=(1-r_n)Z_{n}+r_nX_{n+1}
$$ with $r_n=(a+b+n+1)^{-1}$ and so $\gamma=c=1$.  \\ \indent In the
present work we are interested in the analysis of a system of $N\geq
2$ interacting reinforced stochastic processes $\{X^j=(X_{n,j})_{n\geq
  1}:\, 1 \leq j\leq N\}$ positioned at the vertices of a weighted
directed graph $G=(V,\, E,\, W)$, where $V:=\{1,...,N\}$ denotes the
  set of vertices, $E\underline{\subset} V\times V$ the set of edges
  and $W=[w_{h,j}]_{h,j\in V\times V}$ the weighted adjacency matrix
  with $w_{h,j}\geq 0$ for each pair of vertices.  The presence of the
  edge $(h,j)\in E$ indicates a ``direct influence'' that the vertex
  $h$ has on the vertex $j$ and it corresponds to a strictly positive
  element $w_{h,j}$ of $W$, that represents a weight quantifying this
  influence. We assume the weights to be normalized so that
  $\sum_{h=1}^N w_{h,j}=1$ for each $j\in V$. The interaction between
  the processes $\{X^j:\,j\in V\}$ is explicitly inserted in Equation
  \eqref{reinforced-1} and it is modeled as follows: for any $n\geq
  0$, the random variables $\{X_{n+1,j}:\,j\in V\}$ are conditionally
  independent given ${\mathcal F}_{n}$ with
\begin{equation}\label{interacting-1-intro}
P(X_{n+1,j}=1\, |\, {\mathcal F}_n)=\sum_{h=1}^N w_{h,j} Z_{n,h}
=w_{jj}Z_{n,j}+\sum_{h\neq j} w_{h,j} Z_{n,h}\,,
\end{equation}
where ${\mathcal F}_n:=\sigma(Z_{0,h}: h\in V)\vee \sigma(X_{k,j}:\,
1\leq k\leq n,\,j\in V )$ and, for each $h\in V$, the evolution
dynamics of the single process $(Z_{n,h})_{n\geq 0}$ is the same as in
\eqref{reinforced-2}, that is
\begin{equation}\label{interacting-2-intro}
Z_{n,h}=(1-r_{n-1})Z_{n-1,h}+r_{n-1}X_{n,h}\,,
\end{equation}
with $Z_{0,h}$ a random variable taking values in $[0,1]$ and
$(r_n)_{n\geq 0}$ a sequence of real numbers in $(0,1)$ such that
condition \eqref{ass-r-intro} holds true.  \\ \indent As an example,
we can imagine that $G=(V, E)$ represents a network of $N$ individuals
that at each time-step have to make a choice between two possible
actions $\{0,1\}$. For any $n\geq 1$, the random variables
$\{X_{n,j}:\,j\in V\}$ take values in $\{0,1\}$ and they describe the
actions adopted by the agents of the network along the time-steps;
while each random variable $Z_{n,h}$ takes values in $[0,1]$ and it
can be interpreted as the ``personal inclination'' of the agent $h$ of
adopting ``action 1''. Thus, the probability that the agent $j$ adopts
``action 1'' at time-step $(n+1)$ is given by a convex combination of
$j$'s own inclination and the inclination of the other agents at
time-step $n$, according to the ``influence-weights'' $w_{h,j}$ as in
\eqref{interacting-1-intro}. Note that, from a mathematical point of
view, we can have $w_{jj}\neq 0$ or $w_{jj}=0$. In both cases we have
a reinforcement mechanism for the personal inclinations of the agents:
indeed, by \eqref{interacting-2-intro}, whenever $X_{n,h}=1$, we have
a positive increment in the personal inclination of the agent $h$,
that is $Z_{n,h}\geq Z_{n-1,h}$. However, only in the case $w_{jj}>0$,
this fact results in a greater probability of having $X_{n+1,j}=1$
according to \eqref{interacting-1-intro}. Therefore, if $w_{jj}>0$,
then we have a ``true self-reinforcing'' mechanism; while, in the
opposite case, we have a reinforcement property only in the own
inclination of the single agent, but this does not affect the
probability \eqref{interacting-1-intro}.  \\ \indent The literature
\cite{ale-cri-ghi, cri-dai-lou-min, cri-dai-min, dai-lou-min} focus on
the asymptotic behavior of the stochastic processes of the personal
inclinations $\{Z^j=(Z_{n,j})_n:\, j\in V \}$ of the agents; while
\cite{ale-cri-ghi-MEAN} studies the average of times in which the
agents adopt ``action 1'', i.e. the stochastic processes of the
empirical means $\{\overline{X}_n^j=(\frac{1}{n}\sum_{k=1}^{n}
X_{k,j})_n:\, j\in V\}$.  The results given in
\cite{ale-cri-ghi-MEAN}, together with the resulting statistical
tools, represent a great improvement in any area of application, since
the ``actions'' $X_{n,j}$ adopted by the agents of the network are
much more likely to be observed than their personal inclinations
$Z_{n,j}$ of adopting these actions. More specifically, in that paper,
under suitable assumptions, it is proved that all the empirical means
converge almost surely to the same limit random variable (almost sure
synchronization), which is also the common limit random variable of
the stochastic processes $Z^j=(Z_{n,j})_n$, say $Z_{\infty}$.
Moreover, some Central Limit Theorems (CLTs) for the empirical means
hold true and they lead to the construction of asymptotic confidence
intervals for the common limit random variable $Z_{\infty}$ and of a
statistical test to make inference on the weighted adjacency matrix
$W$ of the network in the case $\gamma=1$.  \\ \indent In the present
paper, we continue in this direction: indeed, we not only extend the
results obtained in \cite{ale-cri-ghi-MEAN} for the empirical means to
the ``weighted empirical means'', but, using a more sophisticated
decomposition, we obtain two improvements: first, we here handle the
two cases, $\gamma<1$ and $\gamma=1$, in the same way (while in
\cite{ale-cri-ghi-MEAN} we use two different arguments) and, second,
we here solve a research question posed in \cite{ale-cri-ghi-MEAN}
and, consequently, we succeed in constructing a test statistics to make
inference on the weighted adjacency matrix $W$ of the network for all
values of the model parameters (not only in the case $\gamma=1$). More
precisely, in this paper we focus on the weighted average of times in
which the agents adopt ``action 1'', i.e. we study the stochastic
processes of the {\em weighted empirical means}
$\{N^j=(N_{n,j})_{n}:\, j\in V\}$ defined, for each $j\in V$, as
$N^j_0:=0$ and, for any $n\geq 1$,
\begin{equation}\label{medie-empiriche-intro}
N_{n,j}:=\sum_{k=1}^{n}q_{n,k} X_{k,j}\,,\quad\hbox{
where }\quad q_{n,k}:=\frac{a_k}{\sum_{l=1}^na_l},
\end{equation}
with $(a_k)_{k\geq 1}$ a suitable sequence of strictly positive
real numbers.  Since $\sum_{k=1}^nq_{n,k}=1$, we have the relation
$$\sum_{k=1}^{n-1}q_{n,k}X_{k,j}\ =\
\frac{\sum_{l=1}^{n-1}a_l}{\sum_{l=1}^na_l}
\left(\sum_{k=1}^{n-1}q_{n-1,k}X_{k,j}\right)\ =\
(1-q_{n,n})N_{n-1,j}
$$
and so we get
\begin{equation}\label{interacting-N}
N_{n,j}=\left(1-q_{n,n}\right)N_{n-1,j}+q_{n,n}X_{n,j}.
\end{equation}
Note that this framework includes as special case the process of
the standard empirical means studied in~\cite{ale-cri-ghi-MEAN},
which corresponds to the case $a_k=1$ for any $k\geq1$ (and hence
$q_{n,k}=1/n$ for any $1\leq k\leq n$).  Furthermore, the above
dynamics \eqref{interacting-1-intro}, \eqref{interacting-2-intro}
and \eqref{interacting-N} can be expressed in a compact form,
using the random vectors
$\mathbf{X}_{n}:=(X_{n,1},\dots,X_{n,N})^{\top}$ for $n\geq 1$,
$\mathbf{N}_{n}:=(N_{n,1},\dots,N_{n,N})^{\top}$ and
$\mathbf{Z}_{n}:=(Z_{n,1},\dots,Z_{n,N})^{\top}$ for $n\geq 0$,
as 
\begin{equation}\label{eq:dynamic-0}
E[\mathbf{X}_{n+1}|\mathcal{F}_{n}]=W^{\top}\,\mathbf{Z}_{n}\,,
\end{equation}
where $W^{\top}\mathbf{1}=\mathbf{1}$ by the normalization of the
weights, and
\begin{equation}\label{eq:dynamic}
\left\{\begin{aligned} &\mathbf{Z}_{n}\ =\
\left(1-r_{n-1}\right)\mathbf{Z}_{n-1}\ +\ r_{n-1}\mathbf{X}_{n},\\
&\mathbf{N}_{n}\ =\ \left(1-q_{n,n}\right)\mathbf{N}_{n-1}\ +\
q_{n,n}\mathbf{X}_{n}.
\end{aligned}\right.
\end{equation}
Under suitable assumptions, we prove the almost sure synchronization
of the stochastic processes $N^j=(N_{n,j})_n$, with $j\in V$, toward
the same limit random variable $Z_\infty$, which is the common limit
random variable of the stochastic processes $Z^j=(Z_{n,j})_n$ and we
provide some CLTs in the sense of stable convergence.  In particular,
we assume
\begin{equation}\label{ipotesi-q-generale}
\lim_n n^{\nu }q_{n,n}=q>0\qquad\mbox{with } 1/2<\nu\leq 1
\end{equation}
and the asymptotic covariances in the provided CLTs depend on the
random variable $Z_{\infty}$, on the eigen-structure of the weighted
adjacency matrix $W$ and on the parameters $\gamma,\, c$ and $\nu,\,q$
governing the asymptotic behavior of the sequence $(r_n)_n$ and
$(q_{n,n})_n$, respectively. We also discuss the possible statistical
applications of these convergence results: asymptotic confidence
intervals for the common limit random variable $Z_{\infty}$ and test
statistics to make inference on the weighted adjacency matrix $W$ of
the network. In particular, as said before, we obtain a statistical
test on the matrix $W$ for all values of the model parameters (not
only in the case $\gamma=\nu=q=1$ as in
\cite{ale-cri-ghi-MEAN}). Moreover, our results give a hint regarding
a possible ``optimal choice'' of $\nu$ and $q$ and so point out the
advantages of employing the weighted empirical means with $\nu<1$,
instead of the simple empirical means.  \\ \indent Finally, we point
out that the existence of joint central limit theorems for the pair
$(\mathbf{Z}_n,\mathbf{N}_n)$ is not obvious because the ``discount
factors'' in the dynamics of the increments
$(\mathbf{Z}_{n}-\mathbf{Z}_{n-1})_n$ and
$(\mathbf{N}_{n}-\mathbf{N}_{n-1})_n$ are generally different. Indeed,
as shown in~\eqref{eq:dynamic}, these two stochastic processes follow
the dynamics
\begin{equation}\label{eq-increments-intro}
\left\{
\begin{aligned}
&\mathbf{Z}_{n}-\mathbf{Z}_{n-1}\ =\
r_{n-1}\left(\mathbf{X}_{n}-\mathbf{Z}_{n-1}\right),\\
&\mathbf{N}_{n}-\mathbf{N}_{n-1}\ =\
q_{n,n}\left(\mathbf{X}_{n}-\mathbf{N}_{n-1}\right),
\end{aligned}
\right.
\end{equation}
and so, when we assume $\nu \neq \gamma$, it could be surprising that
in some cases there exists a common convergence rate for the pair
$(\mathbf{Z}_n,\mathbf{N}_n)$.  It is worthwhile to note that dynamics
similar to~\eqref{eq-increments-intro} have already been considered in
the Stochastic Approximation literature.  Specifically,
in~\cite{mok-pel} the authors established a CLT for a pair of
recursive procedures having two different step-sizes. However, this
result does not apply to our situation.  Indeed, the covariance
matrices $\Sigma_\mu$ and $\Sigma_\theta$ in their main result
(Theorem~1) are deterministic, while the asymptotic covariance
matrices in our CLTs are random (as said before, they depend on the
random variable $Z_\infty$).  This is why we do not use the simple
convergence in distribution, but we employ the notion of stable
convergence, which is, among other things, essential for the
considered statistical applications. Moreover in~\cite{mok-pel}, the
authors find two different convergence rates, depending on the two
different step-sizes, while, as already said, we find a common
convergence rate also in some cases with $\nu\neq \gamma$.\\ \indent
Summing up, this work complete the convergence results obtained in
\cite{ale-cri-ghi-MEAN, ale-cri-ghi} for the stochastic processes of
the personal inclinations $Z^j=(Z_{n,j})_n$ and of the empirical means
$\overline{X}^j=(\overline{X}_{n,j})_n$, and it extend them to the
weighted empirical means $N^j=(N_{n,j})_n$. However the main focus
here concerns the new decomposition employed for the analysis of the
asymptotic behavior of the pair $(\mathbf{Z}_n,\mathbf{N}_n)$, that,
among other things, allows us to solve the research question arisen in
\cite{ale-cri-ghi-MEAN} regarding the statistical test on $W$ in the
case $\gamma< 1$. Thus, in what follows, we will go fast on the point
in common with \cite{ale-cri-ghi-MEAN, ale-cri-ghi}, while we
concentrate on the novelties.  \\ \indent The rest of the paper is
organized as follows. In Section \ref{section_model} we describe the
notation and the assumptions used along the paper.  In Section
\ref{section_main_results} and Section \ref{section_statistics} we
illustrate our main results and we discuss some possible statistical
applications. An interesting example of interacting system is also
provided in order to clarify the statement of the theorems and the
related comments. Section \ref{section_proofs_1} and Section
\ref{section_proofs_statistics} contain the proofs or the main steps
of the proofs of our results, while the technical details have been
gathered in the appendix.  In particular, Subsection
\ref{section_decomposition} contains the main ingredient of the proofs
of the CLTs, that is a suitable decomposition of the joint stochastic
process $(\mathbf{Z}_n,\mathbf{N}_n)$. Finally, for the reader's
convenience, the appendix also supplies a brief review on the notion
of stable convergence and its variants (e.g. see~\cite{crimaldi-2009,
  crimaldi-libro, cri-let-pra-2007, hall-1980, z}).

\section{Notation and assumptions}\label{section_model}
Throughout all the paper, we will assume $N\geq 2$ and adopt the same
notation used in~\cite{ale-cri-ghi-MEAN,ale-cri-ghi}.  In particular,
we denote by ${\mathcal Re}(z)$, ${\mathcal Im}(z)$, $\overline{z}$
and $|z|$ the real part, the imaginary part, the conjugate and the
modulus of a complex number $z$.  Then, for a matrix $A$ with complex
elements, we let $\overline{A}$ and $A^{\top}$ be its conjugate and
its transpose, while we indicate by $|A|$ the sum of the modulus of
its elements. The identity matrix is denoted by $I$, independently of
its dimension that will be clear from the context.  The spectrum of
$A$, i.e.~the set of all the eigenvalues of $A$ repeated with their
multiplicity, is denoted by $Sp(A)$, while its sub-set containing the
eigenvalues with maximum real part is denoted by $\lambda_{\max}(A)$,
i.e. $\lambda^*\in \lambda_{\max}(A)$ whenever ${\mathcal
  Re}(\lambda^*)=\max\{ {\mathcal Re}(\lambda):\, \lambda\in Sp(A)
\}$.  The notation $\hbox{diag}(a_1,\dots,a_d)$ indicates the diagonal
matrix of dimension $d$ with diagonal elements $a_1,\dots,
a_d$. Finally, we consider any vector $\mathbf{v}$ as a matrix with
only one column (so that all the above notations apply to
$\mathbf{v}$) and we indicate by $\|\mathbf{v}\|$ its norm, i.e. $
\|\mathbf{v} \|^2 = \overline{\mathbf{v}}^{\top}\mathbf{v}$.  The
vectors and the matrices whose elements are all ones or zeros are
denoted by $\mathbf{1}$ and $\mathbf{0}$, respectively, independently
of their dimension that will be clear from the context.\\

For the matrix $W$ we make the following assumption:

\begin{ass}\label{ass:W}
The weighted adjacency matrix $W$ is irreducible and
diagonalizable.
\end{ass}

The irreducibility of $W$ reflects a situation in which all the
vertices are connected among each others and hence there are no
sub-systems with independent dynamics (see~\cite{ale-cri-ghi, ale-ghi}
for further details). The diagonalizability of $W$ allows us to find a
non-singular matrix $\widetilde{U}$ such that
$\widetilde{U}^{\top}W(\widetilde{U}^{\top})^{-1}$ is diagonal with
complex elements $\lambda_j\in Sp(W)$.  Notice that each column
$\mathbf{u}_j$ of $\widetilde{U}$ is a left eigenvector of $W$
associated to the eigenvalue $\lambda_j$.  Without loss of
generality, we take $\|\mathbf{u}_j\|=1$.  Moreover, when the
multiplicity of some $\lambda_j$ is bigger than one, we set the
corresponding eigenvectors to be orthogonal.  Then, if we define
$\widetilde{V}=(\widetilde{U}^{\top})^{-1}$, we have that each column
$\mathbf{v}_j$ of $\widetilde{V}$ is a right eigenvector of $W$
associated to the eigenvalue $\lambda_j$ such that
\begin{equation}\label{eq:relazioni-0}
\mathbf{u}_j^{\top}\,\mathbf{v}_j=1,\quad\mbox{ and }\qquad
\mathbf{u}_h^{\top}\,\mathbf{v}_j=0,\ \forall\, h\neq j.
\end{equation}

These constraints combined with the above assumptions on $W$
(precisely, $w_{h,j}\geq 0$, $W^{\top}\mathbf{1}=\mathbf{1}$ and
the irreducibility) imply, by Frobenius-Perron Theorem, that
$\lambda_1:=1$ is an eigenvalue of $W$ with multiplicity one,
$\lambda_{\max}(W)=\{1\}$ and
\begin{equation}\label{eq:relazioni-1}
\mathbf{u}_1=N^{-1/2}\mathbf{1}, \qquad N^{-1/2}{\mathbf
1}^{\top}{\mathbf v}_1=1\qquad \mbox{and}\qquad
v_{1,j}:=[\mathbf{v}_1]_j>0\;\ \forall 1\leq j\leq N.
\end{equation}
Moreover, we recall the relation
\begin{equation}\label{eq-NUM1}
\sum_{j=1}^N\mathbf{u}_j\mathbf{v}_j^{\top}=I.
\end{equation}
Finally, we set $\alpha_j:=1-\lambda_j\in{\mathbb C}$ for each $j\geq
2$, i.e.~for each $\lambda_j$ belonging to $Sp(W)\setminus\{1\}$, and
we denote by $\lambda^*$ an eigenvalue belonging to
$Sp(W)\setminus\{1\}$ such that
\begin{equation*}
{\mathcal Re}(\lambda^*)=\max\left\{ {\mathcal Re}(\lambda_j):\,
\lambda_j\in Sp(W)\setminus\{1\}\right\}.
\end{equation*}

\indent Throughout all the paper, we assume that the two
sequences $(r_n)_{n\geq 0}$ and $(q_{n,n})_{n\geq 1}$, which appear in
\eqref{eq:dynamic}, satisfy the following assumption:
\begin{ass}\label{ass:r_n}
There exist real constants $\gamma,\nu\in (1/2,1]$ and $c,q>0$ such
that
\begin{equation}\label{ass:condition_r_n_gamma}
r_{n-1}=\frac{c}{n^{\gamma}}+O\left(\frac{1}{n^{2\gamma}}\right)\qquad\hbox{
and }\qquad
q_{n,n}=\frac{q}{n^{\nu}}+O\left(\frac{1}{n^{2\nu}}\right).
\end{equation}
\end{ass}
In particular, it follows
\begin{equation*}
\lim_n n^\gamma r_n=c>0\qquad\mbox{and}\qquad 
\lim_n n^{\nu} q_{n,n}=q>0.
\end{equation*}
The following remark will be useful for a certain proof in the sequel.
\begin{rem}
\rm Recalling that $q_{n,n}=a_n/\sum_{l=1}^n a_l$, the second relation
in \eqref{ass:condition_r_n_gamma} implies that
$\sum_{n=1}^{+\infty}a_n=+\infty$. Indeed, the above relation together
with $\sum_{n=1}^{+\infty}a_n=\ell<+\infty$ entails $a_n=q\ell
n^{-\nu}+O(n^{-2\nu})$ and so, since $\nu\leq 1$,
$\sum_{n=1}^{+\infty}a_n=+\infty$, which is a contradiction.
\end{rem}
In the special case considered in~\cite{ale-cri-ghi-MEAN}, where the
random variables $N_{n,j}$ correspond to the standard empirical means
($a_n=1$ for each $n$), we have $\nu=1$ and $q=1$.  Other possible
choices are the following:
\begin{itemize}
\item $\sum_{l=1}^na_l=n^\delta$ with $\delta>0$, which brings to
\begin{equation*}
a_n=n^{\delta}-(n-1)^{\delta}
\end{equation*}
and  
\begin{equation*}
q_{n,n}=1-\frac{\sum_{l=1}^{n-1}a_l}{\sum_{l=1}^na_l}=
1-\left(1-\frac{1}{n}\right)^{\delta}= \delta n^{-1}+O(n^{-2}),
\end{equation*}
so that we have  $\nu=1$ and $q=\delta>0$;
\item $\sum_{l=1}^na_l=\exp(b n^\delta)$ with $b>0$ and
  $\delta\in(0,1/2)$, which brings to
\begin{equation*}
a_n=\exp(b n^\delta)-\exp(b (n-1)^\delta)
\end{equation*}
and
\begin{equation*}
\begin{split}
q_{n,n}&=1-\frac{\sum_{l=1}^{n-1}a_l}{\sum_{l=1}^na_l}=
1-\exp\left[b\left((n-1)^{\delta}-n^{\delta}\right)\right]\\
&=
b
n^{\delta}\left( 1- (1-n^{-1})^\delta \right)
+O\left(\, n^{2\delta}( 1-(1-n^{-1})^\delta )^2 \,\right)
=b n^{\delta}\left( \delta n^{-1}+O(n^{-2}) \right)+O(n^{-(2-2\delta)})\\
&=b \delta n^{-(1-\delta)}+O(n^{-(2-\delta)})+O(n^{-(2-2\delta)})=
b \delta n^{-(1-\delta)}+O(n^{-2(1-\delta)}),
\end{split}
\end{equation*}
so that $\nu=(1-\delta)\in (1/2,1)$ and  $q=b\delta >0$.
\end{itemize}

To ease the notation, we set $\widehat{r}_{n-1}:=cn^{-\gamma}$ and
$\widehat{q}_{n,n}:=qn^{-\nu}$, so that condition
\eqref{ass:condition_r_n_gamma} can be rewritten as
\begin{equation*}
r_{n-1}=\widehat{r}_{n-1}+O\left(\frac{1}{n^{2\gamma}}\right)\qquad\hbox{
and }\qquad
q_{n,n}=\widehat{q}_{n,n}+O\left(\frac{1}{n^{2\nu}}\right).
\end{equation*}

For the CLTs provided in the sequel, we make also the following assumption:
\begin{ass}\label{ass:case1}
When $\gamma=1$, we assume the condition
$c>1/[2(1-\mathcal{R}e(\lambda^*))]$, i.e~${\mathcal
  R}e(\lambda^*)<1-(2c)^{-1}$. When $\nu=1$, we assume $q>1/2$.
\end{ass}

Note that in Assumption \ref{ass:r_n} condition
\eqref{ass:condition_r_n_gamma} for the sequence $(r_n)_n$ is slightly
more restrictive than the one assumed in \cite{ale-cri-ghi-MEAN,
  ale-cri-ghi}. However, it is always verified in the applicative
contexts we have in mind. The reason behind this choice is that we
want to avoid some technical complications in order to focus on the
differences brought by the use of the weighted empirical means,
specially on the relationship between the pair $(\gamma,\nu)$ and the
asymptotic behaviors of the considered stochastic processes. For the
same reason, in the CLTs for the case $\nu=\gamma$, we add also the
following assumption:
\begin{equation}\label{ass-q}
q\neq c\alpha_j\; \forall j\geq 2\,.
\end{equation}
We think that this condition is not necessary. Indeed, if there exists
$j\geq 2$ such that $q=c\alpha_j$, we conjecture that our proofs still
work (but changing the asymptotic expression adopted for a certain
quantity, see the proof of Lemma \ref{conv-general-sum}) and they lead
to exactly the same asymptotic covariances provided in the CLTs under
the above condition \eqref{ass-q}. Our conjecture is motivated by the
fact that this is what happens in \cite{ale-cri-ghi-MEAN} for the
simple empirical means. Moreover, the expressions obtained for the
asymptotic covariances in the following CLTs do not require condition
\eqref{ass-q}. However, as told before, we do not want to make the
following proofs even heavier and so, when $\nu=\gamma$, we will work
under condition \eqref{ass-q}.

\section{Main results on the joint stochastic
process}\label{section_main_results}
The first achievement concerns the {\em almost sure
  synchronization} of all the involved stochastic processes, that is
\begin{equation}\label{as-synchro}
\mathbf{Y}_{n}:=
\begin{pmatrix}
\mathbf{Z}_{n}\\
\mathbf{N}_{n}
\end{pmatrix}
\stackrel{a.s.}{\longrightarrow} Z_{\infty}\mathbf{1},
\end{equation}
where $Z_{\infty}$ is a random variable with values in $[0,1]$. This
fact means that all the stochastic processes $Z^j=(Z_{n,j})_n$
and $N^j=(N_{n,j})_n$ positioned at different vertices $j\in V$
of the graph converge almost surely to the same random variable
$Z_\infty$.  \\
\indent The synchronization for the first component of ${\mathbf
  Y}_n$, that is
\begin{equation}\label{synchro-first-component} 
[{\mathbf Y}_n]_1=\mathbf{Z}_{n}\stackrel{a.s.}{\longrightarrow}
Z_{\infty}\mathbf{1}\,,
\end{equation}
is the result contained in \cite[Theorem 3.1]{ale-cri-ghi}, while for
the second component, we prove in the present work the following result:
\begin{theo}\label{th-as-synchro}
Under Assumptions \ref{ass:W} and \ref{ass:r_n}, we have 
\begin{equation}\label{synchro-second-component}
[\mathbf{Y}_{n}]_2=\mathbf{N}_{n}
\stackrel{a.s.}{\longrightarrow} Z_{\infty}\mathbf{1}\,.
\end{equation}
\end{theo}
Regarding the distribution of~$Z_{\infty}$, we recall that
\cite[Theorems 3.5 and 3.6]{ale-cri-ghi} state the following two
properties:
\begin{itemize}
\item[(i)] $P(Z_\infty=z)=0$ for any $z\in (0,1)$.
\item[(ii)] If we have
  $P(\bigcap_{j=1}^N\{Z_{0,j}=0\})+P(\bigcap_{j=1}^N\{Z_{0,j}=1\})<1$,
  then $P(0<Z_\infty<1)>0$.
\end{itemize}
In particular, these facts entail that the asymptotic covariances in
the following CLTs are ``truly'' random. Indeed, their random part
$Z_\infty(1-Z_{\infty})$ is different from zero with probability
greater than zero and almost surely different from a constant in
$(0,1)$. \\ \indent Furthermore, it is interesting to note that the
almost sure synchronization holds true without any assumptions on the
initial configuration $\mathbf{Z}_{0}$ and for any choice of the
weighted adjacency matrix $W$ with the required assumptions. Finally,
note that the synchronization is induced along time independently of
the fixed size $N$ of the network, and so it does not require a
large-scale limit (i.e. the limit for $N\to +\infty$), which is usual
in statistical mechanics for the study of interacting particle
systems.  \\ \indent Regarding the convergence rate and the
second-order asymptotic distribution of
$(\mathbf{Y}_n-Z_{\infty}\mathbf{1})$, setting for each $\gamma\in
(1/2,1]$
\begin{equation}\label{def:gamma_0}
\gamma_0\ :=\ \max\left\{\frac{1}{2},2\gamma-1\right\}
\in [1/2,\, 1]\,,
\end{equation}
\begin{equation}\label{eq-Sigma-tilde}
\widetilde{\Sigma}_\gamma:=
\widetilde{\sigma}_{\gamma}^2\mathbf{1}\mathbf{1}^{\top}
\qquad\hbox{with}\quad \widetilde{\sigma}_{\gamma}^2 :=
\frac{\|\mathbf{v}_1\|^2c^2}{N(2\gamma-1)}
\end{equation}
and
\begin{equation}\label{def:Utilde}
\widetilde{U}=
\begin{pmatrix}
\mathbf{u}_{1} & \mathbf{u}_{2} & ... & \mathbf{u}_{N}
\end{pmatrix}=
\begin{pmatrix}
N^{-1/2}\mathbf{1} & U
\end{pmatrix}
\quad\mbox{with } 
U:=
\begin{pmatrix}
\mathbf{u}_{2} & ... & \mathbf{u}_{N}
\end{pmatrix}
\,,
\end{equation}
we obtain the following result:

\begin{theo}\label{clt-Y}
Under all the assumptions stated in Section \ref{section_model}, the
following statements hold true:
\begin{itemize}

\item[(a)] If $1/2<\nu<\gamma_0$, then
\begin{equation}
n^{\nu/2}(\mathbf{Y}_n-Z_{\infty}\mathbf{1}) {\longrightarrow}\
\mathcal{N} \left(\ \mathbf{0}\ ,\
Z_{\infty}(1-Z_{\infty})\begin{pmatrix}
\mathbf{0} & \mathbf{0}\\
\mathbf{0} & \widetilde{U}S^{(q)}\widetilde{U}^{\top}
\end{pmatrix} \
\right)\qquad\hbox{stably},
\end{equation}
where, for $1\leq j_1,j_2\leq N$,
\begin{equation}\label{def:S_q}
[S^{(q)}]_{j_1 j_2}:=
\frac{q}{2}\mathbf{v}_{j_1}^{\top}\mathbf{v}_{j_2}\,.
\end{equation}

\item[(b)] If $\gamma_0<\nu<1$, then
\begin{equation}
n^{\gamma-\frac{1}{2}}(\mathbf{Y}_n-Z_{\infty}\mathbf{1})
{\longrightarrow}\ \mathcal{N} \left(\ \mathbf{0}\ ,\
Z_{\infty}(1-Z_{\infty}) \widetilde{\Sigma}_{\gamma}\
\right)\quad\hbox{stably}.
\end{equation}

\item[(c)] If $\nu=\gamma_0<1$, then
\begin{equation}
n^{\gamma-\frac{1}{2}}(\mathbf{Y}_n-Z_{\infty}\mathbf{1})
{\longrightarrow}\mathcal{N}\left(\mathbf{0}\ ,\
Z_{\infty}(1-Z_{\infty}) \left(\widetilde{\Sigma}_{\gamma}\ +\
\begin{pmatrix}
\mathbf{0} & \mathbf{0}\\
\mathbf{0} & \widetilde{U}S^{(q)}\widetilde{U}^{\top}
\end{pmatrix}\right)\right)\;\hbox{stably},
\end{equation}
where $S^{(q)}$ is the same matrix defined in (a) by \eqref{def:S_q}.

\item[(d)] If $\nu=\gamma_0=1$ (that is $\nu=\gamma=1$), then
\begin{equation}
\sqrt{n}(\mathbf{Y}_n-Z_{\infty}\mathbf{1})
{\longrightarrow}
\mathcal{N} \left(\mathbf{0}\ ,\ Z_{\infty}(1-Z_{\infty}) \left(
\widetilde{\Sigma}_{1}+
\begin{pmatrix}
\widetilde{U}S^{11}\widetilde{U}^{\top} & \widetilde{U}S^{12}\widetilde{U}^{\top} \\
\widetilde{U}S^{21}\widetilde{U}^{\top} & \widetilde{U}S^{22}\widetilde{U}^{\top}
\end{pmatrix}
\right)\right)\;\hbox{stably},
\end{equation}
where $S^{21}=(S^{12})^{\top}$ and, for $2\leq j_1,j_2,j\leq N$,
\begin{equation*}
\begin{split}
&[S^{11}]_{11} =\ [S^{11}]_{j_11}\ =\ [S^{11}]_{1j_2}\ :=\ 0,\\
&[S^{11}]_{j_1 j_2} :=\ \frac{c^2}{c(\alpha_{j_1}+\alpha_{j_2})-1}
\mathbf{v}_{j_1}^{\top}\mathbf{v}_{j_2},\\
&[S^{12}]_{11}\ =\ [S^{12}]_{1j_2}\ :=\ 0,\\
&[S^{12}]_{j_11}\ :=\ 
\frac{c(q-c)}{c\alpha_{j_1}+q-1}\mathbf{v}_{j_1}^{\top}\mathbf{v}_{1},\\
&[S^{12}]_{j_1 j_2} :=\
\frac{cq(c\alpha_{j_1}+c-1)}{(c\alpha_{j_1}+c\alpha_{j_2}-1)(c\alpha_{j_1}+q-1)}
\mathbf{v}_{j_1}^{\top}\mathbf{v}_{j_2},\\
&[S^{22}]_{11} :=\ \frac{(q-c)^2}{2q-1}\|\mathbf{v}_{1}\|^2,\\
&[S^{22}]_{j1}=[S^{22}]_{1j} :=\ 
\frac{q(q-c)(c+q-1)}{(c\alpha_{j}+q-1)(2q-1)}
\mathbf{v}_{j}^{\top}\mathbf{v}_{1},\\
&[S^{22}]_{j_1 j_2} :=\ 
q^2\frac{c^3(\alpha_{j_1}+\alpha_{j_2})+2c^2q(\alpha_{j_1}\alpha_{j_2}+1)
-c^2(\alpha_{j_1}\alpha_{j_2}+\alpha_{j_1}+\alpha_{j_2}+2)}
{(2q-1)(c(\alpha_{j_1}+\alpha_{j_2})-1)(c\alpha_{j_1}+q-1)(c\alpha_{j_2}+q-1)}
\mathbf{v}_{j_1}^{\top}\mathbf{v}_{j_2}\\
& +\ q^2\frac{c(q-1)^2(\alpha_{j_1}+\alpha_{j_2})-(2c+q-1)(q-1)}
{(2q-1)(c(\alpha_{j_1}+\alpha_{j_2})-1)(c\alpha_{j_1}+q-1)(c\alpha_{j_2}+q-1)}
\mathbf{v}_{j_1}^{\top}\mathbf{v}_{j_2}.
\end{split}
\end{equation*}

\item[(e)] If $\gamma_0<\nu=1$, then
\begin{equation}\label{risultato-e}
n^{\gamma-\frac{1}{2}}(\mathbf{Y}_n-Z_{\infty}\mathbf{1})
{\longrightarrow}\ \mathcal{N} \left(\ \mathbf{0}\ ,\
Z_{\infty}(1-Z_{\infty}) \left(\widetilde{\Sigma}_{\gamma} +
\frac{\|\mathbf{v}_{1}\|^2 c^2}{N[2q-(2\gamma-1)]}
\begin{pmatrix}
\mathbf{0} & \mathbf{0}\\
\mathbf{0} & \mathbf{1}\mathbf{1}^{\top}\\
\end{pmatrix} \right) \
\right) \quad\hbox{stably}.
\end{equation}

\end{itemize}
\end{theo}

\begin{rem}
\rm Looking at the asymptotic covariance matrices in the different
cases of the above theorem, note that in case (a) the convergence rate
of the first component is bigger then the one of the second
component. Indeed, from our previous work \cite{ale-cri-ghi}, we know
that it is $n^{\gamma_0/2}$. On the other hand, there are cases (see
(b), (c) and (e)) in which the convergence rates of the two components
are the same, although the discount factors $r_n\sim c n^{-\gamma}$
and $q_{n,n}\sim q n^{-\nu}$ in \eqref{eq:dynamic} have different
convergence rates.
\end{rem}

\begin{rem}
\rm Recall that we have
$$
1 \leq 1 + \|\mathbf{v}_1 -\mathbf{u}_1\|^2=
\|\mathbf{v}_1\|^2\leq N.
$$ 
Therefore we obtain the following lower and upper bounds (that do
not depend on $W$) for $\widetilde{\sigma}^2_\gamma$ and for the
second term in the asymptotic covariance of relation
\eqref{risultato-e}:
$$
\frac{c^2}{N(2\gamma-1)}\leq
\widetilde{\sigma}^2_\gamma\leq 
\frac{c^2}{2\gamma-1}
$$
and
$$
\frac{c^2}{N[2q-(2\gamma-1)]}\leq 
\frac{\|\mathbf{v}_{1}\|^2 c^2}{N[2q-(2\gamma-1)]}\leq
\frac{c^2}{2q-(2\gamma-1)}.
$$ 
Notice that the lower bound is achieved when $\mathbf{v}_1 =\mathbf{u}_1= 
N^{-1/2} \mathbf{1}$, i.e. when $W$ is doubly stochastic, 
which means $W\mathbf{1}= W^{\top}\mathbf{1} = \mathbf{1}$.
\end{rem}

\begin{rem}\label{rem:compare_old_paper}
\rm The results of Theorem~\ref{clt-Y} extend those presented
in~\cite{ale-cri-ghi-MEAN}, since they are valid only for $q=\nu=1$
which here corresponds to a special situation in case (d) and (e) of
Theorem~\ref{clt-Y}. Indeed, when $q=\nu=1$ and $\gamma<1$ we have
that~\cite[Theorem 3.2]{ale-cri-ghi-MEAN} coincides with the result of
case (e) while, when $q=\nu=1$ and $\gamma=1$, we have that~\cite[Theorem
3.4]{ale-cri-ghi-MEAN} coincides with the result of case (d),
because we have 
\begin{equation*}
\begin{pmatrix}
\widetilde{U}S^{11}\widetilde{U}^{\top} & \widetilde{U}S^{12}\widetilde{U}^{\top} \\
\widetilde{U}S^{21}\widetilde{U}^{\top} & \widetilde{U}S^{22}\widetilde{U}^{\top}
\end{pmatrix}\ =\
\begin{pmatrix}
U\widehat{S}_{ZZ}U^{\top} & U\widehat{S}_{ZN}\widetilde{U}^{\top} \\
\widetilde{U}\widehat{S}_{ZN}^{\top}U^{\top} & 
\widetilde{U}\widehat{S}_{NN}\widetilde{U}^{\top}
\end{pmatrix},
\end{equation*}
where the matrices $\widehat{S}_{ZZ},\, \widehat{S}_{ZN}$ and $\widehat{S}_{NN}$
are defined in \cite{ale-cri-ghi-MEAN}.
\end{rem}

\begin{rem}\label{rem-N=1}
\rm The main goal of this work is to provide results for a system of
$N\geq 2$ interacting reinforced stochastic processes. However, it is
worth to note that Theorem \ref{th-as-synchro}, and the consequent
limit \eqref{as-synchro}, hold true also for $N=1$. Moreover,
statements (d) and (e) of Theorem \ref{clt-Y} with $N=1$ are true and
they correspond to \cite[Theorems 3.2 and
  3.3]{ale-cri-ghi-MEAN}. Finally statements (a), (b) and (c) of
Theorem \ref{clt-Y} with $N=1$ (and so without the condition on
$\lambda^*$) can be proven with the same proof provided in the sequel
(see the following Remark \ref{rem-N=1-bis}).
\end{rem}

We conclude this section with the example of the ``mean-field''
interaction.

\begin{example}\label{esempio}
\rm 
The mean-field interaction can be expressed in terms of a particular
weighted adjacency matrix $W$ as follows: for any $1 \leq j_1, j_2 \leq N$
\begin{equation}\label{medione} 
w_{j_1,j_2} = \frac{\alpha}{N} + (1-\alpha)\delta_{j_1,j_2} \qquad
\mbox{with }\alpha\in [0, 1],
\end{equation} 
where $\delta_{j_1,j_2}$ is equal to $1$ when $j_1 = j_2$ and to $0$
otherwise. Note that $W$ in \eqref{medione} is irreducible for $\alpha
> 0$ and so we are going to consider this case. Since $W$ is doubly
stochastic, we have $\mathbf{v}_1 = \mathbf{u}_1 =
N^{-1/2}\mathbf{1}$. Moreover, since $W$ is also symmetric, we have
$\widetilde{U} = \widetilde{V}$ and so $\widetilde{U}
\widetilde{U}^{\top}=I$ and $\widetilde{V}^{\top} \widetilde{V} =
I$. Finally, we have $\lambda_j =1-\alpha$ for all $j\geq 2$ and,
consequently, we obtain
\begin{equation*}
\begin{split}
&S^{(q)}=\frac{q}{2}I, \quad
\{[S^{11}]_{j_1j_2}:\, 2\leq j_1,j_2\leq N\}=\frac{c^2}{2c\alpha -1}I,
\\
&[S^{12}]_{j_11}=0\; \mbox{for } 2\leq j_{1}\leq N,
\quad
\{[S^{12}]_{j_1j_2}:\, 2\leq j_1,j_2\leq N\}=
\frac{qc(c\alpha+c-1)}{(2c\alpha -1)(c\alpha+q-1)}I,\\
&[S^{22}]_{11}=\frac{(q-c)^2}{2q-1},\quad 
[S^{22}]_{j1}=[S^{22}_{1j}]=0\;\mbox{for } 2\leq j\leq N,\\
&\{[S^{22}]_{j_1j_2}:\, 2\leq j_1,j_2\leq N\}=\\
&\frac{
(qc)^2[(\alpha^2+1)(2q-1)+2\alpha (c-1)-1+(2\alpha -c^{-1})(q-1)^2-2c^{-1}(q-1)]}
{(2q-1)(2c\alpha-1)(c\alpha+q-1)^2}I,
\end{split}
\end{equation*}
and the condition ${\mathcal R}e(\lambda^*)<1-(2c)^{-1}$ when $\gamma
= 1$ becomes $2c\alpha >1$.
\end{example}

\section{Useful results for statistical applications}
\label{section_statistics}
The first convergence result provided in this section can be used for
the construction of asymptotic confidence intervals for the limit
random variable $Z_\infty$, that requires to know the following quantities: 
\begin{itemize}
\item $N$: the number of agents in the network; 
\item $\mathbf{v}_1$: the right eigenvector of $W$ associated to
  $\lambda_1=1$ (note that it is not required to know the whole
  weighted adjacency matrix $W$, e.g. we have
  $\mathbf{v}_1=\mathbf{u}_1= N^{-1/2}\mathbf{1}$ for any doubly
  stochastic matrix); 
\item $\gamma$ and $c$: the parameters that describe the
  first-order asymptotic approximation of the sequence $(r_n)_n$;
\item $\nu$ and $q$: the parameters that describe the first-order
  asymptotic approximation of the sequence $(q_{n,n})_n$ (recall that
  the weights $q_{n,k}$ are chosen and so $\nu$ and $q$ are always
  known and, moreover, they can be optimally chosen).
\end{itemize}
We point out that it is not required the observation of the random
variables $Z_{n,j}$, nor the knowledge of the initial random variables
$\{Z_{0,j}:\, j\in V\}$ and nor of the exact expression of the
sequence $(r_n)_n$. They are based on the weighted empirical means of
the random variables $X_{n,j}$, that are typically observable.
\\ \indent The second result stated in this section can be employed
for the construction of asymptotic critical regions for statistical
tests on the weighted adjacency matrix $W$ based on the weighted
empirical means (given the values of $\gamma,\,\nu,\,c$ and $q$). In
particular, we point out that in our previous work
\cite{ale-cri-ghi-MEAN} we succeeded to provide a testing procedure
based on the standard empirical means only for the case $\gamma=1$;
while we announced further future investigation for the case
$1/2<\gamma<1$. In the present work we face and solve this issue,
providing a test statistics for all the values of the
parameters. Indeed the following Theorem \ref{thm:N_apice} covers all 
the cases for the pair $(\gamma, \nu)$.  \\

\indent Let us consider the decomposition
$\mathbf{N}_{n}=\mathbf{1}\widetilde{N}_n+{\mathbf{N}}'_{n}$,
where
\begin{equation}\label{eq-NUM4}
\mathbf{1}\widetilde{N}_n:=\mathbf{u}_1\mathbf{v}_1^{\top}\mathbf{N}_n
=
N^{-1/2}\mathbf{1}\mathbf{v}_{1}^{\top}\mathbf{N}_{n}
\qquad\mbox{ and }\qquad
{\mathbf{N}}'_{n}:=\mathbf{N}_n-\mathbf{1}\widetilde{N}_n=
(I-\mathbf{u}_{1}\mathbf{v}_{1}^{\top})\mathbf{N}_{n}.
\end{equation}

Concerning the first term, by \eqref{eq:relazioni-1} and the almost
sure synchronization \eqref{as-synchro}, we immediately obtain
\begin{equation*}
\widetilde{N}_n \stackrel{a.s.}{\longrightarrow} Z_{\infty}.
\end{equation*}
Moreover, under all the assumptions stated in Section
\ref{section_model}, setting
\begin{equation}\label{eq:def-sigma-tilde}
\widetilde{\sigma}^2:=\frac{\|\mathbf{v}_1\|^2}{N}\times
\left\{
\begin{aligned}
& \frac{q}{2}&\mbox{ if } \nu<\gamma_0\quad\mbox{or}\quad \nu=\gamma_0<1,\\
& \frac{(q-c)^2}{2q-1}&\mbox{ if }  \nu=\gamma_0=1\;  
(\mbox{that is } \nu=\gamma=1),\\
& \frac{c^2}{2q-(2\gamma-1)}&\mbox{ if }  \gamma_0<\nu=1,
\end{aligned}
\right.
\end{equation}
where $\gamma_0$ and $\widetilde{\sigma}_{\gamma}^2$ are defined in
\eqref{def:gamma_0} and in~\eqref{eq-Sigma-tilde}, respectively, we
have the following result:

\begin{theo}\label{thm:N_tilde}
Under all the assumptions stated in Section \ref{section_model}, the
following statements hold true:
\begin{itemize}
\item[(a)] If $\nu<\gamma_0$, then
\begin{equation*}
n^{\nu/2}(\widetilde{N}_n-Z_{\infty})\ \longrightarrow\
\mathcal{N} \left(\ 0\ ,\
Z_{\infty}(1-Z_{\infty})\widetilde{\sigma}^2 \
\right)\quad\hbox{stably}.
\end{equation*}

\item[(b)] If $\gamma_0<\nu<1$, then
\begin{equation*}
n^{\gamma-\frac{1}{2}}(\widetilde{N}_n-Z_{\infty})\
\longrightarrow\ \mathcal{N} \left(\ 0\ ,\
Z_{\infty}(1-Z_{\infty})\widetilde{\sigma}_{\gamma}^2 \
\right)\quad\hbox{stably}.
\end{equation*}

\item[(c)] If $\nu=\gamma_0$ or $\nu=1$ (i.e. $\nu=\gamma_0<1$ or
  $\nu=\gamma_0=1$ or $\gamma_0<\nu=1$), then
\begin{equation*}
n^{\gamma-\frac{1}{2}}(\widetilde{N}_n-Z_{\infty})\
{\longrightarrow}\ \mathcal{N} \left(\ 0\ ,\
Z_{\infty}(1-Z_{\infty})\left(\widetilde{\sigma}_{\gamma}^2
+\widetilde{\sigma}^2\right) \ \right)\quad\hbox{stably}.
\end{equation*}

\end{itemize}
\end{theo}
Note that $\widetilde{\sigma}^2$ has not been defined in the case
$\gamma_0<\nu<1$, i.e. in the case (b) of the above result, because in
this case it does not appear in the asymptotic covariance
matrix.\\

\indent In the following remark, we point out the advantages of
employing the weighted empirical means with $\nu<1$, instead of the
simple empirical means (for which we have $\nu=q=1$), providing a
brief discussion on the possible ``optimal choice'' of $\nu$ and $q$:

\begin{rem}\label{remark-ottimo-1}
\rm The convergence rates and the asymptotic variances expressed in
the cases of the above Theorem \ref{thm:N_tilde} allows us to make
some considerations on the existence of an ``optimal'' choice of the
parameters $\nu$ and $q$ in order to ``maximize the convergence'' of
$\widetilde{N}_n$ towards the random limit $Z_\infty$.  Indeed, first
note that the convergence rate in case (a) is slower than the rates of
the other two cases, and, moreover, the asymptotic variance in case
(c) is strictly larger than the variance in case (b). Hence, the
interval $\gamma_0<\nu<1$ in case (b) provides an ``optimal'' range
of values where the parameter $\nu$ should be chosen. In addition,
looking into the proof of Theorem \ref{thm:N_tilde}, it is possible
to investigate more deeply into the behavior of $\widetilde{N}_n$ and
so derive more accurate optimality conditions on the values of $\nu$
and $q$ (see the following Remark \ref{remark-ottimo-2}).
\end{rem}

\indent Analogously, concerning the term
${\mathbf{N}}'_{n}=(I-\mathbf{u}_{1}\mathbf{v}_{1}^{\top})\mathbf{N}_{n}$,
from \eqref{eq:relazioni-1} and the almost sure synchronization
\eqref{as-synchro}, we obtain
\begin{equation*}
{\mathbf{N}}'_{n} \stackrel{a.s.}{\rightarrow} \mathbf{0}.
\end{equation*}

Moreover, setting
\begin{equation}\label{def:U_minus_1}
\widetilde{U}_{-1}:=
\begin{pmatrix}
\mathbf{0} & \mathbf{u}_{2} & ... & \mathbf{u}_{N}
\end{pmatrix}
=
\begin{pmatrix}
\mathbf{0} & U
\end{pmatrix}
\,,
\end{equation}
we get the following theorem:

\begin{theo}\label{thm:N_apice}
Under all the assumptions stated in Section \ref{section_model}, the
following statements hold true:
\begin{itemize}
\item[(a)] If $\nu<\gamma$, then
\begin{equation*}
n^{\frac{\nu}{2}}{\mathbf{N}}'_{n} {\longrightarrow}\ \mathcal{N}
\left(\ \mathbf{0}\ ,\
Z_{\infty}(1-Z_{\infty})\widetilde{U}_{-1}S^{(q)}\widetilde{U}_{-1}^{\top}
\ \right)\qquad\hbox{stably},
\end{equation*}
where $S^{(q)}$ is defined in \eqref{def:S_q}.

\item[(b)] If $\nu=\gamma$, then
\begin{equation*}
n^{\frac{\nu}{2}}{\mathbf{N}}'_{n} {\longrightarrow}\
\mathcal{N} 
\left(\ \mathbf{0}\ ,\ 
Z_{\infty}(1-Z_{\infty})\widetilde{U}_{-1}S_\gamma^{22}\widetilde{U}_{-1}^{\top}
\ \right)\quad\hbox{stably},
\end{equation*}
where, for any $2\leq j_1,j_2\leq N$, we have that
$[S^{22}_{\gamma}]_{11}$, $[S^{22}_{\gamma}]_{1 j_2}$ and
$[S^{22}_{\gamma}]_{j_1 1}$ are not needed to be defined since the
first column of $\widetilde{U}_{-1}$ is $\mathbf{0}$, while the 
remaining elements $[S^{22}_{\gamma}]_{j_1j_2}$ are defined as
\begin{equation}\label{eq-S_gamma_22}
\begin{aligned}
&q^2\frac{c^3(\alpha_{j_1}+\alpha_{j_2})+
2c^2q(\alpha_{j_1}\alpha_{j_2}+1)-
\ind_{\{\gamma=1\}}c^2(\alpha_{j_1}\alpha_{j_2}+\alpha_{j_1}+\alpha_{j_2}+2)}
{(2q-\ind_{\{\gamma=1\}})(c(\alpha_{j_1}+\alpha_{j_2})-
\ind_{\{\gamma=1\}})(c\alpha_{j_1}+q-\ind_{\{\gamma=1\}})
(c\alpha_{j_2}+q-\ind_{\{\gamma=1\}})}
\mathbf{v}_{j_1}^{\top}\mathbf{v}_{j_2}\\
& +\ q^2
\frac{c(q-\ind_{\{\gamma=1\}})^2(\alpha_{j_1}+\alpha_{j_2})-
\ind_{\{\gamma=1\}}(2c+q-1)(q-1)}
{(2q-\ind_{\{\gamma=1\}})(c(\alpha_{j_1}+\alpha_{j_2})-
\ind_{\{\gamma=1\}})
(c\alpha_{j_1}+q-\ind_{\{\gamma=1\}})(c\alpha_{j_2}+q-\ind_{\{\gamma=1\}})}
\mathbf{v}_{j_1}^{\top}\mathbf{v}_{j_2}\,.
\end{aligned}
\end{equation}

\item[(c)] If $\gamma<\nu$, then
\begin{equation*}
n^{\frac{\nu}{2}}{\mathbf{N}}'_{n} {\longrightarrow}\ \mathcal{N}
\left(\ \mathbf{0}\ ,\
Z_{\infty}(1-Z_{\infty})\widetilde{U}_{-1}S\widetilde{U}_{-1}^{\top}
\ \right)\qquad\hbox{stably},
\end{equation*}
where, for any $2\leq j_1,j_2\leq N$, we have that $[S]_{11}$, $[S]_{1
  j_2}$ and $[S]_{j_1 1}$ are not needed to be defined since the first
column of $\widetilde{U}_{-1}$ is $\mathbf{0}$, while the remaining
elements $[S]_{j_1j_2}$ are defined as
\begin{equation*}
q^2\left(\!
\left(\frac{\lambda_{j_1}\lambda_{j_2}}{\alpha_{j_1}\alpha_{j_2}}\right)
\frac{1}{2q-\ind_{\{\nu=1\}}(2\gamma-1)}
+\left(\frac{\lambda_{j_1}}{\alpha_{j_1}}+
\frac{\lambda_{j_2}}{\alpha_{j_2}}
\right)\frac{1}{2q-\ind_{\{\nu=1\}}\gamma}
+\frac{1}{2q-\ind_{\{\nu=1\}}} \!\right)
\mathbf{v}_{j_1}^{\top}\mathbf{v}_{j_2}.
\end{equation*}

\end{itemize}
\end{theo}
Note that the convergence rate for $(\mathbf{N}'_n)$ is always
$n^{\nu/2}$.
\\

\indent In the following example we go on with the analysis of the mean-field
interaction.

\begin{example} 
\rm If we consider again the mean-field interaction (see
\eqref{medione}), we have
$\mathbf{N}_n'=(I-N^{-1}\mathbf{1}\mathbf{1}^{\top})\mathbf{N}_n$
(because $\mathbf{v}_1=\mathbf{u}_1=N^{-1/2}\mathbf{1}$).  Moreover,
since $\widetilde{U}=\widetilde{V}$ and so ${\widetilde
  V}^{\top}\widetilde{V}=I$, we find $S^{(q)}=\frac{q}{2} I$,
\begin{equation*}
\begin{split}
&\{[S^{22}_{\gamma}]_{j_1j_2}:\, 2\leq j_1,j_2\leq N\}=
s^{22}_\gamma I\qquad\mbox{with}
\\
&s^{22}_\gamma\!:=\!\frac{q^2\! [
c^2(\alpha^2+1)(2q-\ind_{\{\gamma=1\}})\!+\!
2c^2\alpha( c-\ind_{\{\gamma=1\}})\!-\!
\ind_{\{\gamma=1\}} c^2\!+\!
2\alpha c(q-\ind_{\{\gamma=1\}})^2\!-\!
\ind_{\{\gamma=1\}}(2c+q-1)(q-1)]}
{(2q-\ind_{\{\gamma=1\}})(2c\alpha-\ind_{\{\gamma=1\}})
(c\alpha+q-\ind_{\{\gamma=1\}})^2}
\end{split}
\end{equation*}
and
$$
S=sI\quad \mbox{with } s:=q^2\left(\!
\left(\frac{1-\alpha}{\alpha}\right)^2
\frac{1}{2q-\ind_{\{\nu=1\}}(2\gamma-1)}
+
2\frac{(1-\alpha)}{\alpha}\frac{1}{2q-\ind_{\{\nu=1\}}\gamma}
+\frac{1}{2q-\ind_{\{\nu=1\}}} \!\right).
$$ 
Hence, since ${\widetilde U}_{-1}I{\widetilde
  U}_{-1}^{\top}=UU^{\top}=I-N^{-1}\mathbf{1}\mathbf{1}^{\top}$, we
get that
$$
n^{\nu/2}(I-N^{-1}\mathbf{1}\mathbf{1}^{\top})\mathbf{N}_n
\longrightarrow 
\mathcal{N}\left(
0, Z_\infty(1-Z_\infty)s^*(I-N^{-1}\mathbf{1}\mathbf{1}^{\top})
\right)\quad\mbox{stably},
$$ 
where $s^*$ is equal to $q/2$ or $s^{22}_\gamma$ or $s$, according to
the values of $\nu$ and $\gamma$. Finally, using the relations
$U^{\top}U=I$ and $UU^{\top}= I-N^{-1}\mathbf{1}\mathbf{1}^{\top}$ and
employing $\widetilde{N}_n$ as a strong consistent estimator of
$Z_\infty$, we get
$$
\frac{n^{\nu/2}}{ \left[\widetilde{N}_n(1-\widetilde{N}_n)s^*\right]^{1/2}}
U^{\top}\mathbf{N}_n\stackrel{d}\sim\mathcal{N}(0,I)
$$ 
and  
$$
\frac{n^{\nu}}{\widetilde{N}_n(1-\widetilde{N}_n)s^*}
\mathbf{N}_n^{\top}(I-N^{-1}\mathbf{1}\mathbf{1}^{\top})\mathbf{N}_n
\stackrel{d}\sim \chi^2_{N-1}.
$$ 
Given the values of $\gamma,\,\nu,\,c$ and $q$, this result can be used
in order to perform a statistical test on the parameter $\alpha$ in
the definition of $W$ (see \eqref{medione}).
\end{example}

\section{Proof of the results on the joint stochastic process}
\label{section_proofs_1}

Here we prove the convergence results stated in Section
\ref{section_main_results}.

\subsection{Proof of Theorem \ref{th-as-synchro}}
As already recalled (see \eqref{synchro-first-component}), we have
$\mathbf{Z}_n\stackrel{a.s.}\to Z_{\infty}$. Hence, since the
condition $W^\top\mathbf{1}=\mathbf{1}$ and the equality
\eqref{eq:dynamic-0}, we get
$E[\mathbf{X}_n|\mathcal{F}_{n-1}]\stackrel{a.s.}\to
Z_{\infty}\mathbf{1}$. Therefore, the convergence
$\mathbf{N}_n\stackrel{a.s.}\to Z_{\infty}\mathbf{1}$ follows from
\cite[Lemma B.1]{ale-cri-ghi-MEAN} with $c_k=k^{\nu}$,
$v_{n,k}=c_kq_{n,k}$ and $\eta=1$. Note that the assumptions on the
weights $q_{n,k}=a_k/\sum_{l=1}^na_l$, easily implies that $c_k$ and
$v_{n,k}$ satisfy the conditions required in the employed lemma: 
indeed, by definition, we have $\sum_{k=1}^n q_{n,k}=1$ and from the
second relation in \eqref{ass:condition_r_n_gamma} we get
$\sum_{n=1}^{+\infty}a_n=+\infty$ and
\begin{equation*}
n^\nu a_n=q\sum_{l=1}^n a_l+O\left(n^{-\nu}\sum_{l=1}^n a_l\right)=
q\sum_{l=1}^n a_l+O\left( a_n(n^{\nu}q_{n,n})^{-1}\right)
=q\sum_{l=1}^n a_l+O(a_n),
\end{equation*}
and so we obtain 
\begin{equation*}
\lim_n v_{n,k}=c_ka_k\lim_n \frac{1}{\sum_{l=1}^n a_l}=0,
\quad \lim_n  v_{n,n}=\lim_{n}c_nq_{n,n}=q,\quad
  \lim_n\sum_{k=1}^n\frac{v_{n,k}}{c_k}=\lim_n\sum_{k=1}^n q_{n,k}=1
 \end{equation*}
and 
\begin{equation*}
\begin{split}
\sum_{k=1}^n|v_{n,k}-v_{n,k-1}|&=
\frac{1}{\sum_{l=1}^na_l}\sum_{k=1}^n k^{\nu} a_k-(k-1)^{\nu} a_{k-1}
\\
&=\frac{1}{\sum_{l=1}^na_l}
\left[\sum_{k=1}^n q\left(\sum_{l=1}^k a_l-\sum_{l=1}^{k-1}a_l\right)
+O\left(\sum_{k=1}^n a_k\right)\right]
\\
&=q\frac{\sum_{k=1}^n a_k}{\sum_{l=1}^na_l}+O(1)=O(1).
\end{split}
\end{equation*}
\qed

\subsection{Decomposition of the joint stochastic process}
\label{section_decomposition}
In this section we describe the main tool used in the following
proofs, that is a suitable decomposition of the joint stochastic
process $\mathbf{Y}:=(\mathbf{Y}_n)_n$. Indeed, in order to determine
the convergence rate and the second-order asymptotic distribution of
$(\mathbf{Y}_n - Z_{\infty}\mathbf{1})$ for any values of the
parameters, we need to decompose $\mathbf{Y}$ into a sum of
``primitive'' stochastic processes, and then establish the asymptotic
behavior for each one of them. As we will see, they converge at
different rates. \\

\indent Let us express the dynamics~\eqref{eq:dynamic} of the
stochastic processes $(\mathbf{Z}_{n})_n$ and $(\mathbf{N}_{n})_n$ as
follows:
\begin{equation}\label{eq:dynamic_SA}
\left\{\begin{aligned} &\mathbf{Z}_{n}-\mathbf{Z}_{n-1}\ =\
-\widehat{r}_{n-1}\left(I-W^{\top}\right)\mathbf{Z}_{n-1}\
+\ \widehat{r}_{n-1}\Delta\mathbf{M}_{n}\ +\ \Delta\mathbf{R}_{Z,n},\\
&\mathbf{N}_{n}-\mathbf{N}_{n-1}\ =\
-\widehat{q}_{n,n}\left(\mathbf{N}_{n-1}-W^{\top}\mathbf{Z}_{n-1}\right)\ +\
\widehat{q}_{n,n}\Delta\mathbf{M}_{n}\ +\ \Delta\mathbf{R}_{N,n},
\end{aligned}\right.
\end{equation}
where
$\Delta\mathbf{M}_{n}:=(\mathbf{X}_{n}-W^{\top}\mathbf{Z}_{n-1})$ is a
martingale increment with respect to the filtration~${\mathcal
  F}:=({\mathcal F}_{n})_n$, while $\Delta\mathbf{R}_{Z,n}:=
(r_{n-1}-\widehat{r}_{n-1})(\mathbf{X}_n-\mathbf{Z}_{n-1})$ and
$\Delta\mathbf{R}_{N,n}:=(q_{n,n}-\widehat{q}_{n,n})(\mathbf{X}_n-\mathbf{Z}_{n-1})$
are two remainder terms.  Hence, by means of~\eqref{eq:dynamic_SA},
the dynamics of the stochastic process $\mathbf{Y}$ can be expressed as 
\begin{equation}\label{eq:dynamic_SA_Y}
\mathbf{Y}_{n}= (I-Q_n)\mathbf{Y}_{n-1}+
R_n\Delta\mathbf{M}_{Y,n}+\Delta\mathbf{R}_{Y,n},
\end{equation}
where
$\Delta\mathbf{M}_{Y,n}:=(\Delta\mathbf{M}_n,\Delta\mathbf{M}_n)^{\top}$,
$\Delta\mathbf{R}_{Y,n}:=(\Delta\mathbf{R}_{Z,n},\Delta\mathbf{R}_{N,n})^{\top}$,
\begin{equation}\label{def:Q_R}
Q_n\ :=\ \begin{pmatrix}
\widehat{r}_{n-1}(I-W^{\top}) & \mathbf{0}  \\
-\widehat{q}_{n,n}W^{\top} & \widehat{q}_{n,n}I
\end{pmatrix}
\qquad\mbox{ and }\qquad R_n\ :=\ \begin{pmatrix}
\widehat{r}_{n-1}I & \mathbf{0} \\
\mathbf{0} & \widehat{q}_{n,n}I
\end{pmatrix}.
\end{equation}
Now, we want to decompose the stochastic process $\mathbf{Y}$ in a sum
of stochastic processes, whose dynamics are of the same types of
\eqref{eq:dynamic_SA_Y}, but more tractable. To this purpose, we set
$U_{j}:=(\mathbf{u}_{j(1)},\mathbf{u}_{j(2)})$ for each $j=1,\dots,N$,
and we impose the following relations:
\begin{equation}\label{eq:UP}
U_{j}=U^{*}_jP_j \qquad\mbox{ with }\qquad
U^{*}_{j}:=\begin{pmatrix}
\mathbf{u}_j & \mathbf{0} \\
\mathbf{0} & \mathbf{u}_j
\end{pmatrix}
\qquad\mbox{ and }\qquad P_{j}:=
\begin{pmatrix}
1 & 0 \\
g(\lambda_j) & 1
\end{pmatrix},
\end{equation}
and, for any $n\geq 1$,
\begin{equation}\label{eq:QU=UD}
Q_nU_j=U_jD_{Q,j,n},\qquad\mbox{ where}\qquad
D_{Q,j,n}\ :=\ \begin{pmatrix}
\widehat{r}_{n-1}(1-\lambda_j) & 0 \\
-\lambda_jh_n(\lambda_j) & \widehat{q}_{n,n}
\end{pmatrix}\,.
\end{equation}
We recall that $\lambda_j$ and $\mathbf{u}_j$ denote the eigenvalues
and the left eigenvectors of $W$, respectively. The above functions
$g$ and $h_n$ will be suitable defined later on. In particular, we
will define $h_n$ in such a way that the sequence $(h_n(\lambda_j))_n$
converges to zero at the biggest possible rate.  In order to solve the
above system of equations, we firstly observe that, by \eqref{eq:UP}, 
we have
\begin{equation}\label{eq-uji}
\mathbf{u}_{j(1)}=\begin{pmatrix}
\mathbf{u}_j \\
g(\lambda_j)\mathbf{u}_j
\end{pmatrix}
,\qquad
\mathbf{u}_{j(2)}=\begin{pmatrix}
\mathbf{0} \\
\mathbf{u}_j
\end{pmatrix}\,,
\end{equation}
\begin{equation}\label{eq:1}
Q_nU_j=Q_nU^{*}_jP_j=U^{*}_j
\begin{pmatrix}
\widehat{r}_{n-1}(1-\lambda_j) & 0 \\
-\widehat{q}_{n,n}\lambda_j & \widehat{q}_{n,n}
\end{pmatrix}
P_j=U^{*}_j
\begin{pmatrix}
\widehat{r}_{n-1}(1-\lambda_j) & 0 \\
-\widehat{q}_{n,n}\lambda_j+\widehat{q}_{n,n}g(\lambda_j) & \widehat{q}_{n,n}
\end{pmatrix}
\end{equation}
and
\begin{equation}\label{eq:2}
U_jD_{Q,j,n}=U^{*}_jP_jD_{Q,j,n} =U^{*}_j
\begin{pmatrix}
\widehat{r}_{n-1}(1-\lambda_j) & 0 \\
\widehat{r}_{n-1}(1-\lambda_j)g(\lambda_j)-\lambda_jh_n(\lambda_j)
& \widehat{q}_{n,n}
\end{pmatrix}.
\end{equation}
Then, combining together~\eqref{eq:1} and~\eqref{eq:2} in order to
satisfy~\eqref{eq:QU=UD}, we obtain
$$-\widehat{q}_{n,n}\lambda_j+\widehat{q}_{n,n}g(\lambda_j)\ =
\widehat{r}_{n-1}(1-\lambda_j)g(\lambda_j)-\lambda_jh_n(\lambda_j),$$ from
which we get the equality 
\begin{equation}\label{eq:3}
\lambda_j[\widehat{q}_{n,n}-h_n(\lambda_j)]\ =\
g(\lambda_j)[\widehat{q}_{n,n}-\widehat{r}_{n-1}(1-\lambda_j)].
\end{equation}
Now, for all values of $\gamma$, $\nu$ and $j\in\{1,\dots,N\}$, we
want to define $g(\lambda_j)$ and $h_n(\lambda_j)$ in such a way
that~\eqref{eq:3} is verified for any $n$ and $h_n(\lambda_j)$
vanishes to zero with the biggest possible rate.  To this end, we note
that by~\eqref{eq:3} we have the following two facts:
\begin{itemize}
\item If $\lambda_j=0$, we can set $g(\lambda_j)=g(0)=0$ and
  $h_n(\lambda_j)=h_n(0)$ is not relevant.
\item If $\lambda_j\neq 0$, $g(\lambda_j)$ does not depend on $n$ only
  if $h_n(\lambda_j)=\widehat{r}_{n-1}(1-\lambda_j)$, which implies
  $g(\lambda_j)=\lambda_j$, or if $h_n(\lambda_j)=\widehat{q}_{n,n}$,
  which implies $g(\lambda_j)=0$.
\end{itemize}
Hence, since $\widehat{r}_{n-1}$ and $\widehat{q}_{n,n}$ have
convergence rates $n^{\gamma}$ and $n^{\nu}$, respectively, we choose
  to set
\begin{equation}\label{eq:def-h}
h_n(x)\ :=\ \left\{\begin{aligned}
&\widehat{r}_{n-1}(1-x)\ &\mbox{ if }\nu<\gamma,\\
&\widehat{q}_{n,n}\ind_{\{x\neq1\}}\ &\mbox{ if }\nu\geq\gamma
\end{aligned}\right.
\end{equation}
and 
\begin{equation}\label{eq:def-g}
g(x)\ :=\ \left\{\begin{aligned}
&x\ &\mbox{ if }\nu<\gamma,\\
&\ind_{\{x=1\}}\ &\mbox{ if }\nu\geq\gamma.
\end{aligned}\right.
\end{equation}
Note that, since $\lambda_1=1$, we have $g(\lambda_1)=g(1)=1$ and
$h_n(\lambda_1)=h_n(1)=0$ regardless the values of $\nu$ and $\gamma$.
\\

\indent Now, recalling that $\mathbf{v}_j$, for $j=1,\dots, N$, denote
the right eigenvectors of $W$, we set
$V_{j}:=(\mathbf{v}_{j(1)},\mathbf{v}_{j(2)})$, for each $j=1,\dots,
N$, with the condition
\begin{equation*}
V_{j}=V_j^{*}P_j^{-\top} \qquad\mbox{ where }\qquad
V^{*}_{j}:=\begin{pmatrix}
\mathbf{v}_j & \mathbf{0} \\
\mathbf{0} & \mathbf{v}_j
\end{pmatrix}
\qquad\mbox{ and }\qquad P^{-\top}_{j}:=
\begin{pmatrix}
1 & -g(\lambda_j) \\
0 & 1
\end{pmatrix},
\end{equation*}
so that we have 
\begin{equation}\label{eq-vji}
\mathbf{v}_{j(1)}=\begin{pmatrix}
\mathbf{v}_j \\
\mathbf{0}
\end{pmatrix}
\qquad\hbox{and}\qquad
\mathbf{v}_{j(2)}=\begin{pmatrix}
-g(\lambda_j)\mathbf{v}_j \\
\mathbf{v}_j
\end{pmatrix}
\,.
\end{equation}
Note that, we also have  
\begin{equation}\label{eq:VQ=DV}
V_j^{\top}Q_n=D_{Q,j,n}V_j^{\top}.
\end{equation}
Moreover, by \eqref{eq:relazioni-0}, we have
\begin{equation}\label{eq:relazioni-0-bis}
\mathbf{u}_{j(i)}^{\top}\,\mathbf{v}_{j(i)}=1,\quad\mbox{ and
}\qquad \mathbf{u}_{h(l)}^{\top}\,\mathbf{v}_{j(i)}=0,\ \forall\,
h\neq j \mbox{ or }l\neq i.
\end{equation}

Finally, since $\{\mathbf{u}_{j(i)}:\, j=1,..,N;\, i=1,2\}$ and
$\{\mathbf{v}_{j(i)}:\,j=1,\dots,N;\, i=1,2\}$ satisfy, for any
$j\in\{1,\dots,N\}$, the relation
\begin{equation}
U_{j}V_{j}^{\top}=
\mathbf{u}_{j(1)}\mathbf{v}_{j(1)}^{\top}+\mathbf{u}_{j(2)}\mathbf{v}_{j(2)}^{\top}=
\begin{pmatrix}
\mathbf{u}_j\mathbf{v}_j^{\top} &\mathbf{0} \\
\mathbf{0} &\mathbf{u}_j\mathbf{v}_j^{\top}
\end{pmatrix}
\end{equation}
and since \eqref{eq-NUM1}, the stochastic process $\{\mathbf{Y}_{n}:\,
n\geq1\}$ can be decomposed as
\begin{equation}\label{eq-NUM2-new}
\mathbf{Y}_{n}\ =\ \sum_{j=1}^N\mathbf{Y}_{j,n}
\qquad\hbox{with } 
\mathbf{Y}_{j,n}\ :=\ U_jV_j^{\top}\mathbf{Y}_{n}\,.
\end{equation}
The dynamics of each term $\mathbf{Y}_{j,n}$ can be deduced
from~\eqref{eq:dynamic_SA_Y} by multiplying this equation by
$U_{j}V_{j}^{\top}=U^{*}_jV_j^{*\top}$ and using \eqref{eq:VQ=DV} and
the relation $V_j^{\top}\mathbf{Y}_{n}=V_j^{\top}\mathbf{Y}_{j,n}$. We
thus obtain
\begin{equation}\label{eq:proof-dynamics-Y_each_j}
\mathbf{Y}_{j,n}\ =\
U_j\left(I-D_{Q,j,n}\right)V_j^{\top}\mathbf{Y}_{j,n-1}\ +
U^{*}_jD_{R,n}V_j^{*\top}\Delta\mathbf{M}_{Y,n}+U_jV_j^{\top}\Delta\mathbf{R}_{Y,n},
\end{equation}
where
\begin{equation}\label{def:matrix-D_Rn}
D_{R,n}\ :=\ \begin{pmatrix}
\widehat{r}_{n-1} & 0 \\
0 & \widehat{q}_{n,n}
\end{pmatrix}.
\end{equation}
For the sequel, it will be useful to decompose $\mathbf{Y}_n$ further 
as 
\begin{equation}\label{ulteriore-decomposizione}
\mathbf{Y}_{n}\ =\ \sum_{j=1}^N\mathbf{Y}_{j,n}\ =\
\sum_{j=1}^N\mathbf{Y}_{j(1),n}\ +\
\sum_{j=1}^N\mathbf{Y}_{j(2),n},
\end{equation}
where, for any $j\in\{1,..,N\}$,
\begin{equation}\label{eq-NUM2}
\mathbf{Y}_{j,n}=\mathbf{Y}_{j(1),n} + \mathbf{Y}_{j(2),n} 
\quad\mbox{ and }\quad
\mathbf{Y}_{j(i),n}:=
\mathbf{u}_{j(i)}\mathbf{v}_{j(i)}^{\top}\mathbf{Y}_{n}=
\mathbf{u}_{j(i)}\mathbf{v}_{j(i)}^{\top}\mathbf{Y}_{j,n},
\quad\mbox{ for } i=1,2.
\end{equation}
and set
\begin{equation}\label{eq-def-Ytilde}
\widetilde{\mathbf{Y}}_n:=\mathbf{Y}_{1(1),n}
=\mathbf{u}_{1(1)}\mathbf{v}_{1(1)}^{\top}\mathbf{Y}_{n}
=\begin{pmatrix}
\mathbf{u}_{1}\mathbf{v}_{1}^{\top}\mathbf{Z}_{n}\\
\mathbf{u}_{1}\mathbf{v}_{1}^{\top}\mathbf{Z}_{n}
\end{pmatrix}
=
\widetilde{Z}_n\begin{pmatrix}
\mathbf{1}\\
\mathbf{1}
\end{pmatrix}
\quad\hbox{with } \widetilde{Z}_n:=N^{-1/2}\mathbf{v}_1^{\top}Z_n,
\end{equation}
and
\begin{equation}\label{decompo-Y-hat}
\begin{split}
\widehat{\mathbf{Y}}_n&:=\mathbf{Y}_n-\widetilde{\mathbf{Y}}_n
=\mathbf{Y}_n-\mathbf{Y}_{1(1),n}
=\sum_{j=2}^N\mathbf{Y}_{j(1),n}\ +\ \sum_{j=1}^N\mathbf{Y}_{j(2),n}
\\
&=\sum_{j=2}^N\mathbf{Y}_{j(1),n}\ +\ \mathbf{Y}_{1(2),n}\ +\
\sum_{j=2}^N\mathbf{Y}_{j(2),n}.
\end{split}
\end{equation}
\begin{rem}\rm 
Note that the random vectors $\widetilde{\mathbf{Y}}_n$ and
$\widehat{\mathbf{Y}}_n$ correspond to the random vectors
$\widetilde{Z}_n(\mathbf{1},\mathbf{1})^{\top}$ and
$(\widehat{\mathbf{Z}}_n,\widehat{\mathbf{N}}_n)^{\top}$,
respectively, considered in~\cite{ale-cri-ghi-MEAN,ale-cri-ghi} in the
case $a_n=1$ for each $n$ (so that we have~$\nu=1$ and $q=1$): indeed,
we have
\begin{equation}\label{eq-NUM3}
\mathbf{1}\widetilde{Z}_n=\mathbf{u}_1\mathbf{v}_1^{\top}\mathbf{Z}_n,
\qquad
\widehat{\mathbf{Z}}_n=(I-\mathbf{u}_1\mathbf{v}_1^{\top})\mathbf{Z}_n=
\mathbf{Z}_n-\mathbf{1}\widetilde{Z}_n,
\qquad
\widehat{\mathbf{N}}_n=\mathbf{N}_n-\mathbf{1}\widetilde{Z}_n\,,
\end{equation}
where $\widetilde{Z}_n$ and $\widehat{\mathbf{Z}}_n$ are exactly the
same stochastic processes considered in
\cite{ale-cri-ghi-MEAN,ale-cri-ghi}, while $\mathbf{N}_n$ (and so
$\widehat{\mathbf{N}}_n$) differs from the stochastic process
considered in \cite{ale-cri-ghi-MEAN} because here the random
variables $N_{n,j}$ are defined in terms of a generic sequence $(a_n)$
(see \eqref{medie-empiriche-intro}) satisfying suitable assumptions.
\end{rem}
\indent Finally, it is worthwhile to point out that the decomposition
of $\widehat{\mathbf{Y}}_n$ in terms of the stochastic processes
$\mathbf{Y}_{j(i),n}$ is a new element with respect to the previous
works and, as we will see in the sequel, it will be the key tool in
order to obtain the exact convergence rate of
$\widehat{\mathbf{Y}}_n$. Indeed, the convergence rate and the
second-order asymptotic distribution of $\widehat{\mathbf{Y}}_n$ will
be the result of the different asymptotic behaviors of the three
quantities in the last term of \eqref{decompo-Y-hat}.

\subsection{Central limit theorems for $\widetilde{\mathbf{Y}}_n$ and 
$\widehat{\mathbf{Y}}_n$}
\label{section_proofs_2}
The convergence rate and the second-order asymptotic distribution of
$(\mathbf{Y}_n-Z_{\infty}\mathbf{1})$ will be obtained by studying
separately and then combining together the second-order convergence of
$\widetilde{\mathbf{Y}}_n$ to $Z_{\infty}\mathbf{1}$ and the
second-order convergence of $\widehat{\mathbf{Y}}_n$ to $\mathbf{0}$.
To this regards, we recall that, by \cite[Theorem 4.2]{ale-cri-ghi},
under Assumptions \ref{ass:W} and \ref{ass:r_n}, we have for
$1/2<\gamma\leq 1$ that
\begin{equation}\label{clt-tilde-Y}
n^{\gamma-\frac{1}{2}}(\widetilde{\mathbf{Y}}_n-Z_{\infty}\mathbf{1})
{\longrightarrow}\ \mathcal{N} \left(\ \mathbf{0}\ ,\
Z_{\infty}(1-Z_{\infty}) \widetilde{\Sigma}_{\gamma} \
\right)\quad\hbox{stably in the strong sense},
\end{equation}
 where $\widetilde{\Sigma}_{\gamma}$ is defined in
 \eqref{eq-Sigma-tilde}. In this work we fully describe the
 second-order convergence of $\widehat{\mathbf{Y}}_n$, proving the
 following theorem:

\begin{theo}\label{clt-hat-Y}
Under all the assumptions stated in Section \ref{section_model}, the
following statements hold true:
\begin{itemize}

\item[(a)] If $\nu<\gamma$, then
\begin{equation*}
n^{\nu/2}\widehat{\mathbf{Y}}_n{\longrightarrow}\ \mathcal{N}
\left(\ \mathbf{0}\ ,\ Z_{\infty}(1-Z_{\infty})\begin{pmatrix}
\mathbf{0} & \mathbf{0}\\
\mathbf{0} & \widetilde{U}S^{(q)}\widetilde{U}^{\top}\\
\end{pmatrix} \
\right)\qquad\hbox{stably},
\end{equation*}
where $\widetilde{U}$ and $S^{(q)}$ are defined in~\eqref{def:Utilde}
and \eqref{def:S_q}, respectively.

\item[(b)] If $\nu=\gamma$, then
\begin{equation*}
n^{\gamma/2}\widehat{\mathbf{Y}}_n{\longrightarrow}\ \mathcal{N}
\left(\ \mathbf{0}\ ,\ Z_{\infty}(1-Z_{\infty})
\begin{pmatrix}
\widetilde{U}S_\gamma^{11}\widetilde{U}^{\top} & 
\widetilde{U}S_\gamma^{12}\widetilde{U}^{\top} \\[3pt]
\widetilde{U}S_\gamma^{21}\widetilde{U}^{\top} & 
\widetilde{U}S_\gamma^{22}\widetilde{U}^{\top}
\end{pmatrix}
\ \right)\qquad\hbox{stably},
\end{equation*}
where $S^{21}_\gamma=(S^{12}_\gamma)^{\top}$ and, for $2\leq j_1,j_2\leq N$,
\begin{equation*}\label{eq-S_gamma}
\begin{split}
&[S^{11}_\gamma]_{11} =\ [S^{11}]_{j_11}\ =\ [S^{11}]_{1j_2}\ :=\ 0,\\
&[S^{11}_\gamma]_{j_1 j_2} :=\ \frac{c^2}{c(\alpha_{j_1}+\alpha_{j_2})-
\ind_{\{\gamma=1\}}}
\mathbf{v}_{j_1}^{\top}\mathbf{v}_{j_2},\\
&[S^{12}_\gamma]_{11}\ =\ [S^{12}]_{1j_2}\ :=\ 0,\\
&[S^{12}_\gamma]_{j_11}\ :=\ 
\frac{c(q-c)}{c\alpha_{j_1}+q-\ind_{\{\gamma=1\}}}
\mathbf{v}_{j_1}^{\top}\mathbf{v}_{1},\\
&[S^{12}_\gamma]_{j_1 j_2} :=\
\frac{cq(c\alpha_{j_1}+c-\ind_{\{\gamma=1\}})}
{(c\alpha_{j_1}+c\alpha_{j_2}-\ind_{\{\gamma=1\}})
(c\alpha_{j_1}+q-\ind_{\{\gamma=1\}})}
\mathbf{v}_{j_1}^{\top}\mathbf{v}_{j_2},\\
&[S^{22}_\gamma]_{11} :=\ \frac{(q-c)^2}{2q-
\ind_{\{\gamma=1\}}}\|\mathbf{v}_{1}\|^2,\\
&[S^{22}_\gamma]_{j1}=[S^{22}_\gamma]_{1j}:=\ 
\frac{q(q-c)(c+q-\ind_{\{\gamma=1\}})}{(c\alpha_{j}+q
-\ind_{\{\gamma=1\}})(2q-\ind_{\{\gamma=1\}})}\mathbf{v}_{j}^{\top}\mathbf{v}_{1},\\
&[S^{22}_\gamma]_{j_1 j_2} :=\ q^2\frac{c^3(\alpha_{j_1}+\alpha_{j_2})+
2c^2q(\alpha_{j_1}\alpha_{j_2}+1)-
\ind_{\{\gamma=1\}}c^2(\alpha_{j_1}\alpha_{j_2}+\alpha_{j_1}+\alpha_{j_2}+2)}
{(2q-\ind_{\{\gamma=1\}})(c(\alpha_{j_1}+\alpha_{j_2})-
\ind_{\{\gamma=1\}})
(c\alpha_{j_1}+q-\ind_{\{\gamma=1\}})(c\alpha_{j_2}+q-\ind_{\{\gamma=1\}})}
\mathbf{v}_{j_1}^{\top}\mathbf{v}_{j_2}\\
& +\ q^2\frac{c(q-\ind_{\{\gamma=1\}})^2(\alpha_{j_1}+\alpha_{j_2})-
\ind_{\{\gamma=1\}}(2c+q-1)(q-1)}
{(2q-\ind_{\{\gamma=1\}})(c(\alpha_{j_1}+\alpha_{j_2})-
\ind_{\{\gamma=1\}})(c\alpha_{j_1}+q-
\ind_{\{\gamma=1\}})(c\alpha_{j_2}+q-\ind_{\{\gamma=1\}})}
\mathbf{v}_{j_1}^{\top}\mathbf{v}_{j_2}\,.
\end{split}
\end{equation*}

\item[(c)] If $\gamma<\nu$, then
\begin{equation*}
n^{\gamma-\frac{\nu}{2}}\widehat{\mathbf{Y}}_n{\longrightarrow}\
\mathcal{N} 
\left(\ \mathbf{0}\ ,\
Z_{\infty}(1-Z_{\infty}) \frac{c^2}{N[2q-\ind_{\{\nu=1\}}(2\gamma-1)]} 
\|\mathbf{v}_{1}\|^2
\begin{pmatrix}
\mathbf{0} & \mathbf{0}\\
\mathbf{0} & \mathbf{1}\mathbf{1}^{\top}\\
\end{pmatrix} \ 
\right)
\qquad\hbox{stably}.
\end{equation*}
\end{itemize}
\end{theo}

\begin{rem}
\rm Note that, when $\nu\neq \gamma$ the convergence rates of the
first and the second component of $\widehat{\mathbf{Y}}_n$ are always
different: indeed from \cite{ale-cri-ghi}, we know that, under our
assumptions, the convergence rate of $\widehat{\mathbf{Z}}_n$ is
always $n^{\gamma/2}$, while the above theorem shows that the convergence
rate of $\widehat{\mathbf{N}}_n$ changes according to the pair
$(\gamma,\nu)$.
\end{rem}

Regarding the proof of Theorem \ref{clt-hat-Y}, we note that, using
the definition \eqref{decompo-Y-hat} of $\widehat{\mathbf{Y}}_n$ given
in Section~\ref{section_decomposition}, we can say that this random
variable can be decomposed in a sum of suitable random variables that
have the form
\begin{equation}\label{general-sum}
\sum_{j\in J}\sum_{i\in I_j}\mathbf{Y}_{j(i),n},
\end{equation}
where $J\underline\subset \{1,\cdots,N\}$, $I_j\underline\subset
\{1,2\}$ for any $j\in J$ and $Y_{j(i),n}$ is defined in
\eqref{eq-NUM2}. Hence, in order to characterize the asymptotic
behavior of $\widehat{\mathbf{Y}}_n$, we first establish the
second-order asymptotic behavior of the above general sum
\eqref{general-sum} under certain specifications of the sets $J$ and
$I_j$ (see Lemma \ref{conv-general-sum}) and then we combine them
together appropriately according to their convergence rates.

\begin{lem}\label{conv-general-sum}
Under all the assumptions stated in Section \ref{section_model},
consider the general sum \eqref{general-sum} in the following cases:
\begin{itemize}
\item[$(i)$] $\nu<\gamma$, $J=\{2,\cdots,N\}$ and
  $I_j=\{1\}$ for all $j\in J$;
\item[$(ii)$] $\nu<\gamma$, $J=\{1,\cdots,N\}$ and
  $I_j=\{2\}$ for all $j\in J$;
\item[$(iii)$] $\nu=\gamma$,
  $J=\{1,\cdots,N\}$, $I_1=\{2\}$ and $I_j=\{1,2\}$ for all $j\in
  J\setminus\{1\}$;
\item[$(iv)$] $\nu>\gamma$,
  $J=\{2,\cdots,N\}$ and $I_j=\{1\}$ for all $j\in J$;
\item[$(v)$] $\nu>\gamma$, 
  $J=\{1\}$ and $I_1=\{2\}$;
\item[$(vi)$] $\nu>\gamma$, 
  $J=\{2,\cdots,N\}$ and $I_j=\{2\}$ for all $j\in J$.
\end{itemize}
Then, in all the above listed cases,  we have
\begin{equation}
t_n(J(I))\sum_{j\in J}\sum_{i\in I_j}\mathbf{Y}_{j(i),n}
\stackrel{stably}\longrightarrow
{\mathcal N}\left(\mathbf{0},
Z_{\infty}(1-Z_{\infty})
\sum_{j_1\in J}\sum_{j_2\in J}\mathbf{v}_{j_1}^{\top}\mathbf{v}_{j_2}
\sum_{i_1\in I_{j_1}}\sum_{i_2\in I_{j_2}}
d^{j_1(i_1),j_2(i_2)}\mathbf{u}_{j_1(i_1)}\mathbf{u}_{j_2(i_2)}^{\top}\right)\,,
\end{equation}
where 
\begin{equation}\label{rate-general-sum}
t_n(J(I)):=
\begin{cases}
n^{\gamma/2} &\hbox{for cases (i), (iii) and (iv)}
\\
n^{\nu/2} &\hbox{for cases (ii) and (vi)}
\\
n^{\gamma-\frac{\nu}{2}} &\hbox{for case (v)}\,.
\end{cases}
\end{equation}
and $d^{j_1(i_1),j_2(i_2)}$ are constants corresponding to the result 
of suitable limits computed in Section \ref{result-limit-d} of the
appendix.
\end{lem}

\noindent{\bf Proof of Theorem \ref{clt-hat-Y}.} From the above lemma,
we immediately get the proof of Theorem \ref{clt-hat-Y}. Indeed, in
case (a) we get
$$ 
n^{\nu/2}\widehat{\mathbf{Y}}_n
= \frac{1}{ n^{(\gamma-\nu)/2} } n^{\gamma/2}\sum_{j=2}^N\mathbf{Y}_{j(1),n}
+ n^{\nu/2}\sum_{j=1}^N\mathbf{Y}_{j(2),n} 
$$ 
where, considering the above cases $(i)$ and $(ii)$, the first term
in the sum converges in probability to zero, while the second term
converges stably to the desired Gaussian kernel, that is the
Gaussian kernel with zero mean and random covariance matrix
$$
Z_{\infty}(1-Z_{\infty})
\sum_{j_1=1}^N\sum_{j_2=1}^N
\mathbf{v}_{j_1}^{\top}\mathbf{v}_{j_2}
d^{j_1(2),j_2(2)}\mathbf{u}_{j_1(2)}\mathbf{u}_{j_2(2)}^{\top},
$$
where
\begin{equation*}
\mathbf{u}_{j_1(2)}\mathbf{u}_{j_2(2)}^{\top}=
\begin{pmatrix}
\mathbf{0} & \mathbf{0}\\
\mathbf{0} & \mathbf{u}_{j_1}\mathbf{u}_{j_2}^{\top}
\end{pmatrix}\,.
\end{equation*}
In case (b) we simply have 
$$
n^{\gamma/2}\widehat{\mathbf{Y}}_n=
n^{\gamma/2}
\left(\sum_{j=2}^N\mathbf{Y}_{j(1),n}+\sum_{j=1}^N\mathbf{Y}_{j(2),n}\right)\,,
$$ where the right-hand term converges stably to the desired Gaussian
kernel (see the above case $(iii)$), that is the Gaussian kernel with
zero mean and random covariance matrix
\begin{equation*}
Z_{\infty}(1-Z_{\infty})
\sum_{j_1=1}^N\sum_{j_2=1}^N\sum_{i_1=1}^2\sum_{i_2=1}^2
(1-\ind_{\{j_1=i_1=1\}}\ind_{\{j_2=i_2=1\}})
\mathbf{v}_{j_1}^{\top}\mathbf{v}_{j_2}
d^{j_1(i_1),j_2(i_2)}\mathbf{u}_{j_1(i_1)}\mathbf{u}_{j_2(i_2)}^{\top}\,,
\end{equation*}
where
\begin{equation*}
\begin{aligned}
\mathbf{u}_{j_1(1)}\mathbf{u}_{j_2(1)}^{\top}&&=&\
\begin{pmatrix}
\mathbf{u}_{j_1}\mathbf{u}_{j_2}^{\top} &
\ind_{\{j_2=1\}}\mathbf{u}_{j_1}\mathbf{u}_{j_2}^{\top}\\
\ind_{\{j_1=1\}}\mathbf{u}_{j_1}\mathbf{u}_{j_2}^{\top} &
\ind_{\{j_1=1\}}\ind_{\{j_2=1\}}\mathbf{u}_{j_1}\mathbf{u}_{j_2}^{\top}
\end{pmatrix},\\
\mathbf{u}_{j_1(1)}\mathbf{u}_{j_2(2)}^{\top}&&=&\
\begin{pmatrix}
\mathbf{0} & \mathbf{u}_{j_1}\mathbf{u}_{j_2}^{\top}\\
\mathbf{0} & \ind_{\{j_1=1\}}\mathbf{u}_{j_1}\mathbf{u}_{j_2}^{\top}
\end{pmatrix},\\
\mathbf{u}_{j_1(2)}\mathbf{u}_{j_2(1)}^{\top}&&=&\
\begin{pmatrix}
\mathbf{0} & \mathbf{0}\\
\mathbf{u}_{j_1}\mathbf{u}_{j_2}^{\top} &
\ind_{\{j_2=1\}}\mathbf{u}_{j_1}\mathbf{u}_{j_2}^{\top}
\end{pmatrix},\\
\mathbf{u}_{j_1(2)}\mathbf{u}_{j_2(2)}^{\top}&&=&\
\begin{pmatrix}
\mathbf{0} & \mathbf{0}\\
\mathbf{0} & \mathbf{u}_{j_1}\mathbf{u}_{j_2}^{\top}
\end{pmatrix}.
\end{aligned}
\end{equation*}
Finally, in case (c), we obtain 
$$
n^{ \gamma-\frac{\nu}{2} }\widehat{\mathbf{Y}}_n=
\frac{1}{ n^{(\nu-\gamma)/2} }
n^{\frac{\gamma}{2}}\sum_{j=2}^N\mathbf{Y}_{j(1),n}
+
n^{ \gamma-\frac{\nu}{2} }\mathbf{Y}_{1(2),n}+
\frac{1}{ n^{(\nu-\gamma)} }
n^{\frac{\nu}{2}}\sum_{j=2}^N\mathbf{Y}_{j(2),n}\,,
$$ 
where, considering the above cases $(iv)$, $(v)$ and $(vi)$, we have that the
first and the third terms in the sum converge in probability to zero,
while the second term converges stably to the desired Gaussian kernel,
that is the Gaussian kernel with zero mean and random covariance matrix
\begin{equation*}
Z_{\infty}(1-Z_{\infty})
\|\mathbf{v}_{1}\|^2
d^{1(2),1(2)}\mathbf{u}_{1(2)}\mathbf{u}_{1(2)}^{\top}\,,
\end{equation*}
where
\begin{equation}\label{equality-matrix}
\mathbf{u}_{1(2)}\mathbf{u}_{1(2)}^{\top}=
\begin{pmatrix}
\mathbf{0} & \mathbf{0}\\
\mathbf{0} & \mathbf{u}_{1}\mathbf{u}_{1}^{\top}
\end{pmatrix}=
\frac{1}{N}
\begin{pmatrix}
\mathbf{0} & \mathbf{0}\\
\mathbf{0} & \mathbf{1}\mathbf{1}^{\top}
\end{pmatrix}.
\end{equation}
\qed

\noindent We now go on with the proof of Lemma \ref{conv-general-sum}.
\\
\noindent{\bf Proof of Lemma \ref{conv-general-sum}.}  Since this
proof is quite long, we split it into various steps and the technical
computations and details are collected in the appendix.
\\
\noindent{\bf First step: decomposition of the general sum
  \eqref{general-sum}.}\\  
\noindent First of all, we observe that, for any set
$J\underline\subset \{1,\cdots,N\}$, the dynamics of $\sum_{j\in
  J}\mathbf{Y}_{j,n}$ can be obtained by summing up
equation~\eqref{eq:proof-dynamics-Y_each_j} for $j\in J$:
\begin{equation*}
\sum_{j\in J}\mathbf{Y}_{j,n}\ =\ \left(\sum_{j\in J}
U_j\left(I-D_{Q,j,n}\right)V_j^{\top}\right)\sum_{j\in
J}\mathbf{Y}_{j,n-1}\ + \left(\sum_{j\in
J}U^{*}_jD_{R,n}V_j^{*\top}\right)\Delta\mathbf{M}_{Y,n}+
\left(\sum_{j\in J}U_jV_j^{\top}\right)\Delta\mathbf{R}_{Y,n}.
\end{equation*}
Then, recalling that $\mathcal{R}e(\alpha_j) > 0$ for each $j\geq 2$
because $\mathcal{R}e(\lambda_j) < 1$ for each $j\geq 2$, and taking
an integer $m_0\geq 2$ large enough such that for $n \geq m_0$ we have
$\mathcal{R}e(\alpha_j)cn^{-\gamma}< 1$ for each $j\geq 2$ and
$qn^{-\nu} < 1$, we can write
\begin{equation}\label{eq:above_equality}
\begin{aligned}
\sum_{j\in J}\mathbf{Y}_{j,n}&= \left(\sum_{j\in J}U_jA^j_{m_0,n-1}
V_j^{\top}\right)\sum_{j\in J}\mathbf{Y}_{j,m_0}\\
&+\sum_{k=m_0}^{n-1}\left(\sum_{j\in J}
U_jA^j_{k+1,n-1}V_j^{\top}U^{*}_jD_{R,k}V_j^{*\top} \right)
\Delta\mathbf{M}_{Y,k+1}\\
&+\sum_{k=m_0}^{n-1}
\left(\sum_{j\in J}U_jA^j_{k+1,n-1}V_j^{\top}\right)\Delta\mathbf{R}_{Y,k+1}
\qquad\hbox{for } n\geq m_0\,,
\end{aligned}
\end{equation}
where, for any $j\in J$,
\begin{equation}\label{def:A_k_n_withoutG}
A^j_{k+1,n-1}=
\begin{cases}
&\prod_{m=k+1}^{n-1}\left(I-D_{Q,j,m}\right)\quad\mbox{for } m_0\leq k\leq n-2\\
& I\quad\mbox{for } k=n-1.
\end{cases}
\end{equation}
Setting for any $x=a_x+ib_x\in\mathbb{C}$ with $a_x>0$ and
$1/2<\delta\leq 1 $,
\begin{equation*}
p^{\delta}_{k}(x):=\prod_{m=m_0}^{k}\left(1-\frac{x}{m^{\delta}}\right)
\;\hbox{for } k\geq m_0
\end{equation*}
and 
\begin{equation*}
F^{\gamma}_{k+1,n-1}(x):=\frac{p^{\gamma}_{n-1}(x)}{p^{\gamma}_{k}(x)}
\;\hbox{for } m_0\leq k\leq n-1,
\end{equation*}
It is easy to see that, for $j=1$, we have 
\begin{equation}\label{eq:A_1_k_n}
A^1_{k+1,n-1}=\begin{pmatrix}
1 & 0 \\
0 & F^{\nu}_{k+1,n-1}(q)
\end{pmatrix}
\quad\mbox{for } m_0\leq k\leq n-1
\end{equation}
and, for $j\geq 2$, after some calculations reported in
Section~\ref{subsection_A_FOGLIO_2} of the appendix, we obtain
\begin{equation}\label{eq:A_j_k_n_withG}
A^j_{k+1,n-1}=\begin{pmatrix}
F^{\gamma}_{k+1,n-1}(c\alpha_j) & 0 \\
\lambda_jG_{k+1,n-1}(c\alpha_j,q) & F^{\nu}_{k+1,n-1}(q)
\end{pmatrix}
\quad\mbox{for } m_0\leq k\leq n-1,
\end{equation}
where
\begin{equation}\label{eq:G}
G_{k+1,n-1}(x,q)\ :=\
\sum_{l=k+1}^{n-1}F^{\gamma}_{l+1,n-1}(x)h_l(1-c^{-1}x)F^{\nu}_{k+1,l-1}(q).
\end{equation}
Then, since $V_j^{\top}U^{*}_j=P_j^{-1}$,
equation~\eqref{eq:above_equality} can be rewritten as
\begin{equation}\label{eq:dynamic-Y-step1}
\sum_{j\in J}\mathbf{Y}_{j,n}= \left(\sum_{j\in J}
U_jA^j_{m_0,n-1} V_j^{\top}\right)\sum_{j\in J}\mathbf{Y}_{j,m_0}+
\sum_{k=m_0}^{n-1}\mathbf{T}^J_{k+1,n-1}+
\sum_{k=m_0}^{n-1}\boldsymbol{\rho}^J_{k+1,n-1} \quad\hbox{for } n\geq
m_0,
\end{equation}
with
\begin{equation*}
\begin{aligned}
\mathbf{T}^J_{k+1,n-1}&&=&\ \left( \sum_{j\in J}
U_jA^j_{k+1,n-1}P_j^{-1}D_{R,k}V_j^{*\top}\right)\Delta\mathbf{M}_{Y,k+1},\\
\boldsymbol{\rho}^J_{k+1,n-1}&&=&\
\left(\sum_{j\in J}U_jA^j_{k+1,n-1}V_j^{\top}\right)\Delta\mathbf{R}_{Y,k+1}.
\end{aligned}
\end{equation*}
In order to get a similar decomposition for the general sum
\eqref{general-sum}, we set, for any $j\in J$, 
\begin{equation}\label{U0}
U^0_{j(1)}:=\begin{pmatrix} \mathbf{u}_{j(1)} & \mathbf{0}
\end{pmatrix}=\begin{pmatrix} \mathbf{u}_{j} & \mathbf{0}\\
g(\lambda_j)\mathbf{u}_{j} & \mathbf{0}\\
\end{pmatrix}\qquad\mbox{ and }\qquad
U^0_{j(2)}:=\begin{pmatrix} \mathbf{0} & \mathbf{u}_{j(2)}
\end{pmatrix}=\begin{pmatrix} \mathbf{0} & \mathbf{0}\\
\mathbf{0} &\mathbf{u}_{j}
\end{pmatrix}
\end{equation}
and taking into account the last relation in \eqref{eq-NUM2}, 
we get
\begin{equation}\label{eq-decomposizione-general-sum}
\sum_{j\in J}\sum_{i\in I_j}\mathbf{Y}_{j(i),n}=
C^{J(I)}_{m_0,n-1}\sum_{j\in J}\mathbf{Y}_{j,m_0}+
\sum_{k=m_0}^{n-1}\mathbf{T}^{J(I)}_{k+1,n-1}+
\sum_{k=m_0}^{n-1}\boldsymbol{\rho}^{J(I)}_{k+1,n-1} \quad\hbox{for }
n\geq m_0,
\end{equation}
with
\begin{equation*}
\begin{aligned}
C^{J(I)}_{m_0,n-1}&&=&\ \sum_{j\in J}\sum_{i\in I_j}
U^0_{j(i)}A^j_{m_0,n-1} V_j^{\top},\\
\mathbf{T}^{J(I)}_{k+1,n-1}&&=&\ \left( \sum_{j\in J}\sum_{i\in I_j}
U^0_{j(i)}A^j_{k+1,n-1}P_j^{-1}D_{R,k}V_j^{*\top}\right)\Delta\mathbf{M}_{Y,k+1},\\
\boldsymbol{\rho}^{J(I)}_{k+1,n-1}&&=&\
\left(\sum_{j\in J}\sum_{i\in I_j}U^0_{j(i)}A^j_{k+1,n-1}V_j^{\top}\right)
\Delta\mathbf{R}_{Y,k+1}.
\end{aligned}
\end{equation*}
In the sequel of the proof, we will establish the asymptotic behavior
of the general sum \eqref{general-sum} by studying separately the
three terms $C_{m_0,n-1}^{J(I)}\sum_{j\in J}\mathbf{Y}_{j,m_0}$,
$\sum_{k=m_0}^{n-1}\mathbf{T}^{J(I)}_{k+1,n-1}$ and
$\sum_{k=m_0}^{n-1}\boldsymbol{\rho}^{J(I)}_{k+1,n-1}$ in the six
cases $(i)$-$(vi)$ specified in the statement of the considered lemma.
\\
\noindent{\bf Second step: asymptotic behavior of
  $C_{m_0,n-1}^{J(I)}\sum_{j\in J}\mathbf{Y}_{j,m_0}$.}\\
\noindent From~\eqref{eq:A_1_k_n}, \eqref{eq:A_j_k_n_withG},
\eqref{eq:G} and~\eqref{U0}, taking into account the fact that  
 in all the considered cases with $1\in J$,
i.e.~$(ii), (iii)$ and $(v)$, we have $1\notin I_1$, we get 
\begin{equation*}
\left|C^{J(I)}_{m_0,n}\right|
\ =\ O(C^{11}_{n}) + O(C^{21}_{n}) + O(C^{22}_{n}),
\end{equation*}
where
\begin{equation*}
\begin{aligned}
C^{11}_{n}\ &&:=&\ \sum_{j\in J,\, j\neq 1}\ind_{\{1\in
I_j\}}|F_{m_0,n-1}^{\gamma}(c\alpha_j)|\,,
\\
C^{21}_{n}\ &&:=&\ \sum_{j\in J,\, j\neq 1} \ind_{\{2\in I_j\}}
|G_{m_0,n-1}(c\alpha_j,q)|\,,
\\
C^{22}_{n}\ &&:=&\ \sum_{j\in J} \ind_{\{2\in I_j\}}
|F_{m_0,n-1}^{\nu}(q)|\,.
\end{aligned}
\end{equation*}
Using~\eqref{affermazione1_F} in the appendix and denoting by $a^*$
the real part of $\alpha^*:=1-\lambda^*$, it is immediate to see that
$$
C^{11}_{n}=\sum_{j\in J,\, j\neq 1}\ind_{\{1\in
I_j\}}
\left\{\begin{aligned}
& O\left(\exp\left(-ca^{*}\frac{n^{1-\gamma}}{1-\gamma}\right)\right)
&\mbox{ if } 1/2<\gamma<1\\
&O(n^{-ca^{*}})
&\mbox{ if }\gamma=1
\end{aligned}\right.
$$
and
$$
C^{22}_{n}=\sum_{j\in J} \ind_{\{2\in I_j\}}
\left\{\begin{aligned}
&\ O\left(\exp\left(-q\frac{n^{1-\nu}}{1-\nu}\right)\right)
&\mbox{ if }1/2<\nu<1
\\
&\ O(n^{-q})&\mbox{ if }\nu=1.
\end{aligned}\right.
$$

For the term $C^{21}_{n}$, we apply
Lemma \ref{lemma-tecnico-G} so that we get:
\begin{description}
\item[Case $\nu<\gamma$] We have $G_{m_0,n-1}(c\alpha_j,q)= O\big(
  n^{-(\gamma-\nu)}|F^{\nu}_{m_0,n-1}(q)|+|F^{\gamma}_{m_0,n-1}(c\alpha_j)|\big)$
  by means of Lemma \ref{lemma-tecnico-G} and so
$$
C^{21}_{n}=
\sum_{j\in J,\,j\neq 1} \ind_{\{2\in I_j\}}
O\left( n^{-(\gamma-\nu)}|F^{\nu}_{m_0,n-1}(q)|+|F^{\gamma}_{m_0,n-1}(c\alpha_j)|
\right),
$$ 
where, as above, by \eqref{affermazione1_F}, we have
$|F^{\nu}_{m_0,n-1}(q)|=O\left(\exp\left(-q\frac{n^{1-\nu}}{1-\nu}\right)\right) 
$ 
and 
$$
|F^{\gamma}_{m_0,n-1}(c\alpha_j)|=\left\{\begin{aligned}
& O\left(\exp\left(-ca^{*}\frac{n^{1-\gamma}}{1-\gamma}\right)\right)
&\mbox{ if } 1/2<\gamma<1\\
&O(n^{-ca^{*}})
&\mbox{ if }\gamma=1.
\end{aligned}\right.
$$

\item[Case $\nu>\gamma$] We have $G_{m_0,n-1}(c\alpha_j,q)=
  O\big(n^{-(\nu-\gamma)}|F^{\nu}_{m_0,n-1}(q)|+
  |F^{\gamma}_{m_0,n-1}(c\alpha_j)|\big)$ by means of Lemma
  \ref{lemma-tecnico-G} and so 
$$
C^{21}_{n}=\sum_{j\in J,\,j\neq 1} \ind_{\{2\in I_j\}}
O\left(n^{-(\nu-\gamma)}|F^{\nu}_{m_0,n-1}(q)|+
|F^{\gamma}_{m_0,n-1}(c\alpha_j)|
\right),
$$
where, as above, by \eqref{affermazione1_F}, we have
$|F^{\gamma}_{m_0,n-1}(c\alpha_j)|=
O\left(\exp\left(-ca^{*}\frac{n^{1-\gamma}}{1-\gamma}\right)\right)$
 and 
$$
|F^{\nu}_{m_0,n-1}(q)|=
\left\{\begin{aligned}
&\ O\left(\exp\left(-q\frac{n^{1-\nu}}{1-\nu}\right)\right)
&\mbox{ if }1/2<\nu<1
\\
&\ O(n^{-q})&\mbox{ if }\nu=1.
\end{aligned}\right.
$$

\item[Case $\nu=\gamma$] By assumption \eqref{ass-q} \footnote{If
  there exists $j$ such that $q=c\alpha_j$, we have to consider the
  other asymptotic expression given in Lemma \ref{lemma-tecnico-G}.} 
  and Lemma \ref{lemma-tecnico-G}, we have $G_{m_0,n-1}(c\alpha_j,q)=
  O\big(|F^{\gamma}_{m_0,n-1}(q)|+
  |F^{\gamma}_{m_0,n-1}(c\alpha_j)|\big)$ and so
$$
C^{21}_{n}=\sum_{j\in J,\,j\neq 1} \ind_{\{2\in I_j\}}
O\left(
|F^{\gamma}_{m_0,n-1}(q)|+|F^{\gamma}_{m_0,n-1}(c\alpha_j)|
\right),
$$
where, as above, by \eqref{affermazione1_F}, we have for 
$x=q$ or $x\in\{c\alpha_j:\, j\in J,\,j\neq 1\}$
$$
|F^{\gamma}_{m_0,n-1}(x)|=
\left\{\begin{aligned}
&\ O\left(\exp\left(-a_x\frac{n^{1-\gamma}}{1-\gamma}\right)\right)
&\mbox{ if }1/2<\nu=\gamma<1
\\
&\ O(n^{-a_x})&\mbox{ if }\nu=\gamma=1
\end{aligned}\right.
$$
and so, setting $x^*:=\min\{q,ca^*\}$, we can write  
$$
C^{21}_n=\sum_{j\in J,\,j\neq 1} \ind_{\{2\in I_j\}}
\left\{\begin{aligned}
& O\left(\exp\left(-x^{*}\frac{n^{1-\gamma}}{1-\gamma}\right)\right)
&\mbox{ if } 1/2<\nu=\gamma<1\\
&O(n^{-x^{*}})
&\mbox{ if }\nu=\gamma=1.
\end{aligned}\right.
$$ 
\end{description}

Summing up, taking into account the conditions $ca^*>1/2$ when
$\gamma=1$ and $q>1/2$ when $\nu=1$, we can conclude that in all the
six cases $(i)$-$(vi)$ we have
$t_n(J(I))\left|C_{m_0,n-1}^{J(I)}\right|\to 0$ and so
$$
t_n(J(I))C_{m_0,n-1}^{J(I)}\sum_{j\in J}\mathbf{Y}_{j,m_0}
\stackrel{a.s.}\longrightarrow \mathbf{0}.
$$
\noindent{\bf Third step: asymptotic behavior of
  $\sum_{k=m_0}^{n-1}\boldsymbol{\rho}^{J(I)}_{k+1,n-1}$ .}\\
\noindent We recall that, by Assumption~\ref{ass:r_n}, we have
$|\Delta\mathbf{R}_{Z,k+1}|=O(k^{-2\gamma})$ and
$|\Delta\mathbf{R}_{N,k+1}|=O(k^{-2\nu})$.  Then,
from~\eqref{eq:A_1_k_n}, \eqref{eq:A_j_k_n_withG}, \eqref{eq:G}
and~\eqref{U0}, taking into account the fact that in all the
considered cases with $1\in J$, i.e.~$(ii), (iii)$ and $(v)$, we have
$1\notin I_1$, we get
\begin{equation*}
\left|\sum_{k=m_0}^{n-1}\boldsymbol{\rho}^{J(I)}_{k+1,n-1}\right|\ =\
O(\rho^{11}_{n}) + O(\rho^{21}_{n}) + O(\rho^{22}_{n}),
\end{equation*}
where
\begin{equation*}
\begin{split}
&\rho^{11}_{n}:=
\sum_{j\in J,\,j\neq 1} \ind_{\{1\in I_j\}}
\sum_{k=m_0}^{n-1}k^{-2\gamma}|F_{k+1,n-1}^{\gamma}(c\alpha_j)|\,,
\\
&\rho^{21}_{n}:=\sum_{j\in J,\,j\neq 1} \ind_{\{2\in I_j\}}
\sum_{k=m_0}^{n-1}k^{-2\gamma}|G_{k+1,n-1}(c\alpha_j,q)|\,,
\\
&\rho^{22}_{n}:=\sum_{j\in J} \ind_{\{2\in I_j\}}
\sum_{k=m_0}^{n-1}(k^{-2\gamma}+k^{-2\nu})|F_{k+1,n-1}^{\nu}(q)|.
\end{split}
\end{equation*}
Using Lemma
\ref{lemma-tecnico-O} (with $\beta=2\gamma>1$, $e=1$ and
$\delta=\gamma$), we get
\begin{equation*}
\rho^{11}_{n}=\sum_{j\in J,\,j\neq 1} \ind_{\{1\in I_j\}}
\begin{cases}
O\left(n^{-\gamma}\right)
\qquad &\hbox{if } 1/2<\gamma<1,\\
O\left(n^{-ca^*}\right)
\qquad &\hbox{if } \gamma=1\;\hbox{and } 1/2< ca^*<1,\\
O\left(n^{-1}\ln(n)\right)
\qquad &\hbox{if } \gamma=1\;\hbox{and } ca^*=1,\\
O\left(n^{-1}\right)
\qquad &\hbox{if } \gamma=1\;\hbox{and } ca^*>1.
\end{cases}
\end{equation*}
For $\rho^{22}_n$, we observe that we have $k^{-2\gamma}=O(k^{-2\nu})$
when $\nu\leq \gamma$ and $k^{-2\nu}=O(k^{-2\gamma})$ when
$\nu>\gamma$. Therefore, using Lemma \ref{lemma-tecnico-O} (with $e=1$
and $\delta=\nu$ and $\beta=2\nu>1$ if $\nu\leq\gamma$ and
$\beta=2\gamma>1$ if $\nu>\gamma$), we obtain for the case $\nu\leq\gamma$
\begin{equation*}
\begin{split}
\rho^{22}_{n}&=
\sum_{j\in J} \ind_{\{2\in I_j\}}
O\left(\sum_{k=m_0}^{n-1}k^{-2\nu}|F_{k+1,n-1}^{\nu}(q)|\right)
\\
&=
\sum_{j\in J} \ind_{\{2\in I_j\}}
\begin{cases}
O\left(n^{-\nu}\right)
\qquad &\hbox{if } 1/2<\nu<1,\\
O\left(n^{-q}\right)
\qquad &\hbox{if } \nu=1\;\hbox{and } 1/2<q<1,\\
O\left(n^{-1}\ln(n)\right)
\qquad &\hbox{if } \nu=1\;\hbox{and } q=1,\\
O\left(n^{-1}\right)
\qquad &\hbox{if } \nu=1\;\hbox{and } q>1\,,
\end{cases}
\end{split}
\end{equation*}
and for the case $\nu>\gamma$
\begin{equation}\label{recall-rho}
\begin{split}
\rho^{22}_{n}&=
\sum_{j\in J} \ind_{\{2\in I_j\}}
O\left(\sum_{k=m_0}^{n-1}k^{-2\gamma}|F_{k+1,n-1}^{\nu}(q)|\right)
\\
&=
\sum_{j\in J} \ind_{\{2\in I_j\}}
\begin{cases}
O\left(n^{-2\gamma+\nu}\right)
\qquad &\hbox{if } 1/2<\nu<1,\\
O\left(n^{-q}\right)
\qquad &\hbox{if } \nu=1\;\hbox{and } 1/2<q<2\gamma-1,\\
O\left(n^{-q}\ln(n)\right)
\qquad &\hbox{if } \nu=1\;\hbox{and } q=2\gamma-1>1/2,\\
O\left(n^{-2\gamma+1}\right)
\qquad &\hbox{if } \nu=1\;\hbox{and } q>\max\{1/2,2\gamma-1\}\,.
\end{cases}
\end{split}
\end{equation}
For the term $\rho^{21}_{n}$, we apply Lemma \ref{lemma-tecnico-O} and
Lemma \ref{lemma-tecnico-G} so that we get:
\begin{description}
\item[Case $\nu<\gamma$] We have $G_{k+1,n-1}(c\alpha_j,q)=O\big(
  n^{-(\gamma-\nu)}|F^{\nu}_{k+1,n-1}(q)|+k^{-(\gamma-\nu)}|F^{\gamma}_{k+1,n-1}(c\alpha_j)|
\big)$
  by means of Lemma \ref{lemma-tecnico-G}, and so we get
$$
\rho^{21}_{n}= \sum_{j\in J\,j\neq 1} \ind_{\{2\in I_j\}}
O\left(n^{-(\gamma-\nu)}
\sum_{k=m_0}^{n-1}
\frac{1}{k^{2\gamma}}|F^{\nu}_{k+1,n-1}(q)|
+
\sum_{k=m_0}^{n-1}
\frac{1}{k^{3\gamma-\nu}}|F^{\gamma}_{k+1,n-1}(c\alpha_j)|
\right)\,,
$$
where, by Lemma \ref{lemma-tecnico-O}, the first term is 
$O(n^{-3\gamma+2\nu})$, while for the second term we have 
$$
\sum_{k=m_0}^{n-1}
\frac{1}{k^{3\gamma-\nu}}|F^{\gamma}_{k+1,n-1}(c\alpha_j)|
=
\begin{cases}
O\left(n^{-2\gamma+\nu}\right)
\qquad &\hbox{if } 1/2<\gamma<1,\\
O\left(n^{-ca^*}\right)
\qquad &\hbox{if } \gamma=1\;\hbox{and } 1/2< ca^*<2-\nu,\\
O\left(n^{-2+\nu}\ln(n)\right)
\qquad &\hbox{if } \gamma=1\;\hbox{and } ca^*=2-\nu,\\
O\left(n^{-2+\nu}\right)
\qquad &\hbox{if } \gamma=1\;\hbox{and } ca^*>2-\nu\,.
\end{cases}
$$

\item[Case $\nu>\gamma$] We have
  $G_{k+1,n-1}(c\alpha_j,q)=O\big(n^{-(\nu-\gamma)}|F^{\nu}_{k+1,n-1}(q)|+
  k^{-(\nu-\gamma)}|F^{\gamma}_{k+1,n-1}(c\alpha_j)|\big)$ by means of
  Lemma \ref{lemma-tecnico-G}, and so we get 
$$
\rho^{21}_{n}=\sum_{j\in J,\,j\neq 1} \ind_{\{2\in I_j\}}
O\left(
n^{-(\nu-\gamma)}
\sum_{k=m_0}^{n-1}
\frac{1}{k^{2\gamma}}|F^{\nu}_{k+1,n-1}(q)|
+
\sum_{k=m_0}^{n-1}
\frac{1}{k^{\gamma+\nu}}|F^{\gamma}_{k+1,n-1}(c\alpha_j)|\right)
\,,
$$ 
where, by Lemma \ref{lemma-tecnico-O}, the second term is
$O(n^{-\nu})$, while the sum in the first term has the asymptotic
behavior given in \eqref{recall-rho}.

\item[Case $\nu=\gamma$] By assumption \eqref{ass-q} \footnote{If
  there exists $j$ such that $q=c\alpha_j$, we have to consider the
  other asymptotic expression given in Lemma \ref{lemma-tecnico-G}.} and
  Lemma \ref{lemma-tecnico-G}, we have $G_{k+1,n-1}(c\alpha_j,q)=
  O\big(|F^{\gamma}_{k+1,n-1}(q)|+|F^{\gamma}_{k+1,n-1}(c\alpha_j)|\big)$,
  and so we get
$$
\rho^{21}_{n}=\sum_{j\in J,\,j\neq 1} \ind_{\{2\in I_j\}}
O\left(\sum_{k=m_0}^{n-1}
\frac{1}{k^{2\gamma}}|F^{\gamma}_{k+1,n-1}(q)|
+
\sum_{k=m_0}^{n-1}
\frac{1}{k^{2\gamma}}|F^{\gamma}_{k+1,n-1}(c\alpha_j)|
\right)\,,
$$ where, by Lemma \ref{lemma-tecnico-O}, we have for $x=q$ or
$x\in\{c\alpha_j:\, j\in J,\,j\neq 1\}$
$$
\sum_{k=m_0}^{n-1}
\frac{1}{k^{2\gamma}}|F^{\gamma}_{k+1,n-1}(x)|
=
\begin{cases}
O\left(n^{-\gamma}\right)
\qquad &\hbox{if } 1/2<\nu=\gamma<1,\\
O\left(n^{-a_x}\right)
\qquad &\hbox{if } \nu=\gamma=1\;\hbox{and } 1/2< a_x<1,\\
O\left(n^{-1}\ln(n)\right)
\qquad &\hbox{if } \nu=\gamma=1\;\hbox{and } a_x=1,\\
O\left(n^{-1}\right)
\qquad &\hbox{if } \nu=\gamma=1\;\hbox{and } a_x>1
\end{cases}
$$ 
and so, setting $x^*:=\min\{q,ca^*\}$, we can write
$$
\rho^{21}_{n}=\sum_{j\in J,\,j\neq 1} \ind_{\{2\in I_j\}}
\begin{cases}
O\left(n^{-\gamma}\right)
\qquad &\hbox{if } 1/2<\nu=\gamma<1,\\
O\left(n^{-x^*}\right)
\qquad &\hbox{if } \nu=\gamma=1\;\hbox{and } 1/2< x^*<1,\\
O\left(n^{-1}\ln(n)\right)
\qquad &\hbox{if } \nu=\gamma=1\;\hbox{and } x^*=1,\\
O\left(n^{-1}\right)
\qquad &\hbox{if } \nu=\gamma=1\;\hbox{and } x^*>1\,.
\end{cases}
$$ 
\end{description}

Summing up, taking into account the conditions $ca^*> 1/2$ when
$\gamma = 1$ and $q > 1/2$ when $\nu = 1$, from the asymptotic
behavior given above we easily obtain that in all the cases
$(i)$-$(v)$ we have
\begin{equation}\label{conclusion-rho}
t_n(J(I))\left|\sum_{k=m_0}^{n-1}\boldsymbol{\rho}^{J(I)}_{n,k}\right|
\stackrel{a.s.}\longrightarrow \mathbf{0}\,.
\end{equation}
In the case $(vi)$, the evaluation of the asymptotic behavior
given in \eqref{recall-rho} for the term $\rho_n^{22}$ is not enough
in order to conclude that $t_n(J(I))\rho_n^{22}\to 0$ a.s. Therefore, 
we need a better evaluation, that we can get applying Lemma
\ref{lemma-tecnico-O} in a different way. Indeed, in the case $(vi)$,
taking $u>1$ and applying Lemma \ref{lemma-tecnico-O} with $e=u$,
$\delta=\nu$ and $\beta=2\gamma u>1$, we find 
\begin{equation*}
\begin{split}
\left(t_n(J(I))\rho_n^{22}\right)^u&=
n^{u \nu/2}
O\left(\sum_{k=m_0}^{n-1}k^{-2\gamma u}|F_{k+1,n-1}^{\nu}(q)|^u\right)
\\
&=
n^{u \nu/2}
\begin{cases}
O\left(n^{-2\gamma u+\nu}\right)
\qquad &\hbox{if } 1/2<\nu<1,\\
O\left(n^{-q u}\right)
\qquad &\hbox{if } \nu=1\;\hbox{and } 1/2<q<2\gamma-u^{-1},\\
O\left(n^{-q u}\ln(n)\right)
\qquad &\hbox{if } \nu=1\;\hbox{and } q=2\gamma-u^{-1}>1/2,\\
O\left(n^{-2\gamma u+1}\right)
\qquad &\hbox{if } \nu=1\;\hbox{and } q>\max\{1/2,2\gamma-u^{-1}\}\,.
\end{cases}
\end{split}
\end{equation*}
Hence, from the above relations we get that it is possible to find
$u>1$ large enough such that $(t_n(J(I))\rho_n^{22})^u\to 0$ a.s, that
trivially implies $t_n(J(I))\rho_n^{22}\to 0$ a.s. Therefore also in
the case $(vi)$, we can conclude that \eqref{conclusion-rho} holds
true.  \\
\noindent{\bf Fourth step: asymptotic behavior of
  $\sum_{k=m_0}^{n-1}\mathbf{T}^{J(I)}_{k+1,n-1}$.}  \\
\noindent We aim at proving that, for each of the cases $(i)-(vi)$,
the quantity $t_n(J(I))\sum_{k=m_0}^{n-1}\mathbf{T}^{J(I)}_{k+1,n-1}$
converges stably to the desired Gaussian kernel.  For this purpose, we
apply Theorem \ref{thm:triangular}. More precisely, we set ${\mathcal
  G}_{k,n}={\mathcal F}_{k+1}$ and, given the fact that condition
$(c1)$ required in this theorem is obviously satisfied, we check only
conditions $(c2)$ and $(c3)$.  
\\
\indent For condition $(c2)$, we have to study the convergence of 
$$
t_n(J(I))^2
\sum_{k=m_0}^{n-1}\mathbf{T}^{J(I)}_{k+1,n-1}(\mathbf{T}^{J(I)}_{k+1,n-1})^{\top}.
$$
To this end, we note that 
\begin{equation*}\begin{aligned}
\sum_{k=m_0}^{n-1}\mathbf{T}^{J(I)}_{k+1,n-1}(\mathbf{T}^{J(I)}_{k+1,n-1})^{\top}&&=&
 \sum_{j_1\in J,\,j_2\in J}\sum_{i_1\in I_{j_1},\,i_2\in I_{j_2}}
U^0_{j_1(i_1)} \left(
\sum_{k=m_0}^{n-1}\mathbf{T}^{j_1}_{k+1,n-1}(\mathbf{T}^{j_2}_{k+1,n-1})^{\top}\right)
U^{0,\top}_{j_2(i_2)}\\
&&=& \sum_{j_1\in J,\,j_2\in J}\sum_{i_1\in I_{j_1},\,i_2\in I_{j_2}}
\left[
\sum_{k=m_0}^{n-1}\mathbf{T}^{j_1}_{k+1,n-1}
(\mathbf{T}^{j_2}_{k+1,n-1})^{\top}
\right]_{(i_1,i_2)}
\mathbf{u}_{j_1(i_1)}\mathbf{u}_{j_2(i_2)}^{\top}\,,
\end{aligned}
\end{equation*} 
where
\begin{equation*}
\mathbf{T}^{j}_{k+1,n-1}:=
A^j_{k+1,n-1}P_j^{-1}D_{R,k}V_j^{*\top}
\Delta\mathbf{M}_{Y,k+1}.
\end{equation*}
Thus, we can focus on the convergence of $t_n(J(I))^2
\sum_{k=m_0}^{n-1}\mathbf{T}^{j_1}_{k+1,n-1}(\mathbf{T}^{j_2}_{k+1,n-1})^{\top}$. 
Regarding to this, we observe that
\begin{equation*}
\mathbf{T}^{j_1}_{k+1,n-1}(\mathbf{T}^{j_2}_{k+1,n-1})^{\top}
=
A^{j_1}_{k+1,n-1}H^{j_1,j_2}_{k+1}(A^{j_2}_{k+1,n-1})^{\top},
\end{equation*}
where
\begin{equation*}
\begin{aligned}
H^{j_1,j_2}_{k+1}&&:=&\ P_{j_1}^{-1}D_{R,k}V_{j_1}^{*\top}
\Delta\mathbf{M}_{Y,k+1}
\Delta\mathbf{M}_{Y,k+1}^{\top}V_{j_2}^{*}D_{R,k}P_{j_2}^{-\top}\\
&&=&\ P_{j_1}^{-1}D_{R,k}V_{j_1}^{*\top}\begin{pmatrix}
I \\
I
\end{pmatrix} \Delta\mathbf{M}_{k+1}
\Delta\mathbf{M}_{k+1}^{\top}\begin{pmatrix} I & I
\end{pmatrix}V_{j_2}^{*}D_{R,k}P_{j_2}^{-\top}\\
&&=&\ P_{j_1}^{-1}D_{R,k}\mathbf{1}
\mathbf{v}_{j_1}^{\top}\Delta\mathbf{M}_{k+1}
\Delta\mathbf{M}_{k+1}^{\top}\mathbf{v}_{j_2}\mathbf{1}^{\top}D_{R,k}P_{j_2}^{-\top}
\\
&&=&\ \mathbf{h}^{j_1}_k
\mathbf{v}_{j_1}^{\top}\Delta\mathbf{M}_{k+1}
\Delta\mathbf{M}_{k+1}^{\top}\mathbf{v}_{j_2}(\mathbf{h}_{k}^{j_2})^{\top}\\
&&=&\ \beta_{k+1}^{j_1,j_2}\mathbf{h}^{j_1}_k
(\mathbf{h}_{k}^{j_2})^{\top},
\end{aligned},
\end{equation*}
with 
\begin{equation*}
\beta_{k+1}^{j_1,j_2}\ :=\
\mathbf{v}_{j_1}^{\top}\Delta\mathbf{M}_{k+1}
\Delta\mathbf{M}_{k+1}^{\top}\mathbf{v}_{j_2}\qquad\mbox{and}
\end{equation*}
\begin{equation*}
\mathbf{h}_{k}^j:= P_{j}^{-1}D_{R,k}\mathbf{1}=
\begin{pmatrix}
\widehat{r}_{k-1} \\
\widehat{q}_{k,k}-\widehat{r}_{k-1}g(\lambda_{j})
\end{pmatrix}.
\end{equation*}
Now, we set $\mathbf{d}^{j}_{k,n}:=A^{j}_{k+1,n-1}\mathbf{h}^{j}_{k}$,
so that we can write
\begin{equation}\label{eq-proof-new}
\sum_{k=m_0}^{n-1}\mathbf{T}^{j_1}_{k+1,n-1}(\mathbf{T}^{j_2}_{k+1,n-1})^{\top}
\ =\  
\sum_{k=m_0}^{n-1}\beta_{k+1}^{j_1,j_2}\mathbf{d}_{k,n}^{j_1}
(\mathbf{d}_{k,n}^{j_2})^{\top}\,.
\end{equation}
Hence, in order to obtain the almost sure convergence of $t_n(J(I))^2
\sum_{k=m_0}^{n-1}\mathbf{T}^{j_1}_{k+1,n-1}(\mathbf{T}^{j_2}_{k+1,n-1})^{\top}$,
by means of~the usual martingale arguments (see \cite[Lemma
  B.1]{ale-cri-ghi-MEAN}) and the technical results collected in
Section \ref{lemmi-tecnici-appendice} of the appendix, it is enough to
prove the convergence of
$t_n(J(I))^2\sum_{k=m_0}^{n-1}\mathbf{d}_{k,n}^{j_1}
(\mathbf{d}_{k,n}^{j_2})^{\top}.$ Indeed, since $\{X_{n,j}:\,
j=1,\dots, N\}$ are conditionally independent given $\mathcal{F}_n$,
we have
\begin{equation*}\label{eq:covar-conditional-M}
E[\Delta{M}_{n,h}\Delta{M}_{n,j}\,|\,{\mathcal F}_{n-1}]=0\quad
\mbox{ for } h\neq j;
\end{equation*}
while, for each $j$, using the normalization $W^{\top}{\mathbf
  1}={\mathbf 1}$, we have
\begin{equation*}\label{eq:var-conditional-M}
E[(\Delta{M}_{n,j})^2\,|\,{\mathcal F}_{n-1}]= \left(\sum_{h=1}^N
w_{h,j}Z_{n-1,h}\right) \left(1-\sum_{h=1}^N
w_{h,j}Z_{n-1,h}\right) \stackrel{a.s.}\longrightarrow
Z_{\infty}(1-Z_{\infty}).
\end{equation*}
Therefore, we get   
\begin{equation*}\label{eq:multidim-limite-conditional-M}
E[(\Delta{\mathbf M}_{n})(\Delta{\mathbf M}_{n})^{\top}
\,|\,{\mathcal F}_{n-1}] \stackrel{a.s.}\longrightarrow
Z_{\infty}(1-Z_{\infty})I
\end{equation*}
and so 
\begin{equation*}
E[\beta_{n+1}^{j_1,j_2}\, |\, {\mathcal F}_{n}] =
\mathbf{v}_{j_1}^{\top}
E[\Delta{\mathbf M}_{n+1} (\Delta{\mathbf M}_{n+1})^{\top}\,|\,{\mathcal F}_{n}]
\mathbf{v}_{j_2}
\ \stackrel{a.s}\longrightarrow\ 
 Z_{\infty}(1-Z_{\infty})\mathbf{v}_{j_1}^{\top}\mathbf{v}_{j_2}\,,
\end{equation*}
from which we finally obtain 
\begin{equation*}
a.s.-\lim_n t_n(J(I))^2
\sum_{k=m_0}^{n-1}\mathbf{T}^{j_1}_{k+1,n-1}(\mathbf{T}^{j_2}_{k+1,n-1})^{\top}
\ =\ Z_{\infty}(1-Z_{\infty})\mathbf{v}_{j_1}^{\top}\mathbf{v}_{j_2}
\lim_n
t_n(J(I))^2\sum_{k=m_0}^{n-1}\mathbf{d}_{k,n}^{j_1}(\mathbf{d}_{k,n}^{j_2})^{\top}.
\end{equation*}
In order to compute the limits in the last term of the above relation,
we observe that, by means of~\eqref{eq:A_1_k_n}
and~\eqref{eq:A_j_k_n_withG}, we have the following analytic
expression of $\mathbf{d}_{k,n}^j$:
\begin{equation}\label{componenti-d-1}
\mathbf{d}_{k,n}^1\ =\ A^{1}_{k+1,n-1}\mathbf{h}_{1,k}\ =\
\begin{pmatrix}
\widehat{r}_{k-1}\\
\left(\widehat{q}_{k,k}-\widehat{r}_{k-1}\right)F^{\nu}_{k+1,n-1}(q)
\end{pmatrix}
\end{equation}
and, for $j\geq 2$,  
\begin{equation}\label{componenti-d}
\mathbf{d}_{k,n}^j\ =\ A^{j}_{k+1,n-1}\mathbf{h}_{j,k}\ =\
\begin{pmatrix}
\widehat{r}_{k-1}F^{\gamma}_{k+1,n-1}(c\alpha_{j}) \\
\lambda_{j}\widehat{r}_{k-1}G_{k+1,n-1}(c\alpha_j,q)+
\left(\widehat{q}_{k,k}-\widehat{r}_{k-1}g(\lambda_{j})\right)F^{\nu}_{k+1,n-1}(q)
\end{pmatrix}.
\end{equation}
Using these equalities, in Section \ref{result-limit-d} of the
appendix, for all the considered cases $(i)-(vi)$, we find the limit
of each component of $t_n(J(I))^2
\sum_{k=m_0}^{n-1}\mathbf{d}_{k,n}^{j_1}(\mathbf{d}_{k,n}^{j_2})^{\top}$,
that is we compute
\begin{equation*}
d^{j_1(i_1),j_2(i_2)}\ :=\
\lim_{n}t_n(J(I))^2\sum_{k=m_0}^{n-1}d_{k,n}^{j_1(i_1)}d_{k,n}^{j_2(i_2)}\,,
\end{equation*}
where $d_{k,n}^{j(1)}$ and $d_{k,n}^{j(2)}$ are, respectively, the
first and the second component of $\mathbf{d}_{k,n}^j$ given
in~\eqref{componenti-d-1} and~\eqref{componenti-d}.  Summing up, we
have
$$
\sum_{k=m_0}^{n-1}\mathbf{T}^{J(I)}_{k+1,n-1}(\mathbf{T}^{J(I)}_{k+1,n-1})^{\top}
\ \stackrel{a.s.}\longrightarrow\ 
Z_{\infty}(1-Z_{\infty})
\sum_{j_1\in J,\,j_2\in J}
\mathbf{v}_{j_1}^{\top}\mathbf{v}_{j_2}
\sum_{i_1\in I_{j_1},\,i_2\in I_{j_2}}
d^{j_1(i_1),j_2(i_2)}
\mathbf{u}_{j_1(i_1)}\mathbf{u}_{j_2(i_2)}^{\top}\,.
$$

\indent For the check of condition (c3) of Theorem
\ref{thm:triangular}, we observe that, by~\eqref{eq:A_1_k_n},
\eqref{eq:A_j_k_n_withG}, \eqref{eq:G} and~\eqref{U0}, taking into
account the fact that in all the considered cases with $1\in J$,
i.e.~$(ii), (iii)$ and $(v)$, we have $1\notin I_1$, we can write
$$
|{\mathbf T}^{J(I)}_{k+1,n-1}|=
O(\Gamma_{k+1,n-1}^{11})+O(\Gamma_{k+1,n-1}^{21})+O(\Gamma_{k+1,n-1}^{22}),
$$ where $\Gamma_{k+1,n-1}^{11},\,\Gamma_{k+1,n-1}^{21}$ and
$\Gamma_{k+1,n-1}^{22}$ are the following deterministic
quantities:
\begin{equation*}
\begin{aligned}
\Gamma_{k+1,n-1}^{11}\ &&:=&\ \sum_{j\in J,\, j\neq 1}\ind_{\{1\in
I_j\}}\widehat{r}_{k-1}|F_{k+1,n-1}^{\gamma}(c\alpha_j)|\,,
\\
\Gamma_{k+1,n-1}^{21}\ &&:=&\ \sum_{j\in J,\, j\neq 1} \ind_{\{2\in I_j\}}
\widehat{r}_{k-1}|G_{k+1,n-1}(c\alpha_j,q)|
\\
\Gamma_{k+1,n-1}^{22}\ &&:=&\ \sum_{j\in J} \ind_{\{2\in I_j\}}
(\widehat{r}_{k-1}+\widehat{q}_{k,k})|F_{k+1,n-1}^{\nu}(q)|\,.
\end{aligned}
\end{equation*}
Therefore, we find for any $u>1$
\begin{equation*}
\begin{split}
&\big(\sup_{m_0\leq k\leq n-1}|t_n(J(I))\mathbf{T}^{J(I)}_{k+1,n-1}|\big)^{2u}
\leq t_n(J(I))^{2u} \sum_{k=m_0}^{n-1} |\mathbf{T}^{J(I)}_{k+1,n-1}|^{2u}
=\\
&
t_n(J(I))^{2u}
\left\{ 
\sum_{k=m_0}^{n-1} O\left((\Gamma^{11}_{k+1,n-1})^{2u}\right)
+
\sum_{k=m_0}^{n-1} O\left((\Gamma^{21}_{k+1,n-1})^{2u}\right)
+
\sum_{k=m_0}^{n-1} O\left((\Gamma^{22}_{k+1,n-1})^{2u}\right)
\right\}\,.
\end{split}
\end{equation*}
We now analyze the last three terms.  For the first one, by Lemma
\ref{lemma-tecnico-O} with $\beta=2\gamma u$, $e=2u$ and
$\delta=\gamma$, we have
\begin{equation*}
\begin{split} 
\sum_{k=m_0}^{n-1} O\left((\Gamma^{11}_{k+1,n-1})^{2u}\right)
&=\sum_{j\in J,j\neq 1} I_{\{1\in I_j\}} 
O\left( 
\sum_{k=m_0}^{n-1}\frac{1}{k^{2\gamma u}}|F^\gamma_{k+1,n-1}(c\alpha_j)|^{2u}
\right)\\
&=
\sum_{j\in J,j\neq 1} I_{\{1\in I_j\}} 
\begin{cases}
O\left(n^{-\gamma(2u-1)}\right)
\qquad &\hbox{if } 1/2<\gamma<1,\\
O\left(n^{-2ca^*u}\right)
\qquad &\hbox{if } \gamma=1\;\hbox{and } 1/2<ca^*<1-(2u)^{-1},\\
O\left(n^{-2u+1}\ln(n)\right)
\qquad &\hbox{if } \gamma=1\;\hbox{and } ca^*=1-(2u)^{-1},\\
O\left(n^{-2u+1}\right)
\qquad &\hbox{if } \gamma=1\;\hbox{and } ca^*>1-(2u)^{-1}.
\end{cases}
\end{split}
\end{equation*}
For the third term, we observe that
$\widehat{r}_{k-1}=O(\widehat{q}_{k,k})$ when $\nu\leq \gamma$ and
$\widehat{q}_{k,k}=O(\widehat{r}_{k-1})$ when $\nu>\gamma$. Hence, by
Lemma \ref{lemma-tecnico-O} with $e=2u$, $\delta=\nu$ and $\beta=2\nu
u$ if $\nu\leq\gamma$ and $\beta=2\gamma u$ if $\nu>\gamma$, we get
for the case $\nu\leq\gamma$
\begin{equation*}
\begin{split}
\sum_{k=m_0}^{n-1} O\left((\Gamma^{22}_{k+1,n-1})^{2u}\right)
&=\sum_{j\in J} I_{\{2\in I_j\}} 
O\left( 
\sum_{k=m_0}^{n-1}\frac{1}{k^{2\nu u}}|F^\nu_{k+1,n-1}(q)|^{2u}
\right)\\
&=
\sum_{j\in J} I_{\{2\in I_j\}} 
\begin{cases}
O\left(n^{-\nu(2u-1)}\right)
\qquad &\hbox{if } 1/2<\nu<1,\\
O\left(n^{-2qu}\right)
\qquad &\hbox{if } \nu=1\;\hbox{and } 1/2<q<1-(2u)^{-1},\\
O\left(n^{-2u+1}\ln(n)\right)
\qquad &\hbox{if } \nu=1\;\hbox{and } q=1-(2u)^{-1},\\
O\left(n^{-2u+1}\right)
\qquad &\hbox{if } \nu=1\;\hbox{and } q>1-(2u)^{-1}\,,
\end{cases}
\end{split}
\end{equation*}
and for the case $\nu>\gamma$
\begin{equation}\label{recall-Gamma}
\begin{split}
\sum_{k=m_0}^{n-1} O\left((\Gamma^{22}_{k+1,n-1})^{2u}\right)
&=\sum_{j\in J} I_{\{2\in I_j\}} 
O\left( 
\sum_{k=m_0}^{n-1}\frac{1}{k^{2\gamma u}}|F^\nu_{k+1,n-1}(q)|^{2u}
\right)\\
&=
\sum_{j\in J} I_{\{2\in I_j\}} 
\begin{cases}
O\left(n^{-2\gamma u+\nu}\right)
\qquad &\hbox{if } 1/2<\nu<1,\\
O\left(n^{-2qu}\right)
\qquad &\hbox{if } \nu=1\;\hbox{and } 1/2<q<\gamma-(2u)^{-1},\\
O\left(n^{-2\gamma u+1}\ln(n)\right)
\qquad &\hbox{if } \nu=1\;\hbox{and } q=\gamma-(2u)^{-1}>1/2,\\
O\left(n^{-2\gamma u+1}\right)
\qquad &\hbox{if } \nu=1\;\hbox{and } q>\max\{1/2,\gamma-(2u)^{-1}\}\,.
\end{cases}
\end{split}
\end{equation}
For the second term, we apply Lemma \ref{lemma-tecnico-O} together
with Lemma \ref{lemma-tecnico-G} so that we get:
\begin{description}
\item[Case $\nu<\gamma$] We have
  $G_{k+1,n-1}(c\alpha_j,q)=O\big(n^{-(\gamma-\nu)}|F^{\nu}_{k+1,n-1}(q)|+
  k^{-(\gamma-\nu)}|F^{\gamma}_{k+1,n-1}(c\alpha_j)|\big)$ by means of
  Lemma \ref{lemma-tecnico-G}, and so we find
\begin{equation*}
\begin{split}
&\sum_{k=m_0}^{n-1} O\left((\Gamma^{21}_{k+1,n-1})^{2u}\right)
=\sum_{j\in J,j\neq 1} I_{\{2\in I_j\}} 
O\left(\sum_{k=m_0}^{n-1}\frac{1}{k^{2\gamma u}}|G_{k+1,n-1}(c\alpha_j,q)|^{2u}\right)
=
\\
&\sum_{j\in J,j\neq 1} I_{\{2\in I_j\}} 
O\left( 
n^{-2(\gamma-\nu)u}
\sum_{k=m_0}^{n-1}\frac{1}{k^{2\gamma u}}|F^\nu_{k+1,n-1}(q)|^{2u}
+
\sum_{k=m_0}^{n-1}\frac{1}{k^{4\gamma u-2\nu u}}|F^\gamma_{k+1,n-1}(c\alpha_j)|^{2u}
\right)\,,
\end{split}
\end{equation*}
where, by Lemma \ref{lemma-tecnico-O}, the first term is 
$O(n^{-4\gamma u +2\nu u+\nu})$, while for the second term we have 
$$
\sum_{k=m_0}^{n-1}\frac{1}{k^{4\gamma u-2\nu u}}|F^\gamma_{k+1,n-1}(c\alpha_j)|^{2u}
=
\begin{cases}
O\left(n^{-4\gamma u+2\nu u+\gamma}\right)
\qquad &\hbox{if } 1/2<\gamma<1,\\
O\left(n^{-2ca^*u}\right)
\qquad &\hbox{if } \gamma=1\;\hbox{and } 1/2< ca^*<2-\nu-(2u)^{-1},\\
O\left(n^{-2ca^*u}\ln(n)\right)
\qquad &\hbox{if } \gamma=1\;\hbox{and } ca^*=2-\nu-(2u)^{-1},\\
O\left(n^{-4u+2\nu u+1}\right)
\qquad &\hbox{if } \gamma=1\;\hbox{and } ca^*>2-\nu-(2u)^{-1}\,.
\end{cases}
$$

\item[Case $\nu>\gamma$] We have
  $G_{k+1,n-1}(c\alpha_j,q)=O\big(n^{-(\nu-\gamma)}|F^{\nu}_{k+1,n-1}(q)|+
  k^{-(\nu-\gamma)}|F^{\gamma}_{k+1,n-1}(c\alpha_j)|\big)$ by means of
  Lemma \ref{lemma-tecnico-G}, and so we find
\begin{equation*}
\begin{split}
&\sum_{k=m_0}^{n-1} O\left((\Gamma^{21}_{k+1,n-1})^{2u}\right)
=
\sum_{j\in J,j\neq 1} I_{\{2\in I_j\}} 
O\left(\sum_{k=m_0}^{n-1}\frac{1}{k^{2\gamma u}}|G_{k+1,n-1}(c\alpha_j,q)|^{2u}\right)
=
\\
&\sum_{j\in J,j\neq 1} I_{\{2\in I_j\}} 
O\left( 
n^{-2(\nu-\gamma)u}
\sum_{k=m_0}^{n-1}\frac{1}{k^{2\gamma u}}|F^\nu_{k+1,n-1}(q)|^{2u}
+
\sum_{k=m_0}^{n-1}\frac{1}{k^{2\nu u}}|F^\gamma_{k+1,n-1}(c\alpha_j)|^{2u}
\right)\,,
\end{split}
\end{equation*}
where, by Lemma \ref{lemma-tecnico-O}, the second term is
$O(n^{-2\nu u+\gamma})$, while the sum in the first term has the
asymptotic behavior given in \eqref{recall-Gamma}.

\item[Case $\nu=\gamma$] By assumption \eqref{ass-q} \footnote{If
  there exists $j$ such that $q=c\alpha_j$, we have to consider the
  other asymptotic expression given in Lemma \ref{lemma-tecnico-G}.} and
  Lemma \ref{lemma-tecnico-G}, we have
  $G_{k+1,n-1}(c\alpha_j,q)=O\big(|F^{\gamma}_{k+1,n-1}(q)|+
  |F^{\gamma}_{k+1,n-1}(c\alpha_j)|\big)$, and so we find
$$
\sum_{k=m_0}^{n-1} O\left((\Gamma^{21}_{k+1,n-1})^{2u}\right)
=\sum_{j\in J,j\neq 1} I_{\{2\in I_j\}} 
O\left( 
\sum_{k=m_0}^{n-1}\frac{1}{k^{2\gamma u}}|F^\gamma_{k+1,n-1}(q)|^{2u}
+
\sum_{k=m_0}^{n-1}\frac{1}{k^{2\gamma u}}|F^\gamma_{k+1,n-1}(c\alpha_j)|^{2u}
\right)\,,
$$ 
where, by Lemma \ref{lemma-tecnico-O}, we have for $x=q$ or
$x\in\{c\alpha_j:\, j\in J,\,j\neq 1\}$
$$
\sum_{k=m_0}^{n-1}
\frac{1}{k^{2\gamma u}}|F^{\gamma}_{k+1,n-1}(x)|^{2u}
=
\begin{cases}
O\left(n^{-\gamma(2u-1)}\right)
\qquad &\hbox{if } 1/2<\nu=\gamma<1,\\
O\left(n^{-2a_xu}\right)
\qquad &\hbox{if } \nu=\gamma=1\;\hbox{and } 1/2< a_x<1-(2u)^{-1},\\
O\left(n^{-2u+1}\ln(n)\right)
\qquad &\hbox{if } \nu=\gamma=1\;\hbox{and } a_x=1-(2u)^{-1},\\
O\left(n^{-2u+1}\right)
\qquad &\hbox{if } \nu=\gamma=1\;\hbox{and } a_x>1-(2u)^{-1}
\end{cases}
$$ 
and so, setting $x^*:=\min\{q,ca^*\}$, we can write
$$
\sum_{k=m_0}^{n-1} O\left((\Gamma^{21}_{k+1,n-1})^{2u}\right)
=\sum_{j\in J,j\neq 1} I_{\{2\in I_j\}} 
\begin{cases}
O\left(n^{-\gamma(2u-1)}\right)
\qquad &\hbox{if } 1/2<\nu=\gamma<1,\\
O\left(n^{-2x^*u}\right)
\qquad &\hbox{if } \nu=\gamma=1\;\hbox{and } 1/2< x^*<1-(2u)^{-1},\\
O\left(n^{-2u+1}\ln(n)\right)
\qquad &\hbox{if } \nu=\gamma=1\;\hbox{and } x^*=1-(2u)^{-1},\\
O\left(n^{-2u+1}\right)
\qquad &\hbox{if } \nu=\gamma=1\;\hbox{and } x^*>1-(2u)^{-1}\,.
\end{cases}
$$ 
\end{description}

Summing up, taking into account the conditions $ca^* >1/2$ when
$\gamma = 1$ and $q > 1/2$ when $\nu = 1$, we can conclude that in all
the six cases $(i)$-$(vi)$, there exists a suitable $u>1$ such that
$$
\big(\sup_{m_0\leq k\leq n-1}|t_n(J(I))\mathbf{T}^{J(I)}_{k+1,n-1}|\big)^{2u}
\stackrel{L^1}\longrightarrow 0.
$$ 
This convergence trivially implies condition $(c3)$ of Theorem
\ref{thm:triangular}.
\qed

\subsection{Proof of Theorem \ref{clt-Y}}
The proof of Theorem \ref{clt-Y} follows by recalling
that
$$
(\widetilde{\mathbf Y}_n-Z_\infty\mathbf{1})\ =\
(\widetilde{\mathbf Y}_n-Z_\infty\mathbf{1})\ +\ \widehat{\mathbf
Y}_n\,,
$$ where the convergence rate for the first term is
$n^{\gamma-\frac{1}{2}}$ for any parameters (see \eqref{clt-tilde-Y}),
while the convergence rate of the second term is $n^e$, with $e$
specified in Theorem \ref{clt-hat-Y} according to the values of the
parameters.  Therefore, we can have three different cases:
\begin{itemize}
\item If $e<\gamma-\frac{1}{2}$, then we have
$$
n^e(\mathbf{Y}_n-Z_\infty\mathbf{1})=
\frac{n^e}{n^{\gamma-\frac{1}{2}}}
n^{\gamma-\frac{1}{2}}(\widetilde{\mathbf
Y}_n-Z_\infty\mathbf{1})+ n^e\widehat{\mathbf Y}_n\,,
$$
where the first term converges in probability to zero and the
second term converges stably to a certain Gaussian kernel. This
occurs only in case (a) with $e=\nu/2$ and $\nu<\gamma_0$.

\item If $e>\gamma-\frac{1}{2}$, then we have
$$
n^{\gamma-\frac{1}{2}}(\mathbf{Y}_n-Z_\infty\mathbf{1})=
n^{\gamma-\frac{1}{2}}(\widetilde{\mathbf
Y}_n-Z_\infty\mathbf{1})+ \frac{n^{\gamma-\frac{1}{2}}}{n^e}
n^e\widehat{\mathbf Y}_n\,,
$$
where the first term converges stably (in the strong sense) to the
Gaussian kernel given in \eqref{clt-tilde-Y} and the second term
converges in probability to zero. This occurs in case (a) with
$e=\nu/2$ and $\gamma_0<\nu<\gamma$, in case (b) with $e=\gamma/2$
and $\nu=\gamma<1$ and in case (c) with $e=\gamma-\nu/2$ and
$\gamma<\nu<1$.

\item  If $e=\gamma-\frac{1}{2}$, then we have
$$
n^{\gamma-\frac{1}{2}}(\mathbf{Y}_n-Z_\infty\mathbf{1})=
n^{\gamma-\frac{1}{2}}(\widetilde{\mathbf
Y}_n-Z_\infty\mathbf{1})+ n^{\gamma-\frac{1}{2}}\widehat{\mathbf
Y}_n\,,
$$
where the first term converges stably in the strong sense to the
Gaussian kernel given in \eqref{clt-tilde-Y} and the second term is
${\mathcal F}_n$-measurable and it converges stably to a certain
Gaussian kernel. Thus, in this case, we can apply Theorem \ref{blocco}
in Appendix. This occurs in case (a) with $e=\nu/2$ and
$\nu=\gamma_0<1$, in case (b) with $e=\gamma/2$ and $\nu=\gamma=1$
(i.e.~$\nu=\gamma_0=1$) and in case (c) with $e=\gamma-\nu/2$ and
$\gamma<\nu=1$ (i.e.~$\gamma_0<\nu=1$).
\end{itemize}
\qed

\begin{rem}\label{rem-N=1-bis}
\rm As told in Remark \ref{rem-N=1}, statements (a), (b) and (c) of
Theorem \ref{clt-Y} with $N=1$ (and so without the condition on
$\lambda^*$) can be proven with the same proof. Specifically, it is
enough to take into account that when $N=1$, we have
$\widehat{Y}_n=Y_{1(2)}$ and ${\widetilde Z}_n=Z_n$.
\end{rem}

\section{Proof of the results for statistical applications}
\label{section_proofs_statistics}
Here we prove the convergence results stated in Section
\ref{section_statistics}. As we will see, the decomposition of
$\mathbf{Y}_n$ given in Section \ref{section_decomposition} is a
fundamental tool also for the proof of these results.

\subsection{Proof of Theorem \ref{thm:N_tilde}}
For the proof of this result, we need the following lemma:

\begin{lem}\label{lemma-1(2)}
Let us set 
\begin{equation}\label{def-beta}
\beta:=\frac{\nu}{2}\ind_{\{\nu\leq\gamma\}}+
\left(\gamma-\frac{\nu}{2}\right)\ind_{\{\gamma<\nu\}}.
\end{equation}
Then, under all the assumptions stated in Section \ref{section_model}, we
have
$$
n^{\beta} \mathbf{Y}_{1(2)} \stackrel{a.s.}\longrightarrow 
\mathcal{N}\left(\mathbf{0}, 
Z_\infty(1-Z_{\infty})\frac{\|\mathbf{v}_1\|^2}{N} 
d^{1(2),1(2)}
\begin{pmatrix}
\mathbf{0} & \mathbf{0}\\
\mathbf{0} & \mathbf{1}\mathbf{1}^{\top}
\end{pmatrix}\,,
\right)
$$
where
$$
d^{1(2),1(2)}=
\begin{cases}
\frac{q}{2}\quad&\hbox{for } \nu<\gamma,\\
\frac{(q-c)^2}{2q-\ind_{\{\nu=1\}}}\quad&\hbox{for } \nu=\gamma,\\
\frac{c^2}{2q-\ind_{\{\nu=1\}}(2\gamma-1)}\quad&\hbox{for } \gamma<\nu.
\end{cases}
$$
\end{lem}

\noindent{\bf Proof.} We observe that $\mathbf{Y}_{1(2)}$ can be written as 
the general sum \eqref{general-sum} with $J=\{1\}$ and
$I_1=\{2\}$. Therefore case $\nu>\gamma$ coincides with the case $(v)$
of Lemma \ref{conv-general-sum}, taking into account the value
$d^{1(2),1(2)}$ computed in Section \ref{result-limit-d} for this case
and equality \eqref{equality-matrix}. The cases $\nu<\gamma$ and
$\nu=\gamma$ follows from the same arguments employed for the proof of
Lemma \ref{conv-general-sum}, setting $t_n(J(I))=n^{\nu/2}$ and using
the value $d^{1(2),1(2)}$ obtained in Section \ref{result-limit-d}
when $\nu \leq \gamma$.
\qed

\begin{rem}\label{remark-d12}
\rm Note that, when $\nu=\gamma$ and $q=c$, we have $d^{1(2),1(2)}=0$
and so we obtain that $n^{\beta}\mathbf{Y}_{1(2)}$ converges to
$\mathbf{0}$ in probability.  This means that in this case the
convergence of $\mathbf{Y}_{1(2)}$ to $\mathbf{0}$ is faster than
$n^{-\beta}=n^{-\gamma/2}$.
\end{rem}

\noindent{\bf Proof of Theorem \ref{thm:N_tilde}.} The convergence
rate and the second-order asymptotic distribution of $\widetilde{N}_n$
can be obtained by combining the second-order convergences of the two
stochastic processes $\widetilde{Z}_n$ and
$(\widetilde{N}_n-\widetilde{Z}_n)$. In order to get the convergence
results for these two last processes, we observe that
\begin{equation*}
\begin{split}
&N^{-1/2}\mathbf{u}_1^{\top}
\begin{pmatrix} \mathbf{0} & I
\end{pmatrix}
\widetilde{\mathbf{Y}}_n\ =\
\widetilde{Z}_nN^{-1/2}\mathbf{u}_1^{\top}\mathbf{1}
=\widetilde{Z}_n
\qquad\mbox{and}
\\
&N^{-1/2}\mathbf{u}_1^{\top}
\begin{pmatrix} \mathbf{0} & I
\end{pmatrix}
\mathbf{Y}_{1(2),n}\ =\ 
N^{-1/2}\mathbf{u}_1^{\top}
\begin{pmatrix} \mathbf{0} & I
\end{pmatrix}
\mathbf{u}_{1(2)}\mathbf{v}_{1(2)}^{\top}\mathbf{Y}_n\ =\
N^{-1/2}\mathbf{u}_1^{\top}
\begin{pmatrix}
-\mathbf{u}_1\mathbf{v}_1^{\top} & \mathbf{u}_1\mathbf{v}_1^{\top}
\end{pmatrix}
\mathbf{Y}_n
\ =\ 
\\
&(\widetilde{N}_n-\widetilde{Z}_n)N^{-1/2}\mathbf{u}_1^{\top}\mathbf{1}
=\widetilde{N}_n-\widetilde{Z}_n
\end{split}
\end{equation*}
(where we have used \eqref{eq-def-Ytilde} for the first equality and
relations \eqref{eq-NUM2}, \eqref{eq-uji}, \eqref{eq-vji},
\eqref{eq-NUM4} and \eqref{eq:relazioni-1} for the other
equalities). Hence, from the convergence result stated
in~\eqref{clt-tilde-Y} and Lemma \ref{lemma-1(2)}, together with
Remark \ref{remark-d12}, we obtain that $\widetilde{Z}_n$ converges in
probability to the random variable $Z_{\infty}$ with rate
$n^{\gamma-1/2}$ and $(\widetilde{N}_n-\widetilde{Z}_n)$ converges in
probability to zero with at least rate $n^\beta$ defined in
\eqref{def-beta}. Then, since
$\widetilde{N}_n=\widetilde{Z}_n+(\widetilde{N}_n-\widetilde{Z}_n)$,
it is possible to follow analogous arguments to those used in the
proof of Theorem~\ref{clt-Y} to combine the asymptotic behaviors of
$\widetilde{Z}_n$ and $(\widetilde{N}_n-\widetilde{Z}_n)$. More
precisely:
\begin{itemize}
\item[(a)] in the case $\nu<\gamma_0$, we necessarily have
  $\gamma_0=2\gamma-1\leq \gamma$ (since $\gamma\leq 1$) and so we have
  $\beta=\nu/2<(\gamma-1/2)$. Thus $\widetilde{N}_n$ has the same
  convergence rate and the same asymptotic variance as
  $(\widetilde{N}_n-\widetilde{Z}_n)=
N^{-1/2}\mathbf{u}_1^{\top}
\begin{pmatrix}
\mathbf{0} & I
\end{pmatrix}
\mathbf{Y}_{1(2),n}$, that is (see Lemma \ref{lemma-1(2)}) we get
\begin{equation*}
n^{\nu/2}(\widetilde{N}_n-Z_{\infty})\ \longrightarrow\
\mathcal{N} \left(\ 0\ ,\
Z_{\infty}(1-Z_{\infty})\widetilde{\sigma}^2 \
\right)\quad\hbox{stably}
\end{equation*}
with $\widetilde{\sigma}^2 =q/2$;
\item[(b)] in the case $\gamma_0<\nu<1$, we have
  $\beta>(\gamma-1/2)$ and hence $\widetilde{N}_n$ has the same
  asymptotic behavior as $\widetilde{Z}_n=
N^{-1/2}\mathbf{u}_1^{\top}
\begin{pmatrix} \mathbf{0} & I
\end{pmatrix}
\widetilde{\mathbf{Y}}_n$, that is (see~\eqref{clt-tilde-Y})
\begin{equation*}
n^{\gamma-\frac{1}{2}}(\widetilde{N}_n-Z_{\infty})\ {\longrightarrow}
\mathcal{N}\left(\ 0\ ,\
Z_{\infty}(1-Z_{\infty})\widetilde{\sigma}_{\gamma}^2 \ \right)
\qquad\hbox{stably};
\end{equation*}

\item[(c)] If $\nu=\gamma_0$ (i.e. $\nu=2\gamma-1\leq\gamma$) or
  $\nu=1$, we have $\beta=(\gamma-1/2)$ and hence the asymptotic
  behavior of $\widetilde{N}_n$ follows by combining the convergence
  results for $(\widetilde{N}_n-\widetilde{Z}_n)$ and $\widetilde{Z}_n$
  as done in the proof of Theorem~\ref{clt-Y}, and so we get
\begin{equation*}
n^{\gamma-\frac{1}{2}}(\widetilde{N}_n-Z_{\infty})\
{\longrightarrow}\ \mathcal{N} \left(\ 0\ ,\
Z_{\infty}(1-Z_{\infty})\left(\widetilde{\sigma}_{\gamma}^2
+\widetilde{\sigma}^2\right) \ \right)\quad\hbox{stably}\,,
\end{equation*}
where $\widetilde{\sigma}^2$ is defined in \eqref{eq:def-sigma-tilde}. 
\end{itemize}
\qed

\begin{rem}\label{remark-ottimo-2}
\rm Returning to Remark \ref{remark-ottimo-1}, we observe that in the
proof of Theorem \ref{thm:N_tilde} the asymptotic behavior of
$\widetilde{N}_n$ is obtained as the combination of the asymptotic
behaviors of $\widetilde{N}_n-\widetilde{Z}_n$ and
$\widetilde{Z}_n$. In case (b), $\widetilde{Z}_n$ converges slower
than $\widetilde{N}_n-\widetilde{Z}_n$, and so only the rate and the
asymptotic variance of $\widetilde{Z}_n$ appear in the statement of
the result. However, if we look at an higher level of approximation,
we should also consider the process $\widetilde{N}_n-\widetilde{Z}_n$,
that converges to zero with at least rate $n^\beta$. Then, we can note
that $\beta$ as a function of $\nu$ has its maximum in $\nu=\gamma$,
which hence provides the ``optimal value'' of $\nu$. In addition, in
this case the quantity $d^{1(2)1(2)}$ as a function of $q$ has its
minimum in $q=c$, which hence gives the ``optimal value'' of $q$. Note
that, as told in the previous Remark \ref{remark-d12}, when
$\nu=\gamma$ and $q=c$, we have $n^{\beta}\mathbf{Y}_{1(2)}\to
\mathbf{0}$ in probability and so also
$n^{\beta}(\widetilde{N}_n-\widetilde{Z}_n)\to 0$ in probability. This
means that in this case the convergence of
$\widetilde{N}_n-\widetilde{Z}_n$ to zero is faster then
$n^{-\beta}=n^{-\gamma/2}$.
\end{rem}

\subsection{Proof of Theorem \ref{thm:N_apice}}
Recalling \eqref{eq-NUM4}, together with \eqref{eq-NUM1} and the fact that
$$
U_j^*V_j^{*\top}=\begin{pmatrix}
\mathbf{u}_j\mathbf{v}_j^{\top} &\mathbf{0}\\
\mathbf{0} & \mathbf{u}_j\mathbf{v}_j^{\top}
\end{pmatrix},
$$
we can write
${\mathbf{N}}'_{n}=
\sum_{j=2}^N\mathbf{u}_{j}\mathbf{v}_{j}^{\top}\mathbf{N}_{n}=
\begin{pmatrix}
\mathbf{0} & I
\end{pmatrix}
\sum_{j=2}^NU^{*}_{j}V_{j}^{*\top}\mathbf{Y}_{n}$.  Now we can use the
decomposition $\mathbf{Y}_{n}=(\widetilde{\mathbf
  Y}_n+\widehat{\mathbf Y}_n)$ and the fact that
$U^{*}_{j}V_{j}^{*\top}\widetilde{\mathbf{Y}}_{n}= \mathbf{0}$ for any
$2\leq j\leq N$ (by \eqref{eq:relazioni-0} and \eqref{eq-def-Ytilde})
 in order to obtain the equality
\begin{equation*}
{\mathbf{N}}'_{n}\ =\
\begin{pmatrix} \mathbf{0} & I
\end{pmatrix}
\sum_{j=2}^NU^{*}_{j}V_{j}^{*\top}\widehat{\mathbf Y}_n\,.
\end{equation*}
Hence, the convergence rate and the second-order asymptotic
distribution of ${\mathbf{N}}'_{n}$ can be obtained by using the
convergences stated in Theorem \ref{clt-hat-Y} or in Lemma
\ref{conv-general-sum}. Specifically, case (a) follows from Theorem
\ref{clt-hat-Y}(a), observing that (by \eqref{eq:relazioni-0}) we have
$$
\begin{pmatrix}
\mathbf{0} & I
\end{pmatrix}\sum_{j=2}^NU^{*}_{j}V_{j}^{*\top}\begin{pmatrix}
\mathbf{0} & \mathbf{0}\\
\mathbf{0} & \widetilde{U}
\end{pmatrix}=
\begin{pmatrix}
\mathbf{0} & I
\end{pmatrix}\begin{pmatrix}
\mathbf{0} & \mathbf{0}\\
\mathbf{0} & \widetilde{U}_{-1}
\end{pmatrix}=\begin{pmatrix}
\mathbf{0} & \widetilde{U}_{-1}
\end{pmatrix}.
$$
Case (b) follows from Theorem \ref{clt-hat-Y}(b), observing that (by
\eqref{eq:relazioni-0}) we have
$$
\begin{pmatrix}
\mathbf{0} & I
\end{pmatrix}\sum_{j=2}^NU^{*}_{j}V_{j}^{*\top}\begin{pmatrix}
\widetilde{U} & \mathbf{0}\\
\mathbf{0} & \widetilde{U}\\
\end{pmatrix}
=\begin{pmatrix} \mathbf{0} & I
\end{pmatrix}\begin{pmatrix}
\widetilde{U}_{-1} & \mathbf{0}\\
\mathbf{0} & \widetilde{U}_{-1}\\
\end{pmatrix}
=\begin{pmatrix} \mathbf{0} & \widetilde{U}_{-1}
\end{pmatrix}.
$$
Finally, case (c) cannot be obtained directly by using the
convergences stated in Theorem \ref{clt-hat-Y} since in this case we
have (by \eqref{eq:relazioni-0})
$$
\begin{pmatrix}
\mathbf{0} & I
\end{pmatrix}\sum_{j=2}^NU^{*}_{j}V_{j}^{*\top}
\begin{pmatrix}
\mathbf{0}\\
\mathbf{1}
\end{pmatrix}=
N^{1/2}
\begin{pmatrix}
\mathbf{0} & I
\end{pmatrix}\sum_{j=2}^NU^{*}_{j}V_{j}^{*\top}
\begin{pmatrix}
\mathbf{0}\\
\mathbf{u}_1
\end{pmatrix}=
\begin{pmatrix}
\mathbf{0} & I
\end{pmatrix}
\begin{pmatrix}
\mathbf{0}\\
\mathbf{0}
\end{pmatrix}=
\mathbf{0}.
$$
Therefore, we need to express ${\mathbf{N}}'_{n}$ in the following
equivalent way:
\begin{equation*}
{\mathbf{N}}'_{n}\ =\
\begin{pmatrix}
\mathbf{0} & I
\end{pmatrix}
\sum_{j=2}^NU^{*}_{j}V_{j}^{*\top}\widehat{\mathbf Y}_n\ =\
\begin{pmatrix}
\mathbf{0} & I
\end{pmatrix}
\left(\sum_{j=2}^N\mathbf{Y}_{j(1),n}\ +\
\sum_{j=2}^N\mathbf{Y}_{j(2),n}\right),
\end{equation*}
where for the last equality we have used the decomposition
\eqref{decompo-Y-hat} of $\widehat{\mathbf{Y}}_n$ and the fact that
$U_{j}^*V_{j}^{*\top}\mathbf{Y}_{1(2),n}=
U_{j}^*V_{j}^{*\top}\mathbf{u}_{1(2)}\mathbf{v}_{1(2)}^{\top}\mathbf{Y}_n=
\mathbf{0}$ for $2\leq j\leq N$. Now, we recall that, in case (c),
that is $\nu>\gamma$, we have $g(\lambda_1)=g(1)=1$ and
$g(\lambda_j)=0$ for $2\leq j\leq N$ and so we get $\begin{pmatrix} \mathbf{0}
  & I
\end{pmatrix}
\mathbf{u}_{j(1)} = \mathbf{0}$ for $2\leq j\leq N$. As a
consequence, since
$\mathbf{Y}_{j(1),n}=\mathbf{u}_{j(1)}\mathbf{v}_{j(1)}^{\top}\mathbf{Y}_{n}$,
we have that $\begin{pmatrix} \mathbf{0} & I
\end{pmatrix}
\sum_{j=2}^N\mathbf{Y}_{j(1),n} = \mathbf{0}$, and the desired
convergence result follows from case $(vi)$ of Lemma
\ref{conv-general-sum}.
\qed

\appendix

\section{Computations for the proof of Lemma \ref{conv-general-sum}}
\label{lemmi-tecnici-appendice}
In all the sequel, given $(z_n)_n, (z'_n)_n$ two sequences of complex
numbers, the notation $z_n=O(z'_n)$ means $|z_n|\leq C |z'_n|$ for a
suitable constant $C>0$ and $n$ large enough. Moreover, if $z'_n\neq
0$, the notation $z_n\sim z z'_n$ with $z\in\mathbb{C}\setminus\{0\}$
means $\lim_n z_n/z'_n=z$ and, finally, the notation $z_n=o(z'_n)$
means $\lim_n z_n/z'_n=0$.  
\\

\indent Given $1/2<\delta\leq 1$, $x=a_x+i\,b_x\in{\mathbb C}$ with
$a_x>0$ and an integer $m_0\geq 2$ such that $a_x m^{-\delta}<1$ for
all $m\geq m_0$, let us set
\begin{equation}
p^{\delta}_{n}(x):=
\prod_{m=m_0}^{n}\left(1-\frac{x}{m^{\delta}}\right)
\quad \hbox{for } n\geq m_0\,.
\end{equation}

\subsection{Some technical results}
We first recall the following result, which has been proved in
\cite{ale-cri-ghi}.

\begin{lem}\cite[Lemma A.4]{ale-cri-ghi}\label{lemma-tecnico_1}
We have
\begin{equation}\label{affermazione1_bis}
|p^{\delta}_{n}(x)|\ =\ \begin{cases}
O\left(\exp \left( -a_x\frac{n^{1-\delta}}{1-\delta} \right)\right)
& \mbox{for } 1/2<\delta<1 \\
O\left(n^{ -a_x}\right) & \mbox{for } \delta=1
\end{cases}
\end{equation}
and
\begin{equation}\label{affermazione1_bis_inversa}
|p^{\delta}_{n}(x)^{-1}|\ =\ \begin{cases}
O\left(\exp \left( a_x\frac{ n^{1-\delta}}{1-\delta} \right)\right)
& \mbox{for } 1/2<\delta<1 \\
O\left(n^{a_x}\right) & \mbox{for } \delta=1\,.
\end{cases}
\end{equation}
\end{lem}
Therefore , if we set
\begin{equation}\label{eq:def-F}
F^{\delta}_{k+1,n}(x):=\frac{p^{\delta}_{n}(x)}{p^{\delta}_{k}(x)}
\qquad\hbox{for } m_0\leq k\leq n\,,
\end{equation}
 we have
\begin{equation}
\label{affermazione1_F}
|F^{\delta}_{k+1,n}(x)|\ =\ \begin{cases}
O\left(
\exp\left(
\frac{a_x }{1-\delta}(k^{1-\delta} - n^{1-\delta})
\right)
\right)
& \mbox{for } 1/2<\delta<1 \\
O\left( \left(\frac{k}{n}\right)^{a_x}
\right)
& \mbox{for } \delta=1.
\end{cases}
\end{equation}

Now, we prove two other results.

\begin{lem}\label{lemma-tecnico-O}
Given $\beta >1$ and $e>0$, we have
\begin{equation}
\sum_{k=m_0}^n \frac{1}{k^\beta}|F_{k+1,n}^{\delta}(x)|^e\ =\
\begin{cases}
O\left(n^{-(\beta-\delta)}\right)
\qquad &\hbox{if } 1/2<\delta<1,\\
O\left(n^{-ea_x}\right)
\qquad &\hbox{if } \delta=1\;\hbox{and } ea_x<\beta-1,\\
O\left(n^{-(\beta-1)}\ln(n)\right)
\qquad &\hbox{if } \delta=1\;\hbox{and } ea_x=\beta-1,\\
O\left(n^{-(\beta-1)}\right)
\qquad &\hbox{if } \delta=1\;\hbox{and } ea_x>\beta-1.
\end{cases}
\end{equation}
\end{lem}

\noindent{\bf Proof.} The desired relations immediately follows from
\eqref{affermazione1_F} using the well-known 
relation
\begin{equation}\label{relazione-nota}
\sum_{k=1}^n\frac{1}{k^{1-a}}\ =\
\begin{cases}
& O(1) \quad\mbox{for } a<0,\\
& \ln(n)+d+O(n^{-1})=\ln(n)+O(1) \quad\mbox{for } a=0,\\
&a^{-1}\, n^a + O(1) \quad\mbox{for } 0<a\leq 1,\\
&a^{-1}\, n^{a} + O(n^{a-1}) \quad\mbox{for } a>1,\,
\end{cases}
\end{equation}
where $d$ is the Euler-Mascheroni constant,
and the relation 
\begin{equation}\label{relazione-nota-2}
\begin{split}
\sum_{k=1}^n \frac{\exp(a k^b/b)}{k^{\beta}}
&=
O\left(\int_{1}^n \frac{\exp(a t^b/b)}{t^{\beta}}\,dt\right)
=O\left(
\left[\frac{\exp(a t^b/b)}{at^{b+\beta-1}}\right]_1^n
+\frac{(b+\beta-1)}{a}\int_{1}^n \frac{\exp(a t^b/b)}{t^{\beta+b}}\,dt
\right)\\
&=
O\left(
\frac{\exp(a n^b/b)}{n^{b+\beta-1}}
\right)
\qquad\hbox{for } a>0,\, b>0,\, \beta>1.
\end{split}
\end{equation}
Indeed, for the case $\delta=1$, it is enough to apply
\eqref{relazione-nota} with $a=ea_x-(\beta-1)$; while, for the case
$1/2<\delta<1$, it is enough to apply \eqref{relazione-nota-2} with
$a=ea_x,\, b=1-\delta$ and $\beta$.
\qed
\\

The following lemma extends \cite[Lemma A.5]{ale-cri-ghi}.

\begin{lem}\label{lemma-tecnico_2-nuovo}
Given $1/2<\delta_1\leq\delta_2\leq 1$, $\beta>\delta_1$ and
$x_1,x_2\in\mathbb{C}$ with $\mathcal{R}e(x_1)>0$,
$\mathcal{R}e(x_2)>0$, let $m_0\geq 2$ be an integer such that
$\max\{\mathcal{R}e(x_1),\, \mathcal{R}e(x_2)\}m^{-\delta_1}<1$ for
all $m\geq m_0$. Then we have
\begin{equation}\label{affermazione2-nuova}
\lim_n\,
n^{\beta-\delta_1} \sum_{k=m_0}^{n}
k^{-\beta} F_{k+1,n}^{\delta_1}(x_1) F^{\delta_2}_{k+1,n}(x_2)\  =\
\begin{cases}
\frac{1}{x_1+x_2 }\;
&\mbox{if } 1/2<\delta_1=\delta_2<1,
\\[3pt]
\frac{1}{x_1+x_2-\beta+1}\;
&\mbox{if } \delta_1=\delta_2=1\;\mbox{and } \mathcal{R}e(x_1+x_2)>\beta-1,
\\[3pt]
\frac{1}{x_1}\;
&\mbox{if } 1/2<\delta_1<\delta_2\leq 1.
\end{cases}
\end{equation}
\end{lem}

\noindent{\bf Proof.} Let us start with observing that, in each considered case,
relation \eqref{affermazione1_bis} implies 
\begin{equation}\label{lim-zero-nuovo}
\lim_n n^{\beta-\delta_1} |p^{\delta_1}_{n}(x_1)|\,|p^{\delta_2}_{n}(x_2)|=0.
\end{equation}
Indeed, in particular, when $\delta_1=\delta_2=1$ we have the
additional condition $\mathcal{R}e(x_1+x_2)>\beta-1$.
\\

Now, fix $k\geq 2$ and let us set $\eta:=\beta-\delta_1$ and
 $\ell^{\delta}_n(x):=1/p^{\delta}_n(x)$ and
define the following quantity
\[\begin{aligned}
D_{k}\ &&:=&\ \frac{1}{k^{\eta}} \ell^{\delta_1}_{k}(x_1)\ell^{\delta_2}_{k}(x_2) -
\frac{1}{(k-1)^{\eta}}\ell^{\delta_1}_{k-1}(x_1)\ell^{\delta_2}_{k-1}(x_2)\\
&&=&\
\left(\frac{1}{k^{\eta}} - \frac{1}{(k-1)^{\eta}} \right)
\ell^{\delta_1}_{k-1}(x_1)\ell^{\delta_2}_{k-1}(x_2)
+\ \frac{1}{k^{\eta}}\left(\ell^{\delta_1}_{k}(x_1)\ell^{\delta_2}_{k}(x_2) -
\ell^{\delta_1}_{k-1}(x_1)\ell^{\delta_2}_{k-1}(x_2)\right)\\
&&=&\
\ell^{\delta_1}_{k}(x_1)\ell^{\delta_2}_{k}(x_2)\left[
\left(\frac{1}{k^{\eta}} - \frac{1}{(k-1)^{\eta}} \right)
\frac{\ell^{\delta_1}_{k-1}(x_1)\ell^{\delta_2}_{k-1}(x_2)}{\ell^{\delta_1}_{k}(x_1)
\ell^{\delta_2}_{k}(x_2)}
+\ \frac{1}{k^{\eta}}
\left(1-\frac{\ell^{\delta_1}_{k-1}(x_1)\ell^{\delta_2}_{k-1}(x_2)}
{\ell^{\delta_1}_{k}(x_1)\ell^{\delta_2}_{k}(x_2)}\right)
\right].
\end{aligned}\]
Then, we observe the following:
\begin{equation}\label{eq:conti_1-nuovi}
\left(\frac{1}{k^{\eta}} - \frac{1}{(k-1)^{\eta}}\right)\ =\
-\frac{\eta}{k^{1+\eta}}
+O\left(\frac{1}{k^{2+\eta}}\right)
\quad\mbox{for } k\to +\infty
\end{equation}
and
\begin{equation}\label{eq:conti_2-nuovi}
\frac{\ell^{\delta_1}_{k-1}(x_1)\ell^{\delta_2}_{k-1}(x_2)}
{\ell^{\delta_1}_{k}(x_1)\ell^{\delta_2}_{k}(x_2)}\ =\
\left(1-\frac{x_1}{k^{\delta_1}}\right)
\left(1-\frac{x_2}{k^{\delta_2}}\right)\ =\
1+\frac{x_1x_2}{k^{(\delta_1+\delta_2)}}
-\frac{x_1}{k^{\delta_1}}-\frac{x_2}{k^{\delta_2}}.
\end{equation}
Now, by using~\eqref{eq:conti_1-nuovi} and \eqref{eq:conti_2-nuovi}
in the above expression of $D_{k}$, we have for $k\to +\infty$
\begin{equation*}
\begin{split}
D_{k}\ &=\
\ell^{\delta_1}_{k}(x_1)\ell^{\delta_2}_{k}(x_2)
\left[
\left(-\frac{\eta}{k^{\eta + 1}}+O(1/k^{\eta+2})\right)
\left(1+\frac{x_1x_2}{k^{(\delta_1+\delta_2)}}
-\frac{x_1}{k^{\delta_1}}-\frac{x_2}{k^{\delta_2}}
\right)
+\frac{1}{k^\eta}
\left(-\frac{x_1x_2}{k^{(\delta_1+\delta_2)}}
+\frac{x_1}{k^{\delta_1}}+\frac{x_2}{k^{\delta_2}}
\right)
\right]
\\
&=\
\begin{cases}
\ell^{\delta}_{k}(x_1)\ell^{\delta}_{k}(x_2)
\left[ \frac{x_1+x_2}{k^{\eta+\delta}}-\frac{\eta}{k^{\eta+1}}
+o(1/k^{\eta+\delta})\right]
\qquad\mbox{if } \delta_1=\delta_2=\delta,
\\
\ell^{\delta_1}_{k}(x_1)\ell^{\delta_2}_{k}(x_2)
\left[ \frac{x_1}{k^{\eta+\delta_1}}+o(1/k^{\eta+\delta_1})\right]
\qquad\mbox{if } \delta_1<\delta_2
\end{cases}
\\
&=\
\begin{cases}
\ell^{\delta}_{k}(x_1)\ell^{\delta}_{k}(x_2)
\frac{x_1+x_2}{k^{\eta+\delta}}
+o(\ell^{\delta}_{k}(x_1)\ell^{\delta}_{k}(x_2)/k^{\eta+\delta})
\qquad\mbox{if } \delta_1=\delta_2=\delta<1,
\\
\ell^{1}_{k}(x_1)\ell^{1}_{k}(x_2)
\frac{x_1+x_2-\eta}{k^{\eta+1}}
+o(\ell^{1}_{k}(x_1)\ell^{1}_{k}(x_2)/k^{\eta+1})
\qquad\mbox{if } \delta_1=\delta_2=1\;
\mbox{and } \mathcal{R}e(x_1+x_2)>\eta,
\\
\ell^{\delta_1}_{k}(x_1)\ell^{\delta_2}_{k}(x_2)
\frac{x_1}{k^{\eta+\delta_1}}
+o(\ell^{\delta_1}_{k}(x_1)\ell^{\delta_2}_{k}(x_2)/k^{\eta+\delta_1})
\qquad\mbox{if } \delta_1<\delta_2.
\end{cases}
\end{split}
\end{equation*}
that is
\begin{equation}\label{conti-nuovi}
D_{k}\sim\
\begin{cases}
\frac{x_1+x_2}{k^{\eta+\delta}} \,\ell^{\delta}_{k}(x_1)\ell^\delta_{k}(x_2)
\; &\mbox{if } 1/2<\delta_1=\delta_2=\delta<1,
\\
\frac{x_1+x_2-\eta}{k^{\eta+1}} \,\ell^{1}_{k}(x_1)\ell^1_{k}(x_2)
\; &\mbox{if } \delta_1=\delta_2=1\;\mbox{and } \mathcal{R}e(x_1+x_2)>\eta,
\\
\frac{x_1}{k^{\eta+\delta_1}}\, \ell^{\delta_1}_{k}(x_1)\ell^{\delta_2}_{k}(x_2)
\;&\mbox{if } 1/2<\delta_1<\delta_2\leq 1.
\end{cases}
\end{equation}
Now, following the same arguments used in the proof of \cite[Lemma
  A.5]{ale-cri-ghi}, in order to conclude, we apply \cite[Corollary
  A.2]{ale-cri-ghi} with
$$
z_n=D_{n},\quad v_n=n^\eta\,p^{\delta_1}_{n}(x_1)p^{\delta_2}_{n}(x_2),\quad
w_n=\frac{\ell_{n,1}\ell_{n,2}}{n^{\eta+\delta_1}D_{n}}, \quad
w=
\begin{cases}
\frac{1}{x_1+x_2}\; &\mbox{if } 1/2<\delta_1=\delta_2=\delta<1
\\
\frac{1}{x_1+x_2-\eta}\; &\mbox{if } \delta_1=\delta_2=1\;\mbox{and }
\mathcal{R}e(x_1+x_2)>\eta
\\
\frac{1}{x_1}\;&\mbox{if } 1/2<\delta_1<\delta_2\leq 1.
\end{cases}
$$
Indeed, $\lim_n v_n=0$ by \eqref{lim-zero-nuovo}, $\lim_n w_n=w\neq 0$
by \eqref{conti-nuovi},
\begin{equation*}
\begin{split}
\lim_nv_n\sum_{k=m_0}^n z_k&=\lim_n n^{\eta}p^{\delta_1}_{n}(x_1)p^{\delta_2}_{n}(x_2)
\sum_{k=m_0}^{n} D_{k}\\
&=
\lim_n n^{\eta}p^{\delta_1}_{n}(x_1)p^{\delta_2}_{2}(x_2)
\left(\frac{\ell^{\delta_1}_{n}(x_1)\ell^{\delta_2}_{n}(x_2)}{n^\eta}-
\frac{\ell^{\delta_1}_{m_0-1}(x_1)\ell^{\delta_2}_{m_0-1}(x_2)}{(m_0-1)^\eta}
\right)=1
\end{split}
\end{equation*}
by \eqref{lim-zero-nuovo} and
$z'_n=z_nw_n=r_n^2\ell_{n,1}\ell_{n,2}$.
\qed

\subsection{Analytic expression of $A^j_{k+1,n-1}$ with $j\geq 2$}
\label{subsection_A_FOGLIO_2}
Let us recall the definition of the following quantities for $j\geq 2$:
\begin{equation*}
A^j_{k+1,n-1}=\prod_{m=k+1}^{n-1}\left(I-D_{Q,j,m}\right),\qquad\mbox{
where }\qquad D_{Q,j,n}=\begin{pmatrix}
\widehat{r}_{n-1}(1-\lambda_j) & 0 \\
-\lambda_jh_n(\lambda_j) & \widehat{q}_{n,n}
\end{pmatrix}
\end{equation*}
\begin{equation*}
\mbox{with } h_n\; \mbox{defined in \eqref{eq:def-h}, that is }
 h_n(x)\ =\ \left\{\begin{aligned}
&\widehat{r}_{n-1}(1-x)\ &\mbox{ if }\nu<\gamma,\\
&\widehat{q}_{n,n}\ind_{\{x\neq1\}}\ &\mbox{ if }\nu\geq\gamma\,.
\end{aligned}\right.
\end{equation*}
The aim of this section is to compute the product above and so
finding the useful expression of $A^j_{k+1,n-1}$ presented
in~\eqref{eq:A_j_k_n_withG}, i.e.
\begin{equation*}
A^j_{k+1,n-1}=\begin{pmatrix}
F^{\gamma}_{k+1,n-1}(c\alpha_j) & 0 \\
\lambda_jG_{k+1,n-1}(c\alpha_j,q) & F^{\nu}_{k+1,n-1}(q)
\end{pmatrix},
\end{equation*}
where
\begin{equation*}
G_{k+1,n-1}(c\alpha_j,q)\ =\
\sum_{l=k+1}^{n-1}F^{\gamma}_{l+1,n-1}(c\alpha_j)h_l(\lambda_j)F^{\nu}_{k+1,l-1}(q).
\end{equation*}

It is straightforward to see that $[A^j_{k+1,n-1}]_{21}=0$,
$$[A^j_{k+1,n-1}]_{11}\ =\ \prod_{m=k+1}^{n-1}(1-\widehat{r}_{m-1}(1-\lambda_j))=
F^{\gamma}_{k+1,n-1}(c\alpha_j),$$
$$
[A^j_{k+1,n-1}]_{22}\ =\ \prod_{m=k+1}^{n-1}(1-\widehat{q}_{m,m})=F^{\nu}_{k+1,n-1}(q),
$$
while it is not immediate to determine $[A^j_{k+1,n-1}]_{12}$. To
this end, let us set $x_{n-1}:=[A^j_{k+1,n-1}]_{21}$ and observe that, since
$A^j_{k+1,n-1}=A^j_{k+1,n-2}(I-D_{Q,j,n-1})$ and
$x_{k+1}=\lambda_jh_{k+1}(\lambda_j)$, we have that
\begin{equation*}
\begin{aligned}
x_{n-1}&=x_{n-2}(1-\widehat{r}_{n-1}(1-\lambda_j)) +
[A^j_{k+1,n-2}]_{22}\lambda_jh_{n-1}(\lambda_j)\\
&=x_{n-2}F^{\gamma}_{n-1,n-1}(c\alpha_j) +
F^{\nu}_{k+1,n-2}(q)\lambda_jh_{n-1}(\lambda_j)\\
&=x_{n-3}F^{\gamma}_{n-2,n-1}(c\alpha_j) +
F^{\nu}_{k+1,n-3}(q)\lambda_jh_{n-2}(\lambda_j)F^{\gamma}_{n-1,n-1}(c\alpha_j)
+F^{\nu}_{k+1,n-2}(q)\lambda_jh_{n-1}(\lambda_j)\\
&=\ \dots\\
&=x_{k+1}F^{\gamma}_{k+2,n-1}(c\alpha_j)+
\sum_{l=k+2}^{n-2}F^{\nu}_{k+1,l-1}(q)
\lambda_jh_l(\lambda_j)F^{\gamma}_{l+1,n-1}(c\alpha_j)
+F^{\nu}_{k+1,n-2}(q)\lambda_jh_{n-1}(\lambda_j)\\
&=\sum_{l=k+1}^{n-1}F^{\nu}_{k+1,l-1}(q)
\lambda_jh_l(\lambda_j)F^{\gamma}_{l+1,n-1}(c\alpha_j)\\
&=\lambda_jG_{k+1,n-1}(c\alpha_j,q)\,.
\end{aligned}
\end{equation*}

\subsection{Asymptotic behavior  of $G_{k+1,n-1}(x,q)$}
\label{subsection_G_FOGLIO_2}
Let us recall the definition
\begin{equation*}
G_{k+1,n-1}(x,q)\ :=\
\sum_{l=k+1}^{n-1}F^{\gamma}_{l+1,n-1}(x)h_l(1-c^{-1}x)F^{\nu}_{k+1,l-1}(q).
\end{equation*}

Here we prove the following result:

\begin{lem}\label{lemma-tecnico-G}
When $\nu=\gamma$, we have for $x\in\mathbb{C}\setminus\{0\}$ 
\begin{equation}\label{lemma-nu-uguale-gamma}
G_{k+1,n-1}(x,q)=
\begin{cases}
\frac{q}{x-q}\left(F^{\gamma}_{k+1,n-1}(q)-F^{\gamma}_{k+1,n-1}(x)\right)
\,&\mbox{if } x\neq q,\\
\frac{q}{1-\gamma} F^{\gamma}_{k+1,n-1}(q)
\left[(n-1)^{1-\gamma}-(k+1)^{1-\gamma}\right]
+O(F^{\gamma}_{k+1,n-1}(q))
\, &\mbox{if } x=q\;\mbox{and } 1/2< \gamma<1,\\
q F^{\gamma}_{k+1,n-1}(q)\ln\left(\frac{n-1}{k+1}\right)
+O(k^{-1}F^{\gamma}_{k+1,n-1}(q))
\, &\mbox{if } x=q\;\mbox{and } \gamma=1.
\end{cases}
\end{equation}
When $\nu\neq \gamma$, we have for $x\in\mathbb{C}\setminus\{0\}$
\begin{equation}\label{lemma-nu-diverso-gamma}
G_{k+1,n-1}(x,q)=
C(x,q)
\left(
\frac{F^{\nu}_{k+1,n-1}(q)}{(n-1)^{\mu}}
-
\frac{F^{\gamma}_{k+1,n-1}(x)}{k^{\mu}}
\right)
+
O\left(\frac{|F^{\nu}_{k+1,n-1}(q)|}{n^{2\mu}}+
\frac{|F^{\gamma}_{k+1,n-1}(x)|}{k^{2\mu}}\right)\,,
\end{equation}
where $\mu:=|\gamma-\nu|$ and 
\begin{equation*}
C(x,q)\ :=
\begin{cases}
-\frac{x}{q}\qquad&\mbox{if } \nu<\gamma,\\
\frac{q}{x}\qquad&\mbox{if } \nu>\gamma.
\end{cases}
\end{equation*}
\end{lem}

\noindent{\bf Proof.}
Recalling the definition \eqref{eq:def-F}, we can write
\begin{equation}\label{espressione-G}
\begin{aligned}
G_{k+1,n-1}(x,q)&=
\sum_{l=k+1}^{n-1}F^{\gamma}_{l+1,n-1}(x)h_l(1-c^{-1}x)F^{\nu}_{k+1,l-1}(q)
=
\sum_{l=k+1}^{n-1}\frac{p^{\gamma}_{n-1}(x)}{p^{\gamma}_{l}(x)}
h_l(1-c^{-1}x)\frac{p^{\nu}_{l-1}(q)}{p^{\nu}_{k}(q)}\\
&=
\frac{p^{\gamma}_{n-1}(x)}{p^{\nu}_{k}(q)}\sum_{l=k+1}^{n-1}
\frac{h_l(1-c^{-1}x)}{(1-\widehat{q}_{l,l})}X_l,
\qquad\mbox{ where }\qquad
X_l\ :=\ \frac{p^{\nu}_{l}(q)}{p^{\gamma}_{l}(x)}.
\end{aligned}
\end{equation}
Moreover, recalling the definition \eqref{eq:def-h}, we have for $x\neq 0$
$$
h_l(1-c^{-1}x)=
\begin{cases}
\widehat{r}_{l-1}c^{-1}x=xl^{-\gamma}\qquad&\mbox{if }\nu<\gamma,\\
\widehat{q}_{l,l}=ql^{-\nu}\qquad&\mbox{if }\nu\geq \gamma.
\end{cases}
$$
Let us start with the case $\nu=\gamma$. In this case, we have 
\begin{equation*}
\begin{split}
\Delta X_l&:= X_{l}-X_{l-1}=
\left(1-\frac{X_{l-1}}{X_l}\right)X_l=
\left(\frac{x-q}{q}\frac{q}{l^{\gamma}(1-ql^{-\gamma})}\right)X_l\\
&=
\frac{x-q}{q}\frac{\widehat{q}_{l,l}}{1-\widehat{q}_{l,l}}X_l
=
\frac{x-q}{q}\frac{h_l(1-c^{-1}x)}{1-\widehat{q}_{l,l}}X_l.
\end{split}
\end{equation*}
It follows that
$$ 
\frac{x-q}{q}\sum_{l=k+1}^{n-1}
\frac{h_l(1-c^{-1}x)}{1-\widehat{q}_{l,l}}X_l=X_{n-1}-X_k.
$$
Since 
\begin{equation}\label{intermedio}
\begin{aligned}
&\frac{p^{\gamma}_{n-1}(x)}{p^{\nu}_{k}(q)}X_{n-1}\ =\
\frac{p^{\gamma}_{n-1}(x)}{p^{\nu}_{k}(q)}
\frac{p^{\nu}_{n-1}(q)}{p^{\gamma}_{n-1}(x)}\
=\
F^{\nu}_{k+1,n-1}(q)\;\mbox{and }\\
&\frac{p^{\gamma}_{n-1}(x)}{p^{\nu}_{k}(q)}X_{k}\ =\
\frac{p^{\gamma}_{n-1}(x)}{p^{\nu}_{k}(q)}
\frac{p^{\nu}_{k}(q)}{p^{\gamma}_{k}(x)}\ =\
F^{\gamma}_{k+1,n-1}(x),
\end{aligned}
\end{equation}
we find by \eqref{espressione-G}
$$
\frac{x-q}{q}G_{k+1,n-1}(x,q)=
\left(F^{\gamma}_{k+1,n-1}(q)-F^{\gamma}_{k+1,n-1}(x)\right)
$$
and so for $x\neq q$ we get 
$$
G_{k+1,n-1}(x,q)=
\frac{q}{x-q}\left(F^{\gamma}_{k+1,n-1}(q)-F^{\gamma}_{k+1,n-1}(x)\right).
$$ When $\nu=\gamma$ and $x=q$, we have $X_l=1$ and so we obtain (by 
\eqref{espressione-G} together with \eqref{relazione-nota})
\begin{equation*}
\begin{split}
G_{k+1,n-1}(x,q)&=q F^{\gamma}_{k+1,n-1}(q)
\sum_{l=k+1}^{n-1}\frac{1}{l^{\gamma}(1-ql^{-\gamma})}=
q F^{\gamma}_{k+1,n-1}(q)
\sum_{l=k+1}^{n-1}\frac{1}{l^{\gamma}}+O\left(\sum_{l\geq k+1}l^{-2\gamma}\right)
\\
&=q F^{\gamma}_{k+1,n-1}(q)
\begin{cases}
\frac{(n-1)^{1-\gamma}}{1-\gamma}-\frac{(k+1)^{1-\gamma}}{1-\gamma}
+O(1)+O(k^{-(2\gamma-1)})
\quad&\mbox{if } 1/2< \gamma<1
\\
\ln(n-1)-\ln(k+1)+O(n^{-1})+O(k^{-1})\quad&\mbox{if } \gamma=1,
\end{cases}
\end{split}
\end{equation*} 
which implies the two different asymptotic behavior in
\eqref{lemma-nu-uguale-gamma} according to the value of $\gamma$.
\\

\indent Now, let us consider the case $\nu\neq \gamma$ and introduce the
sequence $\{y_l;l\geq 1\}$ defined as $y_l:=l^{-\mu}$, with
$\mu=|\gamma-\nu|$.  Then, we have
\begin{equation*}
\begin{split}
y_lX_l-y_{l-1}X_{l-1}&=
\Delta y_l X_l + y_{l-1}\Delta X_l=
\left(
\frac{1}{l^\mu}-\frac{1}{(l-1)^{\mu}}
\right) 
X_l+
\left(
\frac{1}{l^\mu}+O\left(\frac{1}{l^{1+\mu}}\right)
\right)
\Delta X_l
\\
&=
\left(
\frac{-\mu}{l^{1+\mu}}+O\left(\frac{1}{l^{2+\mu}}\right)
\right)X_l+
\left(
\frac{1}{l^{\mu}}+O\left(\frac{1}{l^{1+\mu}}\right)
\right)
\Delta X_l\,,
\end{split}
\end{equation*}
where
\begin{equation*}
\Delta X_l:= X_{l}-X_{l-1}
=\left(1-\frac{X_{l-1}}{X_l}\right)X_l
=R_lX_l
\end{equation*}
with 
$$R_l:=\left(1-\frac{X_{l-1}}{X_l}\right)
=\frac{xl^{-\gamma}-ql^{-\nu}}{1-ql^{-\nu}}
=
\frac{\widehat{r}_{l-1}c^{-1}x-\widehat{q}_{l,l}}{1-\widehat{q}_{l,l}}
=O\left(\frac{1}{l^{\min\{\gamma,\nu\}}}\right).
$$ 
Taking into account that $\mu+\min\{\gamma,\nu\}<1+\mu$ for
$\nu\neq \gamma$, we obtain that
\begin{equation*}
y_lX_l-y_{l-1}X_{l-1}= 
\left[
\frac{R_l}{l^{\mu}}+O\left(\frac{1}{l^{1+\mu}}\right)
\right]X_l
=
K(x,q)
\frac{h_l(1-c^{-1}x)}{(1-\widehat{q}_{l,l})}X_l
+Q_lX_l,
\end{equation*}
where 
$$K(x,q)\ :=\
\left(-\frac{q}{x}\right)\ind_{\{\nu<\gamma\}}
+ 
\left(\frac{x}{q}\right)\ind_{\{\nu>\gamma\}}
=C(x,q)^{-1}
$$
and
$$
Q_l:=
\begin{cases}
\frac{x l^{-(2\gamma-\nu)}}{1-\widehat{q}_{l,l}}+O(l^{-(1+\mu)})
\qquad&\mbox{if } \nu<\gamma,\\
-\frac{ql^{-(2\nu-\gamma)}}{1-\widehat{q}_{l,l}}+O(l^{-(1+\mu)})
\qquad&\mbox{if } \nu>\gamma.\\
\end{cases}
$$ 
Note that $Q_l\sim \kappa l^{-(2\mu+\min\{\gamma,\nu\})}$ with a
suitable $\kappa\neq 0$.  The above expression implies that
\begin{equation}\label{aggiunta-numero}
\begin{split}
\frac{X_{n-1}}{(n-1)^{\mu}}-\frac{X_{k}}{k^{\mu}}=
\sum_{l=k+1}^{n-1}(y_lX_l-y_{l-1}X_{l-1})
=K(x,q)\sum_{l=k+1}^{n-1}
\frac{h_l(1-c^{-1}x)}{(1-\widehat{q}_{l,l})}X_l+
\sum_{l=k+1}^{n-1}Q_lX_l.
\end{split}
\end{equation}
With similar computations, setting 
$$
R_l^*:=1-\frac{|X_{l-1}|}{|X_l|}=
\frac{|1-ql^{-\nu}|-|1-xl^{-\gamma}|}{|1-ql^{-\nu}|}
$$ 
and taking into account that $R_l^*l^{-2\mu}\sim \kappa'
l^{-(2\mu+\min\{\gamma,\nu\})}$ with a suitable $\kappa'\neq 0$ and
$\min\{\gamma,\nu\}<1$ for $\nu\neq \gamma$, we find
$$
\frac{|X_l|}{l^{2\mu}}-\frac{|X_{l-1}|}{(l-1)^{2\mu}}=
\left[\frac{R_l^*}{l^{2\mu}}+
O\left(\frac{1}{l^{1+2\mu}}\right)\right]|X_l|
=Q_l^*|X_l|.
$$ Then, since $Q_l\sim \kappa'' Q_l^*$ with a suitable
$\kappa''\neq 0$, 
$$
\sum_{l=k+1}^{n-1}Q_lX_l
=O\left(\sum_{l=k+1}^{n-1}Q_l^*|X_l|\right)=
O\left(\frac{|X_{n-1}|}{(n-1)^{2\mu}}-\frac{|X_{k}|}{k^{2\mu}}\right).
$$ 
Finally, by \eqref{espressione-G}, \eqref{intermedio},
\eqref{aggiunta-numero} and the last above relations, we obtain for
$x\neq 0$
\begin{equation*}
G_{k+1,n-1}(x,q)
=
C(x,q)
\left(
\frac{F^{\nu}_{k+1,n-1}(q)}{(n-1)^{\mu}}
-
\frac{F^{\gamma}_{k+1,n-1}(x)}{k^{\mu}}
\right)
+
O\left(
\frac{|F^{\nu}_{k+1,n-1}(q)|}{n^{2\mu}}
+
\frac{|F^{\gamma}_{k+1,n-1}(x)|}{k^{2\mu}}
\right).
\end{equation*}
\qed

\subsection{Computation of the limit $d^{j_1(i_1),j_2(i_2)}$}
\label{result-limit-d}
Recall that we have 
\begin{equation*}
\begin{split}
&d_{k,n}^{j(1)}=
\begin{cases}
&\widehat{r}_{k-1}\quad\hbox{for } j=1\\
&\widehat{r}_{k-1}F^{\gamma}_{k+1,n-1}(c\alpha_{j})\quad\hbox{for } j\geq 2 
\end{cases}
\\
&\hbox{and}
\\
&d_{k,n}^{j(2)}=
\begin{cases}
&\left(\widehat{q}_{k,k}-\widehat{r}_{k-1}\right)F^{\nu}_{k+1,n-1}(q)
\;\hbox{for } j=1
\\
&\lambda_{j}\widehat{r}_{k-1}G_{k+1,n-1}(c\alpha_j,q)+
\left(\widehat{q}_{k,k}-\widehat{r}_{k-1}g(\lambda_{j})\right)F^{\nu}_{k+1,n-1}(q)
\;\hbox{for } j\geq 2\,,
\end{cases}
\end{split}
\end{equation*}
where $g$ is defined in \eqref{eq:def-g}, and so, for each $j\geq 2$,
we have $g(\lambda_j)=\lambda_j$ when $\nu<\gamma$, while $g(\lambda_j)=0$
when $\nu\geq\gamma$.
\\

Here, for each of the six cases $(i)-(vi)$ listed in Lemma
\ref{conv-general-sum}, we compute the limit
\begin{equation*}
d^{j_1(i_1),j_2(i_2)}\ =\ 
\lim_{n}t_n(J(I))^2
\sum_{k=m_0}^{n-1}d_{k,n}^{j_1(i_1)}d_{k,n}^{j_2(i_2)}\,.
\end{equation*}
For all the computations, we make the assumptions stated in Section 
\ref{section_model} and we use Lemma \ref{lemma-tecnico_2-nuovo} and
Lemma \ref{lemma-tecnico-G}.

\begin{description}
\item[Case $(i)$] 
Take $\nu<\gamma$, $j_1,j_2\in\{2,\dots, N\}$ and $i_1=i_2=1$. We have
\begin{equation*}
\begin{split}
\lim_n t_n(J(I))^2
\sum_{k=m_0}^{n-1}d_{k,n}^{j_1(i_1)}d_{k,n}^{j_2(i_2)}
&=
\lim_n
n^{\gamma}\sum_{k=m_0}^{n-1}\widehat{r}_{k-1}^2
F^{\gamma}_{k+1,n-1}(c\alpha_{j_1})F^{\gamma}_{k+1,n-1}(c\alpha_{j_2})
\\
&= 
c^2
\lim_n
n^{\gamma}\sum_{k=m_0}^{n-1} 
k^{-2\gamma} F^{\gamma}_{k+1,n-1}(c\alpha_{j_1}) F^{\gamma}_{k+1,n-1}(c\alpha_{j_2})
\\
&=
\frac{c^2}{c(\alpha_{j_1}+\alpha_{j_2})-\ind_{\{\gamma=1\}}}\,. 
\end{split}
\end{equation*}

\item[Case $(ii)$] 
Take $\nu<\gamma$, $j_1,j_2\in\{1,\dots, N\}$ and $i_1=i_2=2$. 
For $j_1=j_2=1$, we have
\begin{equation*}
\begin{split}
\lim_n t_n(J(I))^2
\sum_{k=m_0}^{n-1}d_{k,n}^{j_1(i_1)}d_{k,n}^{j_2(i_2)}
&=
\lim_n n^{\nu} \sum_{k=m_0}^{n-1}
\left(\widehat{q}_{k,k}-\widehat{r}_{k-1}\right)^2F^{\nu}_{k+1,n-1}(q)^2
\\
&=\lim_n
n^{\nu}\sum_{k=m_0}^{n-1}\widehat{q}_{k,k}^2 F^{\nu}_{k+1,n-1}(q)^2
\\
&=q^2
\lim_n n^{\nu}\sum_{k=m_0}^{n-1} 
k^{-2\nu} F^{\nu}_{k+1,n-1}(q)^2 
=\frac{q}{2}.
\end{split}
\end{equation*}
(Note that the above second equality is due to the fact that some
terms are $o(n^{-\nu})$ and so we can cancel them.)  Similarly, for
the cases $j_1\geq 2,\,j_2\geq 2$ and $j_1=1,\,j_2\geq 2$ and $j_1\geq
2,\,j_2=1$, using Lemma \ref{lemma-tecnico-G}, which allows us to
replace in the computation of the desired limit the quantity
$G_{k+1,n-1}(c\alpha_j,q)$ by
\begin{equation*}
-\frac{c\alpha_j}{q}
\left(
\frac{F^{\nu}_{k+1,n-1}(q)}{(n-1)^{\gamma-\nu}}
-
\frac{F^{\gamma}_{k+1,n-1}(c\alpha_j)}{k^{\gamma-\nu}}
\right),
\end{equation*}
and removing the terms which are $o(n^{-\nu})$, we obtain
\begin{equation*}
\begin{split}
\lim_n t_n(J(I))^2
\sum_{k=m_0}^{n-1}d_{k,n}^{j_1(i_1)}d_{k,n}^{j_2(i_2)}
&=
\lim_n
n^{\nu}\sum_{k=m_0}^{n-1}\widehat{q}_{k,k}^2 F^{\nu}_{k+1,n-1}(q)^2
\\
&=q^2
\lim_n n^{\nu}\sum_{k=m_0}^{n-1} 
k^{-2\nu} F^{\nu}_{k+1,n-1}(q)^2 
=\frac{q}{2}.
\end{split}
\end{equation*}

\item[Case $(iii)$] Take $\nu=\gamma$, $j_1,j_2\in\{1,\dots,N\}$ and
  $i_1,i_2\in\{1,2\}$ with $i_h\neq 1$ if $j_h=1$. Recall assumption
  \eqref{ass-q} \footnote{In the case $q=c\alpha_j$ for some $j\geq 2$
    the computations are similar, but we have to consider the other
    asymptotic expression given in Lemma
    \ref{lemma-tecnico-G}.}. Therefore, for $j_1=j_2=1$ and
  $i_1=i_2=2$, we have
\begin{equation*}
\begin{split}
\lim_n t_n(J(I))^2
\sum_{k=m_0}^{n-1}d_{k,n}^{j_1(i_1)}d_{k,n}^{j_2(i_2)}
&=
\lim_n
n^{\gamma}\sum_{k=m_0}^{n-1}
\left(\widehat{q}_{k,k}-\widehat{r}_{k-1}\right)^2F^{\gamma}_{k+1,n-1}(q)^2
\\&=
\lim_n (q-c)^2
n^{\gamma}\sum_{k=m_0}^{n-1}
\frac{1}{k^{2\gamma}}F^{\gamma}_{k+1,n-1}(q)^2
=\frac{(q-c)^2}{2q-\ind_{\{\gamma=1\}}}.
\end{split}
\end{equation*}
For $j_1=1$, $j_2\geq 2$, $i_1=2$ and $i_2=1$, we have 
\begin{equation*}
\begin{split}
&\lim_n t_n(J(I))^2
\sum_{k=m_0}^{n-1}d_{k,n}^{j_1(i_1)}d_{k,n}^{j_2(i_2)}
=\\
&\lim_n
n^{\gamma}\sum_{k=m_0}^{n-1}
\left(\widehat{q}_{k,k}-\widehat{r}_{k-1}\right)F^{\gamma}_{k+1,n-1}(q)
\widehat{r}_{k-1}F^{\gamma}_{k+1,n-1}(c\alpha_{j_2})
=
\\
&
(q-c)c\lim_n
n^{\gamma}\sum_{k=m_0}^{n-1}
\frac{1}{k^{2\gamma}}F^{\gamma}_{k+1,n-1}(q)F^{\gamma}_{k+1,n-1}(c\alpha_{j_2})
=
\frac{c(q-c)}{c\alpha_{j_2}+q-\ind_{\{\gamma=1\}}}.
\end{split}
\end{equation*}
By symmetry, for $j_1\geq 2$, $j_2=1$, $i_1=1$ and $i_2=2$, we have 
$$
d^{j_1(i_1),j_2(i_2)}=\frac{c(q-c)}{c\alpha_{j_1}+q-\ind_{\{\gamma=1\}}}.
$$
For $j_1=1$, $j_2\geq 2$ and $i_1=i_2=2$, we observe that, by
means of Lemma \ref{lemma-tecnico-G}, in the computation of the
considered limit, we can replace $G_{k+1,n-1}(c\alpha_j,q)$ by
$$q(c\alpha_j-q)^{-1}\big(F^{\gamma}_{k+1,n-1}(q)
-F^{\gamma}_{k+1,n-1}(c\alpha_j)\big),$$ that is we can replace 
$d_{k,n}^{j(2)}$, with $j\geq 2$ by 
$$
\widehat{q}_{k,k}
\left(\frac{(c\alpha_j-c)F^{\gamma}_{k+1,n-1}(c\alpha_j) - 
(q-c)F^{\gamma}_{k+1,n-1}(q)}{c\alpha_j-q}
\right).
$$
Therefore, we have 
\begin{equation*}
\begin{split}
&\lim_n t_n(J(I))^2
\sum_{k=m_0}^{n-1}d_{k,n}^{j_1(i_1)}d_{k,n}^{j_2(i_2)}
=\\
&\lim_n
n^{\gamma}\sum_{k=m_0}^{n-1}
\left(\widehat{q}_{k,k}-\widehat{r}_{k-1}\right)F^{\gamma}_{k+1,n-1}(q)
\widehat{q}_{k,k}
\left(\frac{(c\alpha_{j_2}-c)F^{\gamma}_{k+1,n-1}(c\alpha_{j_2}) - 
(q-c)F^{\gamma}_{k+1,n-1}(q)}{c\alpha_{j_2}-q}
\right)
=
\\
&\frac{q(q-c)(c\alpha_{j_2}-c)}{c\alpha_{j_2}-q}
\lim_n
n^{\gamma}\sum_{k=m_0}^{n-1}\frac{1}{k^{2\gamma}}
F^{\gamma}_{k+1,n-1}(c\alpha_{j_2})F^{\gamma}_{k+1,n-1}(q)
\\
&\qquad\qquad+
\frac{q(q-c)^2}{c\alpha_{j_2}-q}
\lim_n
n^{\gamma}\sum_{k=m_0}^{n-1}\frac{1}{k^{2\gamma}}F^{\gamma}_{k+1,n-1}(q)^2
=\\
&\frac{q(q-c)(c+q-\ind_{\{\gamma=1\}})}
{(c\alpha_{j_2}+q-\ind_{\{\gamma=1\}})(2q-\ind_{\{\gamma=1\}})}.
\end{split}
\end{equation*}
By symmetry,  for $j_1\geq 2$, $j_2=1$ and $i_1=i_2=2$, we get
\begin{equation*}
d^{j_1(i_1),j_2(i_2)}=
\frac{q(q-c)(c+q-\ind_{\{\gamma=1\}})}
{(c\alpha_{j_1}+q-\ind_{\{\gamma=1\}})(2q-\ind_{\{\gamma=1\}})}.
\end{equation*}
Similarly, for $j_1\geq 2$, $j_2\geq 2$, $i_1=i_2=1$, we have 
\begin{equation*}
\begin{split}
\lim_n t_n(J(I))^2
\sum_{k=m_0}^{n-1}d_{k,n}^{j_1(i_1)}d_{k,n}^{j_2(i_2)}
&=
\lim_n n^{\gamma}\sum_{k=m_0}^{n-1}
\widehat{r}_{k-1}^2
F^{\gamma}_{k+1,n-1}(c\alpha_{j_1})F^{\gamma}_{k+1,n-1}(c\alpha_{j_2})
\\&=
\frac{c^2}{c(\alpha_{j_1}+\alpha_{j_2})-\ind_{\{\gamma=1\}}}.
\end{split}
\end{equation*}
For $j_1\geq 2$, $j_2\geq 2$, $i_1=1$ and $i_2=2$, we have 
\begin{equation*}
\begin{split}
&\lim_n t_n(J(I))^2
\sum_{k=m_0}^{n-1}d_{k,n}^{j_1(i_1)}d_{k,n}^{j_2(i_2)}
=\\
&\lim_n n^{\gamma}\sum_{k=m_0}^{n-1}
\widehat{r}_{k-1}
F^{\gamma}_{k+1,n-1}(c\alpha_{j_1})\widehat{q}_{k,k}
\left(\frac{(c\alpha_{j_2}-c)F^{\gamma}_{k+1,n-1}(c\alpha_{j_2}) - 
(q-c)F^{\gamma}_{k+1,n-1}(q)}{c\alpha_{j_2}-q}\right)
=
\\
&\frac{cq(c\alpha_{j_1}+c-\ind_{\{\gamma=1\}})}
{(c\alpha_{j_1}+c\alpha_{j_2}-\ind_{\{\gamma=1\}})(c\alpha_{j_1}+q-\ind_{\{\gamma=1\}})}.
\end{split}
\end{equation*}
By symmetry, for $j_1\geq 2$, $j_2\geq 2$, $i_1=2$ and $i_2=1$, we get
$$
d^{j_1(i_1),j_2(i_2)}=\frac{cq(c\alpha_{j_2}+c-\ind_{\{\gamma=1\}})}
{(c\alpha_{j_1}+c\alpha_{j_2}-\ind_{\{\gamma=1\}})(c\alpha_{j_2}+q-\ind_{\{\gamma=1\}})}.
$$
Finally, for $j_1\geq 2$, $j_2\geq 2$ and $i_1=i_2=2$, we have 
\begin{equation*}
\begin{split}
&\lim_n t_n(J(I))^2
\sum_{k=m_0}^{n-1}d_{k,n}^{j_1(i_1)}d_{k,n}^{j_2(i_2)}
=\\
&\lim_n n^{\gamma}\!\!\sum_{k=m_0}^{n-1}\!\!
\widehat{q}_{k,k}^2\!
\left(\!\frac{(c\alpha_{j_1}-c)F^{\gamma}_{k+1,n-1}(c\alpha_{j_1}) - 
(q-c)F^{\gamma}_{k+1,n-1}(q)}{c\alpha_{j_1}-q}
\!\right)\\
& \qquad\qquad
\left(\!\frac{(c\alpha_{j_2}-c)F^{\gamma}_{k+1,n-1}(c\alpha_{j_2}) - 
(q-c)F^{\gamma}_{k+1,n-1}(q)}{c\alpha_{j_2}-q}
\!\right)=
\\
&q^2\frac{c^3(\alpha_{j_1}+\alpha_{j_2})
+2c^2q(\alpha_{j_1}\alpha_{j_2}+1)
-\ind_{\{\gamma=1\}}c^2(\alpha_{j_1}\alpha_{j_2}+\alpha_{j_1}+\alpha_{j_2}+2)}
{(2q-\ind_{\{\gamma=1\}})(c(\alpha_{j_1}+\alpha_{j_2})
-\ind_{\{\gamma=1\}})
(c\alpha_{j_1}+q-\ind_{\{\gamma=1\}})(c\alpha_{j_2}+q-\ind_{\{\gamma=1\}})}
\\
&\qquad+ 
q^2\frac{c(q-\ind_{\{\gamma=1\}})^2(\alpha_{j_1}+\alpha_{j_2})
-\ind_{\{\gamma=1\}}(2c+q-1)(q-1)}
{(2q-\ind_{\{\gamma=1\}})(c(\alpha_{j_1}+\alpha_{j_2})
-\ind_{\{\gamma=1\}})
(c\alpha_{j_1}+q-\ind_{\{\gamma=1\}})(c\alpha_{j_2}+q-\ind_{\{\gamma=1\}})}.
\end{split}
\end{equation*}

\item[Case $(iv)$] Take $\gamma<\nu$, $j_1,j_2\in\{2,\dots,
  N\}$ and $i_1=i_2=1$. We have
\begin{equation*}
\begin{split}
\lim_n t_n(J(I))^2
\sum_{k=m_0}^{n-1}d_{k,n}^{j_1(i_1)}d_{k,n}^{j_2(i_2)}
&=
\lim_n n^{\gamma}\sum_{k=m_0}^{n-1}\widehat{r}_{k-1}^2 
F^{\gamma}_{k+1,n-1}(c\alpha_{j_1})
F^{\gamma}_{k+1,n-1}(c\alpha_{j_2})
\\
&=
\lim_n n^{\gamma}\sum_{k=m_0}^{n-1}k^{-2\gamma}F^{\gamma}_{k+1,n-1}(c\alpha_{j_1})
F^{\gamma}_{k+1,n-1}(c\alpha_{j_2})
=\frac{c}{\alpha_{j_1}+\alpha_{j_2}}\,.
\end{split}
\end{equation*}
The difference with the computations in the case $\nu<\gamma$
concerns only the fact that here it is not possible that
$\gamma=1$ since $\gamma<\nu\leq 1$.

\item[Case $(v)$] Take $\gamma<\nu$, $j_1,j_2=1$ and $i_1=i_2=2$. We
  have
\begin{equation*}
\begin{split}
\lim_n t_n(J(I))^2
\sum_{k=m_0}^{n-1}d_{k,n}^{j_1(i_1)}d_{k,n}^{j_2(i_2)}
&=
\lim_nn^{2\gamma-\nu}\sum_{k=m_0}^{n-1}
\left(\widehat{q}_{k,k}-\widehat{r}_{k-1}\right)^2F^{\nu}_{k+1,n-1}(q)^2
\\
&=
\lim_n n^{2\gamma-\nu}
\sum_{k=m_0}^{n-1}\widehat{r}_{k-1}^2F^{\nu}_{k+1,n-1}(q)^2
=
\frac{c^2}{2q-\ind_{\{\nu=1\}}(2\gamma-1)}\,.
\end{split}
\end{equation*}
(Note that the above second equality is due to the fact that some
terms are $o(n^{- (2\gamma-\nu) })$ and so we can cancel them.)  

\item[Case $(vi)$] Take $\gamma<\nu$, $j_1,j_2\in\{2,\dots,N\}$ and
  $i_1=i_2=2$. Using Lemma \ref{lemma-tecnico-G}, which allows us to
  replace in the computation of the desired limit the quantity  
  $G_{k+1,n-1}(c\alpha_j,q)$ by
$$
\frac{q}{c\alpha_j}
\left(
\frac{F^{\nu}_{k+1,n-1}(q)}{(n-1)^{\nu-\gamma}}
-
\frac{F^{\gamma}_{k+1,n-1}(c\alpha_j)}{k^{\nu-\gamma}}
\right),
$$ and removing the terms which are $o(n^{-\nu})$, we have
\begin{equation*}
\begin{split}
&\lim_n t_n(J(I))^2
\sum_{k=m_0}^{n-1}d_{k,n}^{j_1(i_1)}d_{k,n}^{j_2(i_2)}=
\\
&
\lim_n n^{\nu}\sum_{k=m_0}^{n-1}
\left(\lambda_{j_1}\widehat{r}_{k-1}G_{k+1,n-1}(c\alpha_{j_1},q)+
\widehat{q}_{k,k}F^{\nu}_{k+1,n-1}(q)\right)
\left(
\lambda_{j_2}\widehat{r}_{k-1}G_{k+1,n-1}(c\alpha_{j_2},q)+
\widehat{q}_{k,k}F^{\nu}_{k+1,n-1}(q)\right)
=\\
&\frac{\lambda_{j_1}\lambda_{j_2}}{\alpha_{j_1}\alpha_{j_2}}q^2
\lim_n n^{2\gamma-\nu}\sum_{k=m_0}^{n-1}k^{-2\gamma}F^{\nu}_{k+1,n-1}(q)^2
+
\left(
\frac{\lambda_{j_1}}{\alpha_{j_1}}+\frac{\lambda_{j_2}}{\alpha_{j_2}}
\right)
q^2
\lim_n n^{\gamma}\sum_{k=m_0}^{n-1}k^{-(\gamma+\nu)}F^{\nu}_{k+1,n-1}(q)^2
\\
&+q^2\lim_n n^{\nu}\sum_{k=m_0}^{n-1}k^{-2\nu}F^{\nu}_{k+1,n-1}(q)^2
=
\\
&
\left(\frac{\lambda_{j_1}\lambda_{j_2}}{\alpha_{j_1}\alpha_{j_2}}\right)
\frac{q^2}{2q-\ind_{\{\nu=1\}}(2\gamma-1)} \ +\
\left(\frac{\lambda_{j_1}}{\alpha_{j_1}}+
\frac{\lambda_{j_2}}{\alpha_{j_2}}\right)\frac{q^2}{2q-\ind_{\{\nu=1\}}\gamma}
\ +\ \frac{q^2}{2q-\ind_{\{\nu=1\}}}\,.
\end{split}
\end{equation*}
\end{description}

\section{Stable convergence and its variants}\label{app-B}
This brief appendix contains some basic definitions and results
concerning stable convergence and its variants. For more details,
we refer the reader to \cite{crimaldi-2009, crimaldi-libro,
  cri-let-pra-2007, hall-1980} and the references therein.\\

\indent Let $(\Omega, {\mathcal A}, P)$ be a probability space,
and let $S$ be a Polish space, endowed with its Borel
$\sigma$-field. A {\em
  kernel} on $S$, or a random probability measure on $S$, is a
collection $K=\{K(\omega):\, \omega\in\Omega\}$ of probability
measures on the Borel $\sigma$-field of $S$ such that, for each
bounded Borel real function $f$ on $S$, the map
$$
\omega\mapsto K\!f(\omega)=\int f (x)\, K(\omega)(dx)
$$
is $\mathcal A$-measurable. Given a sub-$\sigma$-field $\mathcal
H$ of $\mathcal A$, a kernel $K$ is said $\mathcal H$-measurable
if all the
above random variables $K\!f$ are $\mathcal H$-measurable.\\

\indent On $(\Omega, {\mathcal A},P)$, let $(Y_n)_n$ be a sequence
of $S$-valued random variables, let $\mathcal H$ be a
sub-$\sigma$-field of $\mathcal A$, and let $K$ be a $\mathcal
H$-measurable kernel on $S$. Then we say that $Y_n$ converges {\em
$\mathcal H$-stably} to $K$, and we write $Y_n\longrightarrow K$
${\mathcal H}$-stably, if
$$
P(Y_n \in \cdot \,|\, H)\stackrel{weakly}\longrightarrow
E\left[K(\cdot)\,|\, H \right] \qquad\hbox{for all } H\in{\mathcal
H}\; \hbox{with } P(H) > 0,
$$where $K(\cdot)$ denotes the random variable
  defined, for each Borel set $B$ of $S$, as $\omega\mapsto
  K\!I_B(\omega)=K(\omega)(B)$.  In the case when ${\mathcal
  H}={\mathcal A}$, we simply say that $Y_n$ converges {\em stably} to
$K$ and we write $Y_n\longrightarrow K$ stably. Clearly, if
$Y_n\longrightarrow K$ ${\mathcal H}$-stably, then $Y_n$ converges
in distribution to the probability distribution $E[K(\cdot)]$.
Moreover, the $\mathcal H$-stable convergence of $Y_n$ to $K$ can
be stated in terms of the following convergence of conditional
expectations:
\begin{equation}\label{def-stable}
E[f(Y_n)\,|\, {\mathcal H}]\stackrel{\sigma(L^1,\,
L^{\infty})}\longrightarrow K\!f
\end{equation}
for each bounded continuous real function $f$ on $S$. \\

\indent In \cite{cri-let-pra-2007} the notion of $\mathcal
H$-stable convergence is firstly generalized in a natural way
replacing in (\ref{def-stable}) the single sub-$\sigma$-field
$\mathcal H$ by a collection ${\mathcal G}=({\mathcal G}_n)_n$
(called conditioning system) of sub-$\sigma$-fields of $\mathcal
A$ and then it is strengthened by substituting the convergence in
$\sigma(L^1,L^{\infty})$ by the one in probability (i.e. in $L^1$,
since $f$ is bounded). Hence, according to
\cite{cri-let-pra-2007}, we say that $Y_n$ converges to $K$ {\em
stably in the strong sense}, with respect to ${\mathcal
G}=({\mathcal G}_n)_n$, if
\begin{equation}\label{def-stable-strong}
E\left[f(Y_n)\,|\,{\mathcal G}_n\right]\stackrel{P}\longrightarrow
K\!f
\end{equation}
for each bounded continuous real function $f$ on $S$.\\

\indent Finally, a strengthening of the stable convergence in the
strong sense can be naturally obtained if in
(\ref{def-stable-strong}) we replace the convergence in
probability by the almost sure convergence: given a conditioning
system ${\mathcal G}=({\mathcal
  G}_n)_n$, we say that $Y_n$ converges to $K$ in the sense of the
{\em almost sure conditional convergence}, with respect to
${\mathcal
  G}$, if
\begin{equation*}
E\left[f(Y_n)\,|\,{\mathcal
G}_n\right]\stackrel{a.s.}\longrightarrow K\!f
\end{equation*}
for each bounded continuous real function $f$ on
  $S$. The almost sure conditional convergence has been introduced in
  \cite{crimaldi-2009} and, subsequently, employed by others in the
  urn model literature.  \\

We now conclude this section recalling two convergence results
that we
need in our proofs. \\

From \cite[Proposition~3.1]{cri-pra}, we can get the following result.

\begin{theo}\label{thm:triangular}
Let $({\mathbf T}_{k,n})_{ 1\leq k\leq k_n,\, n\geq 1}$ be a
triangular array of $d$-dimensional real random vectors, such
that, for each fixed $n$, the finite sequence $({\mathbf
T}_{k,n})_{1\leq k\leq k_n}$ is a martingale difference array with
respect to a given filtration $({\mathcal G}_{k,n})_{k\geq 0}$.
Moreover, let $(t_n)_n$ be a sequence of real numbers and assume
that the following conditions hold:
\begin{itemize}
\item[(c1)] ${\mathcal G}_{k,n}{\underline{\subset}} {\mathcal G}_{k, n+1}$ 
for each $n$ and $1\leq k\leq k_n$;
\item[(c2)] $\sum_{k=1}^{k_n} (t_n{\mathbf
  T}_{k,n})(t_n{\mathbf T}_{k,n})^{\top}=t_n^2\sum_{k=1}^{k_n} {\mathbf
  T}_{k,n}{\mathbf T}_{k,n}^{\top} \stackrel{P}\longrightarrow \Sigma$,
  where $\Sigma$ is a random positive semi\-defi\-ni\-te matrix;
\item[(c3)] $\sup_{1\leq k\leq k_n} |t_n{\mathbf T}_{k,n}|
\stackrel{L^1}\longrightarrow 0$.
\end{itemize}
Then $t_n\sum_{k=1}^{k_n}{\mathbf T}_{k,n}$ converges stably to
the Gaussian kernel ${\mathcal N}(\mathbf{0}, \Sigma)$.
\end{theo}

The following result combines together a stable convergence and a
stable convergence in the strong sense.

\begin{theo}\cite[Lemma 1]{ber-cri-pra-rig}\label{blocco}
Suppose that $C_n$ and $D_n$ are $S$-valued random variables, that
$M$ and $N$ are kernels on $S$, and that ${\mathcal G}=({\mathcal
G}_n)_n$ is a filtration satisfying for all $n$
$$
\sigma(C_n)\underline\subset{\mathcal G}_n\quad\hbox{and }\quad
\sigma(D_n)\underline\subset
\sigma\left({\textstyle\bigcup_n}{\mathcal G}_n\right)
$$

\noindent If $C_n$ stably converges to $M$ and $D_n$ converges to
$N$ stably in the strong sense, with respect to $\mathcal G$, then
$$
(C_n, D_n)\longrightarrow M \otimes N \qquad\hbox{stably}.
$$
(Here, $M\otimes N$ is the kernel on $S\times S$ such that $(M
\otimes N )(\omega) = M(\omega) \otimes N(\omega)$ for all
$\omega$.)
\end{theo}


\begin{thebibliography}{10}

\bibitem{ale-cri-ghi-MEAN}
G.~Aletti, I.~Crimaldi, and A.~Ghiglietti.
\newblock Networks of reinforced stochastic processes: asymptotics for the
  empirical means.
\newblock {\em arXiv: 1705.02126 (2017). Forthcoming in Bernoulli}.

\bibitem{ale-cri-ghi}
G.~Aletti, I.~Crimaldi, and A.~Ghiglietti.
\newblock Synchronization of reinforced stochastic processes with a
  network-based interaction.
\newblock {\em Ann. Appl. Probab.}, 27:3787--3844, 2017.

\bibitem{ale-ghi}
G.~Aletti and A.~Ghiglietti.
\newblock Interacting generalized {F}riedman's urn systems.
\newblock {\em Stochastic Process. Appl.}, 127:2650--2678, 2017.

\bibitem{ale-ghi-ros}
G.~Aletti, A.~Ghiglietti, and W.~F. Rosenberger.
\newblock Nonparametric covariate-adjusted response-adaptive design based on a
  functional urn model.
\newblock {\em Ann. Statist.}, 46(6B):3838--3866, 12 2018.

\bibitem{ale-ghi-vid}
G.~Aletti, A.~Ghiglietti, and A.~N. Vidyashankar.
\newblock Dynamics of an adaptive randomly reinforced urn.
\newblock {\em Bernoulli}, 24(3):2204--2255, 2018.

\bibitem{ben}
M.~Bena{\"{\i}}m, I.~Benjamini, J.~Chen, and Y.~Lima.
\newblock A generalized {P}\'olya's urn with graph based interactions.
\newblock {\em Random Struct. Algor.}, 46(4):614--634, 2015.

\bibitem{ber-cri-pra-rig}
P.~Berti, I.~Crimaldi, L.~Pratelli, and P.~Rigo.
\newblock A central limit theorem and its applications to multicolor randomly
  reinforced urns.
\newblock {\em J. Appl. Probab.}, 48(2):527--546, 2011.

\bibitem{ber-cri-pra-rig-barriere}
P.~Berti, I.~Crimaldi, L.~Pratelli, and P.~Rigo.
\newblock Asymptotics for randomly reinforced urns with random barriers.
\newblock {\em J. Appl. Probab.}, 53(4):1206--1220, 2016.

\bibitem{che-luc}
J.~Chen and C.~Lucas.
\newblock A generalized {P}\'olya's urn with graph based interactions:
  convergence at linearity.
\newblock {\em Electron. Commun. Probab.}, 19:no. 67, 13, 2014.

\bibitem{chen-kuba}
M.-R. Chen and M.~Kuba.
\newblock On generalized {P}\'olya urn models.
\newblock {\em J. Appl. Prob.}, 50:1169--1186, 2013.

\bibitem{cir}
P.~Cirillo, M.~Gallegati, and J.~H{\"u}sler.
\newblock A {P}\'olya lattice model to study leverage dynamics and contagious
  financial fragility.
\newblock {\em Adv. Complex Syst.}, 15(suppl. 2):1250069, 26, 2012.

\bibitem{collevecchio}
A.~Collevecchio, C.~Cotar, and M.~LiCalzi.
\newblock On a preferential attachment and generalized {P}\'olya's urn model.
\newblock {\em Ann. Appl. Prob.}, 23:1219--1253, 2013.

\bibitem{crimaldi-2009}
I.~Crimaldi.
\newblock An almost sure conditional convergence result and an application to a
  generalized {P}\'olya urn.
\newblock {\em Int. Math. Forum}, 4(21-24):1139--1156, 2009.

\bibitem{cri-ipergeom}
I.~Crimaldi.
\newblock Central limit theorems for a hypergeometric randomly reinforced urn.
\newblock {\em J. Appl. Prob.}, 53(3):899--913, 2016.

\bibitem{crimaldi-libro}
I.~Crimaldi.
\newblock {\em Introduzione alla nozione di convergenza stabile e sue varianti
  (Introduction to the notion of stable convergence and its variants)},
  volume~57.
\newblock Unione Matematica Italiana, Monograf s.r.l., Bologna, Italy., 2016.
\newblock Book written in Italian.

\bibitem{cri-dai-min}
I.~Crimaldi, P.~Dai~Pra, and I.~G. Minelli.
\newblock Fluctuation theorems for synchronization of interacting {P}\'olya's
  urns.
\newblock {\em Stochastic Process. Appl.}, 126(3):930--947, 2016.

\bibitem{cri-let-pra-2007}
I.~Crimaldi, G.~Letta, and L.~Pratelli.
\newblock A strong form of stable convergence.
\newblock In {\em S\'eminaire de {P}robabilit\'es {XL}}, volume 1899 of {\em
  Lecture Notes in Math.}, pages 203--225. Springer, Berlin, 2007.

\bibitem{cri-dai-lou-min}
I.~Crimaldi, P.~D. Pra, P.-Y. Louis, and I.~G. Minelli.
\newblock Synchronization and functional central limit theorems for interacting
  reinforced random walks.
\newblock {\em Stochastic Processes and their Applications}, 129(1):70 -- 101,
  2019.

\bibitem{cri-pra}
I.~Crimaldi and L.~Pratelli.
\newblock Convergence results for multivariate martingales.
\newblock {\em Stochastic Process. Appl.}, 115(4):571--577, 2005.

\bibitem{dai-lou-min}
P.~Dai~Pra, P.-Y. Louis, and I.~G. Minelli.
\newblock Synchronization via interacting reinforcement.
\newblock {\em J. Appl. Probab.}, 51(2):556--568, 2014.

\bibitem{egg-pol}
F.~Eggenberger and G.~P\'{o}lya.
\newblock \"{U}ber die statistik verketteter vorg\"{a}nge.
\newblock {\em Z. Angewandte Math. Mech.}, 3:279--289, 1923.

\bibitem{fortini}
S.~Fortini, S.~Petrone, and P.~Sporysheva.
\newblock {On a notion of partially conditionally identically distributed
  sequences}.
\newblock {\em Stoch. Process. Appl.}, 128:819--846, 2018.

\bibitem{ghi-pag14}
A.~Ghiglietti and A.~M. Paganoni.
\newblock Statistical properties of two-color randomly reinforced urn design
  targeting fixed allocations.
\newblock {\em Electron. J. Stat.}, 8(1):708--737, 2014.

\bibitem{ghi-vid-ros}
A.~Ghiglietti, A.~N. Vidyashankar, and W.~F. Rosenberger.
\newblock Central limit theorem for an adaptive randomly reinforced urn model.
\newblock {\em Ann. Appl. Probab.}, 27(5):2956--3003, 2017.

\bibitem{hall-1980}
P.~Hall and C.~C. Heyde.
\newblock {\em Martingale limit theory and its application}.
\newblock Academic Press, Inc. [Harcourt Brace Jovanovich, Publishers], New
  York-London, 1980.
\newblock Probability and Mathematical Statistics.

\bibitem{ieee-paper}
M.~Hayhoe, F.~Alajaji, and B.~Gharesifard.
\newblock A polya urn-based model for epidemics on networks.
\newblock In {\em 2017 American Control Conference (ACC)}, pages 358--363,
  2017.

\bibitem{laru-page}
S.~Laruelle and G.~Pag\`es.
\newblock Randomized urn models revisited using stochastic approximation.
\newblock {\em Ann. Appl. Prob.}, 23:1409--1436, 2013.

\bibitem{lima}
Y.~Lima.
\newblock Graph-based {P}\'olya's urn: completion of the linear case.
\newblock {\em Stoch. Dyn.}, 16(2):1660007, 13, 2016.

\bibitem{mah}
H.~M. Mahmoud.
\newblock {\em {P}\'olya urn models}.
\newblock Texts in Statistical Science Series. CRC Press, Boca Raton, FL, 2009.

\bibitem{mok-pel}
A.~Mokkadem and M.~Pelletier.
\newblock Convergence rate and averaging of nonlinear two-time-scale stochastic
  approximation algorithms.
\newblock {\em Ann. Appl. Probab.}, 16(3):1671--1702, 2006.

\bibitem{pag-sec}
A.~M. Paganoni and P.~Secchi.
\newblock Interacting reinforced-urn systems.
\newblock {\em Adv. in Appl. Probab.}, 36(3):791--804, 2004.

\bibitem{pem}
R.~Pemantle.
\newblock A survey of random processes with reinforcement.
\newblock {\em Probab. Surv.}, 4:1--79, 2007.

\bibitem{z}
L.-X. Zhang.
\newblock A {G}aussian process approximation for two-color randomly reinforced
  urns.
\newblock {\em Electron. J. Probab.}, 19(86):1--19, 2014.

\end{thebibliography}

\end{document}